\documentclass[aos]{imsart}

\RequirePackage{amsthm,amsmath,amsfonts,amssymb}
\RequirePackage[numbers]{natbib}
\usepackage{bm}
\usepackage{enumerate}
\usepackage{graphicx}
\usepackage{mathrsfs}
\usepackage{listings}
\usepackage{multirow}
\usepackage{caption}
\usepackage{subcaption}
\usepackage{bbm}
\usepackage{makecell}
\usepackage{url}
\usepackage{float}
\usepackage{color,soul} 
\usepackage{xcolor} 
\usepackage[linesnumbered,ruled,vlined]{algorithm2e}
\usepackage{mathtools}
\usepackage{accents}
\usepackage[shortlabels]{enumitem}
\usepackage{chngcntr}
\RequirePackage[colorlinks,citecolor=blue,urlcolor=blue]{hyperref}

\startlocaldefs
\theoremstyle{plain}
\newtheorem{theorem}{Theorem}
\newtheorem{lemma}{Lemma}
\newtheorem{corollary}{Corollary}

\newtheorem{condition}{Condition}
\theoremstyle{definition}

\newtheorem{remark}{Remark}

\newcommand{\eps}{\epsilon}
\newcommand{\pr}{\mathrm{Pr}}
\newcommand{\var}{\mathrm{Var}}
\newcommand{\cov}{\mathrm{Cov}}
\newcommand{\E}{\mathbb{E}}

\newcommand{\mr}{\mathrm}

\newcommand{\C}{\mathcal{C}}
\newcommand{\D}{\mathcal{D}}
\newcommand{\N}{\mathcal{N}}
\newcommand{\X}{\mathcal{X}}
\newcommand{\Y}{\mathcal{Y}}

\newcommand{\inner}[1]{\langle #1 \rangle}
\newcommand{\expinner}[1]{e^{-\iota\langle #1, \xi\rangle}}
\newcommand{\hatK}[1]{\hat{K}^{\tau,s}_{#1}}

\renewcommand{\appendix}{
  \setcounter{section}{0} 
  \renewcommand{\thesection}{\Alph{section}} 
  \renewcommand{\section}{
    \@startsection {section}{1}{\z@}
      {-3.5ex \@plus -1ex \@minus -.2ex}
      {2.3ex \@plus.2ex}
      {\normalfont\Large\bfseries}}}
\makeatother
\endlocaldefs

\begin{document}

\begin{frontmatter}
\title{Vecchia Gaussian Processes: on probabilistic and statistical properties}
\runtitle{Vecchia Gaussian processes}

\begin{aug}
\author[A]{\fnms{Botond}~\snm{Szabo}\ead[label=e1]{botond.szabo@unibocconi.it}}
\and
\author[B]{\fnms{Yichen}~\snm{Zhu}\ead[label=e2]{yczhu@hku.hk}\orcid{0000-0003-4050-5874}}
\address[A]{Bocconi Institute for Data Science and Analytics and Department of Decision Sciences,
Bocconi University \printead[presep={ ,\ }]{e1}}
\address[B]{Department of Statistics and Actuarial Science,
the University of Hong Kong \printead[presep={,\ }]{e2}}
\runauthor{B. Szabo and Y. Zhu}
\end{aug}

\begin{abstract}
Gaussian Processes (GPs) are widely used to model dependencies in spatial statistics and machine learning. However, exact inference is computationally intractable for GP regression, with a time complexity of $O(n^3)$. The Vecchia approximation scales up computation by introducing sparsity into the spatial dependency structure, represented by a directed acyclic graph (DAG). Despite its practical popularity, this approach lacks rigorous theoretical foundations, and the choice of DAG structure remains an open problem.

In this paper, we systematically study the Vecchia approximation of the popular, isotropic Mat\'{e}rn GP  as standalone stochastic process and uncover key probabilistic and statistical properties. We propose selecting parent sets as norming sets with fixed cardinality in the Vecchia approximation. On the probabilistic side, we show that the conditional distributions of Matérn GPs, as well as their Vecchia approximations, can be characterized by polynomial interpolations. This enables us to establish several results on small ball probabilities and the Reproducing Kernel Hilbert Spaces (RKHSs) of Vecchia GPs. Building on these probabilistic results, we prove that in the nonparametric regression model, the corresponding posterior contracts around the truth at the optimal minimax rate, both under oracle rescaling and hierarchical tuning of the prior.

We illustrate the theoretical findings through numerical experiments on synthetic datasets. Our core algorithms are implemented in C$++$ with an R interface.

\end{abstract}



\end{frontmatter}


\section{Introduction}
\subsection{Background and related work}\label{sec:background}
Gaussian processes have seen a wide range of applications among others in spatial statistics (e.g., \cite{gelfand2016spatial, cressie2011statistics, gelfand2010handbook}), machine learning (e.g., \cite{williams2006gaussian, lee2017deep, liu2020gaussian}) and epidemiology \cite{bhatt2015effect} due to their flexibility in characterizing dependency structures and their convenience for uncertainty quantification. However, exact statistical inference in the Gaussian Process (GP) regression model -- including maximum likelihood estimation of GP parameters and posterior inference based on GP priors -- suffers from a computational complexity of $O(n^3)$, owing to the need to compute the determinant and inverse of the $n\times n$ covariance matrix.

A substantial body of work has focused on scalable inference methods for Gaussian Processes, including but not limited to inducing variable or reduced rank approximations \citep{quinonero2005unifying,cressie2008fixed,banerjee2008gaussian,Titsias2009InducingVariablesVB}, covariance tapering \citep{furrer2006covariance,kaufman2008covariance,stein2013statistical}, compositional likelihoods \citep{bai2012joint,bevilacqua2015comparing}, distributed methods \citep{deisenroth2015distributed,daxberger2017distributed, szabo2019asymptotic} and Vecchia approximations. 
Vecchia approximations of Gaussian processes (Vecchia GPs), named after Aldo V. Vecchia \cite{vecchia1988estimation}, are constructed as follows. First, consider decomposing the density of a random vector $Z_{\mathcal{X}_n}=(Z_{X_1},...,Z_{X_n})\in\mathbb{R}^n$ into a product of univariate conditional densities, 
$$p(Z_{\X_n}) = p(Z_{X_1})\prod_{i=2}^n p(Z_{X_i}|Z_{X_{<i}}),$$
where $Z_{X_{<i}}=\{Z_{X_1},...,Z_{X_{i-1}}\}$.
The key idea behind the Vecchia approximation is to replace each conditioning set  $Z_{X_{<i}}$ with a subset $Z_{\mr{pa}(X_i)}\subseteq Z_{X_{<i}}$, often referred to as the \textit{parent set}, thereby introducing sparsity into the implied graphical model.
More formally, the joint density is approximated as
$$p(Z_{\X_n}) \approx p(Z_{X_1})\prod_{i=2}^n p(Z_{X_i}|Z_{\mr{pa}(X_i)}).$$
Applying the above approximation to the finite dimensional marginals of a GP results, under mild assumptions, in a new GP, called the \textit{Vecchia GP}.\footnote{For a proof, we refer to \cite{datta2016hierarchical} or \cite{zhang2024fixed}. We provide a more formal description of the process in Section \ref{sec:vecchia GP}.} The approximated GP is referred to as the \textit{mother GP.} The underlying conditional dependence structure encoded in the approximation -- namely the set of parents for each covariate -- is summarized by the DAG.

There is a substantial literature on Vecchia approximations of GPs (e.g., \cite{datta2016hierarchical, katzfuss2020vecchia, katzfuss2021general, peruzzi2022highly, cao2022scalable}). Their popularity can largely be attributed to two practical advantages. First, when all parent sets have bounded cardinality, Vecchia approximations allow evaluation of the joint density in $O(n)$ computational time. This makes Vecchia GPs one of the most scalable methods for GP approximation.
Second, the Vecchia approximation of a GP is itself a valid GP. Consequently, inference based on Vecchia GPs is statistically meaningful regardless of how well the approximation matches the mother GP.

Despite their success in applications and widespread use in the literature, major methodological and theoretical challenges remain. Methodologically, the optimal choice of the DAG structure is unclear, even though it plays a crucial role in the performance of Vecchia GPs. Indeed, several proposed strategies contradict each other. For example, many papers advocate selecting nearby neighbors as parent sets  \citep{datta2016hierarchical,zhu2024radial}, while others propose choosing more distant locations \citep{banerjee2015hierarchical} or even random ordering \citep{guinness2018permutation}. Moreover, while it is intuitively clear that larger parent sets yield better approximations at the expense of slower computation, there is no consensus on an appropriate parent-set size. 
The qualitative results available in the literature are only of polylogarithmic order \citep{schafer2021compression,kang2024asymptotic,zhu2024radial}, and it remains unclear whether further increasing the cardinality provides substantial statistical benefits. We summarize the above challenges in the following two research questions:
\begin{itemize}
\item \textbf{Q1} What is the minimal cardinality of the parent sets in the Vecchia GP needed to achieve optimal statistical inference?\\
 \vspace{-4mm}
\item \textbf{Q2} How should the parent sets be chosen to attain optimal inference?
\end{itemize}

We address these questions from a frequentist perspective. We show that suitably tuned Vecchia GPs can achieve optimal statistical inference in a minimax sense. In particular, 
we investigate the posterior contraction rates of Vecchia GPs around the true functional parameter of interest and demonstrate that  minimax-optimal rates can be achieved for appropriately chosen DAG structures and optimally tuned mother GPs. 

Despite the great popularity of Vecchia GPs in practice,  their theoretical understanding has so far been rather limited. We are aware of only a few recent theoretical contributions. In the univariate case and under certain assumptions on the variance, \cite{zhang2024fixed} derives asymptotic normality of  Vecchia GPs. Under suitable conditions on the DAG structure, covariance function, and dataset, \cite{schafer2021compression} shows that the precision matrix of a Vecchia GP closely approximates that of the corresponding mother GP. Additionally, \cite{zhu2024radial} demonstrates that Vecchia GPs provide accurate approximations to the full GPs in Wasserstein distance. Furthermore, \cite{kang2024asymptotic} shows that the maximum likelihood estimators of the Vecchia GP parameters have the same convergence properties and asymptotic normality properties as those of the mother GP.

Although all the above papers offer solid contributions to the literature of Vecchia GPs, none of them take a frequentist perspective, addressing the contraction rate of the posterior around the truth in nonparametric models. In fact, all previous works tackle the theory of Vecchia GPs from the perspective of ``approximating the mother GP''. They either directly quantify the approximation error between Vecchia GPs and the corresponding mother GPs or study the statistical properties of Vecchia GPs based on error bounds between the two. However, such an indirect strategy might not lead to optimal results. Vecchia GPs are well-defined, practically relevant Gaussian processes in their own right. Therefore, their properties -- both empirical and theoretical -- shall be studied as standalone problems that do not require good approximations to their mother GPs as prerequisites. After all, the main goal is to provide a scalable process with statistical properties comparable to those of the original GP. Restricting ourselves to an overly accurate approximation of the mother GP might prevent us from fully exploring the computational advantages offered by the Vecchia approach.

To investigate the statistical aspects of Vecchia GPs, one has to  first understand their probabilistic properties, which so far have been barely studied.  More concretely, GPs are typically defined through their marginal distributions and the standard probabilistic techniques are tailored to this perspective, see, for example \cite{kuelbs1994gaussian,li:linde:1999,aurzada2009small}. However, Vecchia GPs are defined through conditional distributions encoded in a DAG structure. For such constructions, to the best of our knowledge, no general tools have been developed for computing small deviation bounds or other probabilistic characteristics of the GP. Hence, one of the main challenges of our work was to derive new mathematical techniques for studying the probabilistic properties of Vecchia approximated GPs. To highlight that the mother and Vecchia GPs can have substantially different properties, we note that in the case of the popular Mat\'{e}rn process, the mother GP is stationary, whereas its Vecchia approximation is not. In fact, little is known about the covariance kernel of the Vecchia GP or the corresponding Reproducing Kernel Hilbert Space (RKHS). Therefore, techniques used for isotropic covariance functions --such as the Mat\'{e}rn kernel -- that rely on spectral densities no longer apply in the case of Vecchia approximations.

\subsection{Our contributions}
In this section, we collect our contributions on the probabilistic and statistical properties of Vecchia Gaussian Processes. We focus on isotropic Mat\'ern processes, while extensions beyond this class of GPs are briefly discussed in Section \ref{sec:discussion}.
These results build heavily on one another and extend our understanding of this approximation technique in several directions.

\subsubsection{Probabilistic properties}
Our contributions to the probabilistic aspects are twofold. First, we systematically study the conditional distributions that define Vecchia GPs. We show that the conditional expectations of Vecchia Mat\'ern GPs can be characterized asymptotically via local polynomial interpolations on the parent sets. Furthermore, we prove that the conditional variances coincide with the approximation errors of these local polynomial approximations. Moreover, under suitable DAG structures, the degree of these local polynomials is exactly $\underline{\alpha}$, the largest integer strictly smaller than the regularity of the Mat\'{e}rn process. While the limits of Gaussian interpolation have been studied in the literature on radial basis functions \citep{driscoll2002interpolation,larsson2005theoretical,song2012multivariate,lee2015study}, to the best of our knowledge, we are the first to prove the stronger version with uniform convergence rate and to make this connection to Vecchia GPs. In contrast, the characterization of the variance through the approximation error is an entirely novel contribution.

Our second contribution is the derivation of various probabilistic properties of Vecchia GPs based on their conditional probabilities, rather than following the standard route through marginal probabilities. One important quantity is the small deviation bound (i.e. small ball probability) of the GP, which is typically derived using its relation to the entropy of the unit ball of the associated RKHS \citep{kuelbs1994gaussian, li:linde:1999, aurzada2009small, ghosal2017fundamentals}. However, for Vecchia GPs, the properties of the corresponding RKHS are not well characterized. In our work, we provide a new set of tools tailored to Vecchia GPs, which are also potentially useful for other types of processes.

\subsubsection{Statistical properties}
Building upon the probabilistic properties discussed above, we prove -- within the context of the nonparametric regression model -- that appropriately tuned Vecchia GPs attain minimax posterior contraction rates for H\"{o}lder smooth functions. This demonstrates that Vecchia GPs can substantially reduce the computational costs of their mother GPs without sacrificing statistical accuracy. The optimal choice of the tuning/scaling parameter depends on the smoothness of the underlying function of interest, which is typically not available. Therefore, we also consider a data-driven tuning strategy by endowing the scaling parameter with another layer of prior, forming a two-level, hierarchical prior distribution. We show that the posteriors arising from such hierarchical Vecchia GPs, under mild assumptions, concentrate their masses around the true function at the minimax rate. Alternatively, one may plug-in the marginal maximum likelihood estimators of the hyperparameters into the Vecchia GP. The two approaches possess similar theoretical properties; see, for instance, \cite{szabo:etal:2013,rousseau2017asymptotic}. In the present paper, we focus on the former, hierarchical approach and leave the analysis of the latter for future work.

Our research was partially motivated by the methodological questions Q1 and Q2. We address these questions from a frequentist perspective, focusing on optimal recovery of the underlying true function in the nonparametric regression model. We show, that well-chosen parent sets with cardinality equal to ${\underline{\alpha}+d \choose \underline{\alpha}}$ are sufficient to achieve minimax contraction rate with the Vecchia approximation. More concretely,  we propose selecting a \textit{norming set} as parents with a fixed, universal norming constant, see Section \ref{sec:poly2DAG} for the definition and background. For (approximately) grid-structured data, we provide an explicit construction of such sets, while for more general datasets we propose algorithms that can potentially identify such norming sets when they exist. Numerical studies in Section \ref{sec:num} demonstrate that our proposal achieves high estimation accuracy while maintaining fast computation.

Our perspective and consequently our answers to questions Q1 and Q2 differ from much of the existing literature on Vecchia GPs. We believe that these differences reveal previously hidden features and characteristics, thereby contributing to a deeper understanding of this popular, yet so far theoretically underdeveloped method. In Section \ref{sec:discussion}, we provide a detailed discussion of our results and their implications. Here we focus on one of the most important methodological aspect: the cardinality of the parent sets. The results in \cite{schafer2021compression}, \cite{zhu2024radial} and \cite{kang2024asymptotic} all require the parent sets to grow at a poly-logarithmic rate in $n$, whereas we require only fixed sized parent sets. There is, however, no contradiction here. Previous work typically considers 
Vecchia GP approximations to the mother GPs and the poly-logarithmic growth ensures vanishing approximation error. In contrast, we investigate directly the properties of the posterior induced by the Vecchia GP as standalone prior, and we derive contraction rates in a nonparametric framework. By bypassing the approximation step and taking this more direct approach, one can typically obtain weaker conditions.

\subsection{Notation}\label{sec:nota}
Next we introduce some notation used throughout the paper. For easier readability, we organize them into subsections.

\subsubsection{Gaussian processes}
Let $\X\subset\mathbb{R}^d$ be the domain of the centered mother GP $(Z_x, x\in \X)$ with covariance kernel $K(\cdot,\cdot):\X\times \X \to \mathbb{R}$.
We focus on the Mat\'{e}rn covariance function (Section 4.2.1 of \cite{seeger2004gaussian}) with regularity hyperparameter $\alpha>0$:
\begin{equation}\label{eq:MaternCov o}
K(x_1,x_2) = \frac{2^{1-\alpha}}{\Gamma(\alpha)} \|x_1-x_2\|_2^{\alpha} \mathcal{K}_{\alpha}\left(\|x_1-x_2\|_2\right),
\end{equation}
where $\mathcal{K}_\alpha$ denotes the modified Bessel function of the second kind with parameter $\alpha$. The Mat\'ern kernel also admits a simple Fourier characterization:
\begin{equation}\label{eq:MaternCov}
    K(x_1, x_2) \propto \int_{\mathbb{R}^d} e^{-\iota\langle x_1-x_2,\xi \rangle} (1+\|\xi\|_2^2)^{-\alpha-d/2} d\xi,
\end{equation} 
where the proportionality constant depends only on $\alpha$ and $d$, see equation (11.4) of \cite{ghosal2017fundamentals}. For $s,\tau>0$, the rescaled GP $(Z_x^{\tau,s},x\in \X)$ is defined by
\begin{equation*}
    Z^{\tau,s}_x = s Z_{\tau x}, \;\;\forall x\in \tau^{-1} \X.
\end{equation*}
We will refer to $\tau$ and $s$ as the time and space (re)scaling hyperparameters, respectively. 

Let $(\hat{Z}_x^{\tau,s},x\in \X)$ denote the Vecchia approximation of the mother GP $(Z^{\tau,s}_x,x\in\X)$, defined formally in Section \ref{sec:vecchia GP}. When $\tau=s=1$, we abbreviate this as $(\hat{Z}_x,x\in \X)$ and let $\hat{K}(\cdot, \cdot):\X\times\X\to\mathbb{R}$ denote its covariance kernel, which is typically different from the mother covariance kernel $K(\cdot,\cdot)$. For arbitrary scaling hyperparameters $\tau,s>0$, let $K^{\tau,s}(\cdot,\cdot)$ and $\hat{K}^{\tau,s}(\cdot,\cdot)$ denote the covariance kernels of the corresponding rescaled mother and Vecchia processes, respectively. For notational convenience, we abbreviate $(Z^{\tau,s}_x,x\in \X)$ as $Z^{\tau,s}$, and similarly write $\hat{Z}^{\tau,s}$, $Z$ and $\hat{Z}$. Let $\mathbb{H}$ and $\mathbb{H}^{\tau,s}$ denote the RKHSs associated with the mother GPs $Z$ and $Z^{\tau,s}$, respectively. Likewise, let $\mathbb{H}_n$ and $\mathbb{H}_n^{\tau,s}$ denote the RKHSs for the $n$-dimensional marginals on $\mathcal{X}_n$ of the processes $\hat{Z}$ and $\hat{Z}^{\tau,s}$.

For an arbitrary covariance function $K(\cdot,\cdot):\mathbb{R}^d\times\mathbb{R}^d \to \mathbb{R}$, and multi-indices $k_1,k_2\in\mathbb{N}^d$, we denote by $K^{(k_1,k_2)}(x_1,x_2) = \frac{\partial^{k_1}}{\partial x_1^{k_1}} \frac{\partial^{k_2}}{\partial x_2^{k_2}} K(x_1,x_2)$ the mixed partial derivatives of order $k_1$ and $k_2$ with respect to the first and second arguments, respectively.
For any finite sets $A,B\subset\X$, we denote by $K_{A,B}$ the covariance matrix between $A$ and $B$ under the covariance kernel $K(\cdot,\cdot)$. If $B=\{x\}$ is a singleton, we slightly abuse our notation and write $K_{A,x} = K_{A,\{x\}}$.

\subsubsection{Ordering and indices}
For a $d$-dimensional vector $x\in\mathbb{R}^d$, we denote its coordinates by
\begin{equation}\label{eq:coord}
x = (x[1],x[2],\cdots, x[d])^T.
\end{equation}
Similarly, $M[i,j]$ denotes the entry in the $i$th row and $j$th column of a matrix $M$. Furthermore, for multi-index $k=(k[1],k[2],\cdots,k[d])^T\in\mathbb{N}^d$, we use the conventions $|k|=\sum_{j=1}^d k[j]$ and $k! = \prod_{j=1}^d k[j]!$.

For vectors in $\mathbb{R}^d$, we define the lexicographical ordering ``$\prec$'' as follows: for $x=(x[1],x[2],\cdots,x[d])^T$ and $y=(y[1],y[2],\cdots,y[d])^T$, we say that $x\prec y$ if either $\sum_i x[i] < \sum_i y[i]$ or $\sum_i x[i] = \sum_i y[i]$ and there exists $ d'\le d$ such that $x[d']<y[d']$ and $x[i]=y[i],\;\forall i\le d'-1$.
Let $k_{(i)}\in\mathbb{N}^d$ denote the $i$th multi-index in $\mathbb{N}^d$ with respect to this lexicographical ordering. Note that the first index is $k_{(1)}=(0,0,\cdots,0)^T$.

Finally, let $f:\mathbb{R}^d \to\mathbb{R}$ be an $l$th-order differentiable function. For any multi-index $k\in\mathbb{N}^d$, the corresponding partial derivative of $f$  at $x\in \mathbb{R}^d$ is denoted by
$$f^{(k)}(x) = \prod_{i=1}^d \frac{\partial^{k[i]}}{\partial x[i]^{k[i]}} f(x). $$

\subsubsection{Sets, functions and functionals}\label{sec:notation set}
For a finite set $A$, we denote its cardinality by $|A|$. For $t>0$, we write
$t A=\{t x: x\in A\}$. For a set $A=\{x_1, x_2, \cdots, x_{|A|}\}\subset \X$ and a function $f:\X \to \mathbb{R}$, the vectorized version of $f$ over $A$ is
$ \bm{f}(A) = (f(x_1), f(x_2), \cdots, f(x_{|A|}))^T\in\mathbb{R}^{|A|}.$

 For a set $\Omega\subset\mathbb{R}^d$ and a function $f:\Omega\to\mathbb{R}$, let $\|f\|_{C^\alpha}$ denote the H\"{o}lder-$\alpha$ norm of $f$ and let $C^\alpha(\Omega)_1$ be the unit ball in this norm. Similarly, let $\|f\|_{W^{\alpha}}$ be the Sobolev-$\alpha$ norm (equivalent to the $\|\cdot\|_{2,2,\alpha}$ norm in Definition C.6 of \cite{ghosal2017fundamentals}) and denote by $W^{\alpha}(\Omega)_1$ the corresponding unit ball. When $\Omega$ coincides with the domain $\X$ of the GP, we abbreviate these sets by $C^\alpha_1$ and $W^{\alpha}_1$. The $L_p(\Omega)$ norm of the function $f$ is denoted by $\|f\|_p$.

For a finite set $\X_n\subset\mathbb{R}^d$, the supremum- and $L_2$-norm of a function $f:\mathbb{R}^d\to\mathbb{R}$ restricted to $\X_n$ are defined as $\|f\|_{\infty,n} = \sup_{x\in\X_n} |f(x)|$ and $\|f\|_{2,n}=[[|\X_n|^{-1}\sum_{x\in\X_n} f(x)^2]^{1/2}$, respectively. For sets $\Omega_1, \Omega_2$ , $p\in(0,+\infty]$ and an operator $M: L_p(\Omega_1) \to L_p(\Omega_2)$, the $L_p$ operator norm  is $\|M\|_p = \sup_{f\in L_p(\Omega_1)} \|M(f)\|_p / \|f\|_p$.

\subsubsection{Other notation}
We use $\iota$ to denote imaginary unit. For $\alpha\in\mathbb{R}$, let $\underline{\alpha}$ denote the largest integer strictly smaller than $\alpha$. For two functions $f_1(n)$ and $f_2(n)$ of the sample size $n$,  the notation $f_1(n)\lesssim f_2(n)$ means that there exists a universal constant $c>0$, such that $f_1(n)\le cf_2(n)$ for all $n$. 
Let $\mathbbm{1}_m$ be the $m$-dimensional column vector whose entries are all equal to one. Finally, we denote by  $Z_1\overset{d}{=}Z_2$ that the random variables $Z_1,Z_2$ have the same distribution.

\subsection{Organization of the paper}
The rest of the paper is organized as follows. Section \ref{sec:vecchia gp} introduces Vecchia approximations of GPs.  Section \ref{sec:layered norming dag} defines the layered norming DAG, providing a principled, theory driven choice of the underlying graphical structures. Section \ref{sec:prob} establishes key probabilistic properties of Mat\'{e}rn processes and their Vecchia counterparts. Building on these results, Section \ref{sec:BNP} derives minimax contraction rates for appropriately rescaled Vecchia GPs, using both oracle tuning and hierarchical Bayesian approaches.  Finally, Section \ref{sec:discussion} discusses the limitations and broader implications of our results. Due to space constraints, we defer the numerical analysis, some background material on polynomials and the concept of norming sets, as well as all proofs, to the Supplementary Material.

\section{Vecchia approximations of Gaussian processes}\label{sec:vecchia gp}
We first describe the Gaussian Process regression model, 
and then introduce the Vecchia approximation for GPs. We consider Vecchia GPs both with fixed hyperparameters and with an additional layer of prior, resulting in a hierarchical prior distribution.

\subsection{Gaussian process regression}\label{sec:GP regression}
Assume that we observe $Y_1,...,Y_n$ satisfying
\begin{equation}\label{eq:regression}
Y_i = f(X_i) + \varepsilon_i, \quad \varepsilon_i \overset{\text{i.i.d.}}{\sim} N(0,\sigma^2),\quad i=1,...,n,
\end{equation}
where the covariates $(X_i)_{1\leq i\leq n}$ are deterministic, and belong to a compact subset $\mathcal{X}\subset\mathbb{R}^d$,  $f\in L_2(\mathcal{X})$ is the unknown functional parameter of interest, and $\sigma^2$ is the unknown noise variance. Let us denote by $\D_n=\{(X_i,Y_i)_{1\leq i\leq n}\}$ the observed data. We adopt a Bayesian Gaussian process regression framework by endowing $f$ with a GP prior, i.e.
\begin{equation}\label{eq:prior}
f \sim (Z^{\tau,s}_x, x\in\mathcal{X}),
\end{equation}
where $Z^{\tau,s}$ is a GP  with rescaling parameters $(\tau,s)$. We also place a prior on the noise variance $\sigma^2$, denoted as 
\begin{equation}\label{eq:prior nugg}
\sigma^2 \sim p(\sigma^2).
\end{equation}
The posterior distribution of the regression function given the data, $f|\D_n$, provides a probabilistic solution to the regression problem. If the hyperparameters $\tau, s$ and other GP tuning parameters are fixed, conjugacy ensures that the posterior is also a Gaussian process. In practice, however, choosing these hyperparameters directly is often difficult. Therefore, it is common to place another layer of prior on them, letting the data guide their selection \footnote{There may exist other hyperparameters for the process $Z^{\tau,s}_x$ beyond $(\tau,s)$, such as the regularity parameter $\alpha$ in the Mat\'{e}rn kernel (\ref{eq:MaternCov}). However, fitting these hyperparameters is often computationally more challenging.}, i.e. 
\begin{equation}\label{eq:theta}
(\tau,s) \sim p(\tau,s).
\end{equation}
Equations (\ref{eq:regression}) - (\ref{eq:theta}) jointly define the Gaussian process regression model. Due to the presence of the hyper-prior (\ref{eq:theta}), the conjugacy is typically lost and the posterior is computed via Markov Chain Monte Carlo.

Even in the fixed hyperparameter setting, where the GP posterior has an explicit form, the computation involves the inversion of an $n\times n$ matrix, which has computational complexity of order $O(n^3)$. This complexity is intractable for many practical applications where $n$ can be as large as millions. Therefore, various approximation methods have been proposed in the literature to speed up the computations. In the next section, we describe one of the most popular such approaches: the Vecchia GP.

\subsection{Vecchia Gaussian processes}\label{sec:vecchia GP}
Vecchia approximations of Gaussian processes, or Vecchia GPs, are a class of methods that allow scalable computation for Gaussian processes. Specifically, a Vecchia GP $\hat{Z}^{\tau,s}$ is defined by two components: the mother Gaussian process $Z^{\tau,s}$ on $\mathcal{X}$ and a directed acyclic graph $G=(\mathcal{X},E)$. In most cases, $\mathcal{X}$ is an uncountable set and the directed acyclic graph $G$ can be decomposed into two subgraphs: a graph defined on a finite reference set, and a bipartite graph between the reference set and the remaining nodes. For simplicity, in this paper we take the observed covariates $\X_n=(X_1,...,X_n)$ as the reference set. The graph on $\X_n$ characterizes the dependence structure within  reference set. 
The vertices of the bipartite graph on $\mathcal{X}$ consist of two parts: $\X_n$ and $\mathcal{X}\backslash \X_n$. All directed edges in this bipartite graph originate from vertices in $\X_n$ and terminate at vertices in $\mathcal{X}\backslash \X_n$. In other words, the nodes of $\mathcal{X}\backslash \X_n$ are conditional independent given $\X_n$\footnote{It is worth noting that the nodes of $\X\backslash\X_n$ remain marginally dependent. When the sample size $n$ is large, this marginal dependency can be sufficiently strong to capture spatial patterns. }.

 The marginal densities of the Vecchia GP $\hat{Z}^{\tau,s}$ are defined in two steps. First, for the root vertex $X\in \X_n$, we have
\begin{equation}\label{eq:vecchia 1}
\hat{Z}^{\tau,s}_{X} \overset{d}{=}   Z^{\tau,s}_{X}.
\end{equation}
For any $X\in\mathcal{X}$, we say $X'$ is a parent of $X$ if there is a directed edge from $X'$ to $X$ in the DAG $G$. We denote the parent set of $X$, i.e., the collection of all parents of $X$, by $\mr{pa}(X)$. Then $\hat{Z}^{\tau,s}_X$ conditionally on the process value at its parent locations, follows the Gaussian distribution:
\begin{align}\label{eq:vecchia 2}
\hat{Z}^{\tau,s}_{X}|\hat{Z}^{\tau,s}_{\mr{pa}(X)} \sim N\bigg( & K_{\tau X,\tau \mr{pa}(X)} K_{\tau \mr{pa}(X),\tau \mr{pa}(X)}^{-1} \hat{Z}^{\tau,s}_{\mr{pa}(X)}, \nonumber\\
&\;\;s^2\big[K_{\tau X,\tau X} - K_{\tau X,\tau \mr{pa}(X)} K_{\tau \mr{pa}(X),\tau \mr{pa}(X)}^{-1} K_{\tau \mr{pa}(X),\tau X}\big] \bigg). 
\end{align}
Equation (\ref{eq:vecchia 2}) implies that $\hat{Z}^{\tau,s}_{X}|(\hat{Z}^{\tau,s}_{\mr{pa}(X)}=z) \overset{d}{=} Z^{\tau,s}_X|(Z^{\tau,s}_{\mr{pa}(X)}=z)$. Let $\hat{K}(\cdot,\cdot)$ be the covariance function of the process $\hat{Z}^{\tau,s}$. Then the relationship between $K$ and $\hat{K}$ can be summarized by the following two assertions, which will be used repeatedly in the sections that follow: 
\begin{align}\label{eq:vecchia exp}
 \hat{K}_{\tau X,\tau\mr{pa}(X)} &\hat{K}_{\tau\mr{pa}(X),\tau\mr{pa}(X)}^{-1} z  =
\mathbb{E}(\hat{Z}^{\tau,s}_X|\hat{Z}^{\tau,s}_{\mr{pa}(X)}=z) \nonumber \\
= & \mathbb{E}(Z^{\tau,s}_X|Z^{\tau,s}_{\mr{pa}(X)}=z)
= K_{\tau X,\tau\mr{pa}(X)} K_{\tau\mr{pa}(X),\tau\mr{pa}(X)}^{-1} z, \quad\forall z,  
\end{align}

\begin{align}\label{eq:vecchia var}
 s^2\big[\hat{K}_{\tau X,\tau X}& - \hat{K}_{\tau X,\tau\mr{pa}(X)} \hat{K}_{\tau\mr{pa}(X),\tau\mr{pa}(X)}^{-1} \hat{K}_{\tau\mr{pa}(X),\tau X} \big] \nonumber\\
= & \var\left[\hat{Z}^{\tau,s}_X - \mathbb{E} (\hat{Z}^{\tau,s}_X|\hat{Z}^{\tau,s}_{\mr{pa}(X)})\right]  =\var\left[Z^{\tau,s}_X -\mathbb{E}(Z^{\tau,s}_X|Z^{\tau,s}_{\mr{pa}(X)})\right]  \nonumber\\
= &  s^2\big[ K_{\tau X,\tau X} - K_{\tau X,\tau\mr{pa}(X)} K_{\tau\mr{pa}(X),\tau\mr{pa}(X)}^{-1} K_{\tau\mr{pa}(X),\tau X} \big].
\end{align}

Because $G$ is a directed acyclic graph, equations (\ref{eq:vecchia 1}) and (\ref{eq:vecchia 2}) consistently define the finite-dimensional marginals of $\hat{Z}^{\tau,s}$, and hence define a valid Gaussian process, see for example Section S1 of \cite{zhu2024radial} for a proof. Assuming $X_1$ is the root node, i.e., the node with no parent, we can write the joint density of the process $\hat{Z}^{\tau,s}$ on $\X_n$ as
\begin{equation}\label{eq:dag decomp}
p(\hat{Z}^{\tau,s}_{\X_n}) = p(\hat{Z}^{\tau,s}_{X_1}) \prod_{i=2}^n p(\hat{Z}^{\tau,s}_{X_i} |    \hat{Z}^{\tau,s}_{\mr{pa}(X_i)}).
\end{equation}
Intuitively, the Vecchia approximated process $\hat{Z}^{\tau,s}$ is obtained by removing certain conditional dependencies in the mother process $Z^{\tau,s}$, where the DAG $G$ specifies which conditional dependencies (or equivalently, which directed edges) are retained. As a consequence, the precision matrices  of the marginal distributions of $\hat{Z}^{\tau,s}$  and their Cholesky decompositions are sparse. This sparsity enables Vecchia GPs scalable computation on large data sets. Specifically, if the cardinality of the parent sets $\mr{pa}(X_i)$ for all covariates $X_i\in\mathcal{X}_n$  is bounded by  $m$, then the computational complexity of evaluating the joint density on $\X_n$ using equation (\ref{eq:dag decomp}) is $O(nm^3)$. This linear complexity in $n$ (when $m$ does not increasing with $n$) is significantly smaller than the $O(n^3)$ complexity of directly inverting the covariance matrix.

Endowing the regression function $f$ in equation (\ref{eq:regression}) with the Vecchia GP prior
\begin{equation}\label{eq:prior vecchia}
f\sim (\hat{Z}^{\tau,s}_x, x\in\mathcal{X}),
\end{equation}
the equations (\ref{eq:regression}), \eqref{eq:prior nugg}, (\ref{eq:prior vecchia}) and (\ref{eq:theta}) together specify the Vecchia GP regression model. Our paper focuses on analyzing the probabilistic and statistical properties of this model. It is clear that these properties crucially depend on the underlying DAG structure. Although
many different DAG structures have been proposed
based on various heuristics, we take a theory driven approach to selecting DAGs. Motivated by the flat limit of the conditional distribution of Mat\'{e}rn processes, we propose the \textit{layered norming DAG} in the following section.

\section{Layered norming DAGs}\label{sec:layered norming dag}
The foundation of our probabilistic and statistical results is the previously underexplored connection between local polynomial interpolation and conditional distributions of Mat\'ern processes. Due to space limitations, we present only the key concepts in this section, including norming sets and DAG structures. A brief introduction to this topic is provided in Section
\ref{sec:poly2DAG} of the supplement, where we define Vandermonde matrices and the unisolvency of polynomials, present useful lemmas, and discuss their connections. If the reader is not familiar with these concepts, we recommend consulting the supplement first. 

\subsection{Norming sets for polynomials}\label{sec:norming}
A \textit{norming set} is a key concept that connects the properties of polynomials with the spatial dispersion of certain finite point sets.
Let $\Omega$ be a compact subset of $\mathbb{R}^d$. Equipping the class of $l$th-order polynomials $\mathscr{P}_l(\Omega)$ with the supremum norm
$$\|P\| = \sup_{x\in\Omega} |P(x)|,\;\forall P \in \mathscr{P}_l(\Omega)$$
forms a normed vector space. 
A finite set $A=\{w_1,w_2,\cdots, w_m\}\subset\Omega$ is said to be a norming set for $\mathscr{P}_l(\Omega)$ with \textit{norming constant} $c_N>0$ if
\begin{equation*}
\sup_{x\in\Omega} |P(x)| \le c_N \sup_{x'\in A} |P(x')|, \;\forall P\in\mathscr{P}_l(\Omega).
\end{equation*}
The notion of norming sets was originally introduced by \cite{jetter1999error} for general normed vector spaces and their duals. In this paper, we focus on polynomials restricted to a compact set and point evaluation functionals on this set. Note that for a norming set $A$, the norm $\|\cdot\|$ and the weaker seminorm $\|\cdot\|_{A}$ 
$$\|P\|_A = \sup_{x\in A} |P(x)|,\;\forall P \in \mathscr{P}_l(\Omega),$$
are equivalent on $ \mathscr{P}_l(\Omega)$. Heuristically, this means that the behavior of a polynomial $P\in\mathscr{P}_l(\Omega)$ on a norming set $A\subset\Omega$ characterizes its behavior on the entire compact domain $\Omega$. The existence of such finite set follows from the finite, ${l+d \choose l}$ dimensionality of the vector space $\mathscr{P}_l(\Omega)$.

\subsection{The Layered and Norming conditions}
We begin by introducing the \textit{Layered} and \textit{Norming} conditions, which together define the layered norming DAGs, used in our proposal. For simplicity, we consider the domain $\mathcal{X}=[0,1]^d$.

We first describe the Layered condition.
\begin{condition}[Layered]\label{cond:layered}
We assume that the vertex set of the DAG can be partitioned as $\X_n=\cup_{j=0}^J \N_j\subset [0,1]^d$ satisfying
\begin{itemize}
\item $\mr{pa}(X_i)\subset\cup_{j'=0}^{j-1} \N_{j'}$, $\forall X_i\in\N_j$;
\item $\exists$ $\gamma>1$, $c_d>0$, such that $\|X_{i_1}-X_{i_2}\|_\infty\ge c_d\gamma^{-j}$, $\forall X_{i_1},X_{i_2}\in\cup_{j'=0}^j\N_{j'}$.
\end{itemize}
\end{condition}
Throughout the remainder of the paper, we assume that the indices of the vertices in $\X_n$ are ordered according to their corresponding layers. More precisely, denoting by $\eta(i)$ the layer of $X_i$, i.e. $X_i\in\N_{\eta(i)}$, we assume that for all $X_i, X_{i'}\in\X_n$ with $i<i'$, we have $\eta(i)\le \eta(i')$.

The Layered condition, as the name suggests, partitions the training set $\X_n$ into $J+1$ layers $\N_0, \N_1,\ldots, \N_J$, such that for any $X_i\in\mathcal{N}_j$, its parents can only belong to layers preceding $\mathcal{N}_j$. Moreover, the minimum distance between elements is governed by their respective layers. This layered structure facilitates a multiresolution framework in which the first few layers tend to be relatively spread out over the domain, while the latter layers gradually fill in the gaps. The layered condition is very mild and can be satisfied for an arbitrary training set  $\X_n$. 

Next we introduce the Norming condition for DAGs.
\begin{condition}[Norming]\label{cond:norming}
Consider a DAG satisfying Condition \ref{cond:layered} with parameter $\gamma>1$. Let $m = {l+d \choose l}$ and assume that there exists $i_0\in\mathbb{N}$, such that for all $i\ge i_0$, $|\mr{pa}(X_i)|=m$. Moreover, there exist constants $c_L, c_N>0$, such that for all  $X_i\in\N_j$,
with $ 0\le j\le J$ and $i\ge i_0$, there exists a $d$-dimensional cube $\mathcal{C} \supset \mr{pa}(X_i)\cup\{X_i\}$ with side length no greater than $c_L \gamma^{-j}$, such that $\mr{pa}(X_i)$ is a norming set on $\mathscr{P}_{l}(\mathcal{C})$ with norming constant $c_N$.
\end{condition}

The above Norming condition, built upon the Layered condition, requires that the parent set of each $X_i$ (for $i\ge i_0$) is a norming set on a neighborhood of $X_i$. Moreover, it assumes the existence of a universal norming constant $c_N$ for these norming sets. 

Unlike Condition \ref{cond:layered}, which can be satisfied for any $\X_n$ with a properly constructed DAG, the Norming condition, especially the requirement of a universal norming constant $c_N$, is not always straightforward to verify. In the univariate case it holds automatically, since any set of $m$ distinct points is a norming set for the polynomial space of order $m-1$, as computed in equation \eqref{eq:c_N d=1}. However, the situation is more complicated in the multivariate case, where not every set of size $m={l+d \choose l}$ forms a norming set for $\mathscr{P}_l([0,1]^d)$, let alone with a universal norming constant. Nevertheless, even in the higher dimensional cases there exist specific sets, such as grid design or Fekete points \cite{bos2018fekete}, that do possess uniformly bounded norming constants.

\subsection{Building layered norming DAGs}\label{sec:build dag}
In this section, we provide a constructive algorithm for finding DAG structures that satisfies Conditions \ref{cond:layered} and \ref{cond:norming}. As discussed in the previous subsection, the latter condition can be challenging to meet in high-dimensional spaces. For simplicity, we restrict our attention to grid data, and defer the discussion of general data set  $\X_n$ to Section \ref{sec:non grid} and Section \ref{sec:build dag general} of the Supplementary Material.

For notational convenience, suppose that $n=\tilde{n}^d$ for some $\tilde{n}\in\mathbb{N}$. Then we consider the grid data set $\X_n$ of the form
\begin{align}
\X_n = \big\{x=(x[1], x[2], \cdots, x[d])^T: x[j] \in \{0, 1/(\tilde{n}-1), 2/(\tilde{n}-1), \cdots, 1\}, \forall 1\le j\le d \big\}. \label{def:grid}
\end{align}
In other words, elements in $\X_n$ come from the $d$-dimensional tensor product of the set $\{0, 1/(\tilde{n}-1), 2/(\tilde{n}-1), \cdots, 1\}$. 

We first consider the case $\tilde{n} = 2^r+1, r\in\mathbb{N}$.
The construction of the layered norming DAGs consists of two steps: constructing the layers and assigning the parent sets. 
The layers $\N_j, 0\le j \le r$ are defined as follows. The base (root) layer $\N_0$ is 
$$\N_0 = \big\{x=(x[1], x[2], \cdots, x[d])^T: x[h] \in \{0, 1\}, \forall 1\le h\le d \big\}.$$
In other words, the layer $\N_0$ consists of the $2^d$ vertices of the $d$-dimensional unit cube.
For $1\le j\le r$, the layer $\N_j$ is defined as
\begin{equation}\label{eq:layer grid}
\N_j = \big\{x=(x[1], x[2], \cdots, x[d])^T: x[h] \in \{(2k+1)/2^j: 0\le k \le 2^{j-1}-1\}, \forall 1\le h\le d \big\}.
\end{equation}
In other words, the union $\cup_{j'=0}^j \N_{j'}$ forms a $(2^j+1)^d$ grid on the $d$-dimensional unit cube. 

The layers $\N_j$ defined above satisfy Condition \ref{cond:layered} with parameters $\gamma=2$ and $c_d=1$. It remains to specify the parent sets $\mr{pa}(X_i)$ for all $X_i$ such that Condition \ref{cond:norming} is also satisfied. 
Let $j_0$ be the smallest positive integer such that $2^{j_0-1}+1 \ge l+1$. In other words, $j_0$ is the first layer for which the preceding layers $\N_0,\N_1,\cdots, \N_{j_0-1}$ contain a sufficient number of points to form norming sets. The vertices in $\N_0$ have no parents. For all $X_i\in\N_j, 0<j<j_0$, the parent set is defined as
$\mr{pa}(X_i) = \cup_{j'=0}^{j-1} \N_{j'}$. For vertices in  $\N_j, j\ge j_0$ we employ the  \textit{corner sets} discussed in Section 3.2 of \cite{neidinger2019multivariate}, as described below. 
 
For each $X_i\in\N_j$, the corresponding parent set is a subset of 
$$\cup_{j'=0}^{j-1} \N_{j'} = \big\{x=(x[1], x[2], \cdots, x[d])^T: x[h] \in \{k/2^{j-1}: 0\le k \le 2^{j-1}\}, \forall 1\le h\le d \big\}.$$ 
For each coordinate $1\le h\le d$, we order the values $\{k/2^{j-1}: 0\le k \le 2^{j-1}\}$ by their distance from $X_i[h]$. In the case of  ties, the tied elements may be ordered arbitrarily. Let us denote the resulting ordered set by
$$B_h=(x_{s(1),h},x_{s(2),h},\cdots,x_{s(2^{j-1}+1),h}).$$
The \textit{corner} parent set of $X_i$ is then defined as 
\begin{equation}\label{eq:parent grid}
\mr{pa}(X_i)=\Big\{ (x_{s(i_1),1}, x_{s(i_2),2},\cdots, x_{s(i_d),d} ): \sum_{h=1}^d [i_h-1]\le l \Big\}.
\end{equation}
The term corner refers to the condition $\sum_{h=1}^d [i_h-1]\le l$, which ensures that the selected coordinate indices, ordered by increasing distance from each $X_i[h]$, lie near the ``corner'' of the $d$-dimensional index grid. Examples of the shapes and norming constants of several corner sets are given in Figure \ref{fig:norming d=2} of the supplement. 

Note that the proposed parent sets in \eqref{eq:parent grid} have cardinality ${l+d \choose l}$. Moreover, the corresponding polynomial interpolation problem is unisolvent and by the multivariate divided-difference formula (see Theorem 5 of \cite{neidinger2019multivariate}), admits an explicit representation.
Consequently, in view of Lemma \ref{lem:norming def}, each of such parent set is a norming set. Although the norming constant cannot be written as a simple explicit function of $d$ and $l$, there are only finitely many possible shapes for the corner sets, and each possess a finite norming constant. Taking their maximum yields an absolute norming constant that applies uniformly to all possible corner sets, thereby establishing Condition \ref{cond:norming}.

\begin{figure}
  \centering
\subfloat[Layer \textcolor{red}{$\N_0$}]{\includegraphics[width=0.26\linewidth]{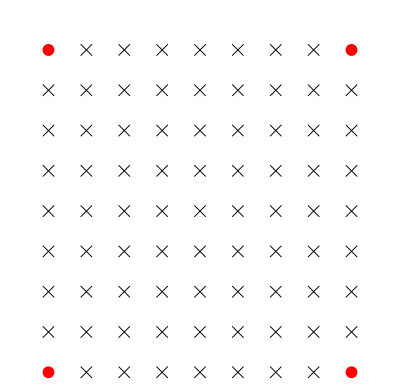}}
\subfloat[Layer \textcolor{blue}{$\N_0$}, \textcolor{red}{$\N_1$}]{\includegraphics[width=0.26\linewidth]{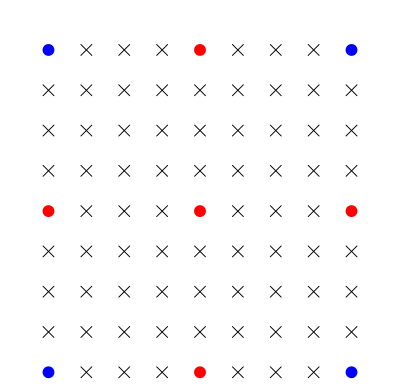}}
\subfloat[Layer \textcolor{blue}{$\N_0$, $\N_1$}, \textcolor{red}{$\N_2$}]{\includegraphics[width=0.26\linewidth]{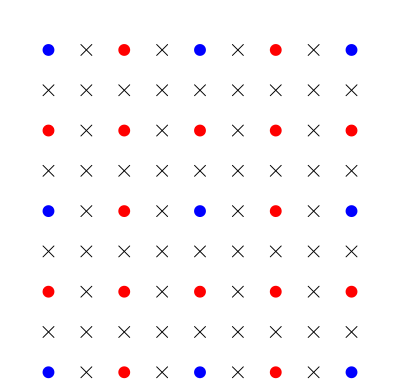}}
\caption{Illustration of layers on a $9\times 9$ grid: red dots: current layer; blue dots: all previous layers; black crosses: all latter layers.}\label{fig:layers}
\end{figure}

\begin{figure}
    \centering    \includegraphics[width=0.79\linewidth]{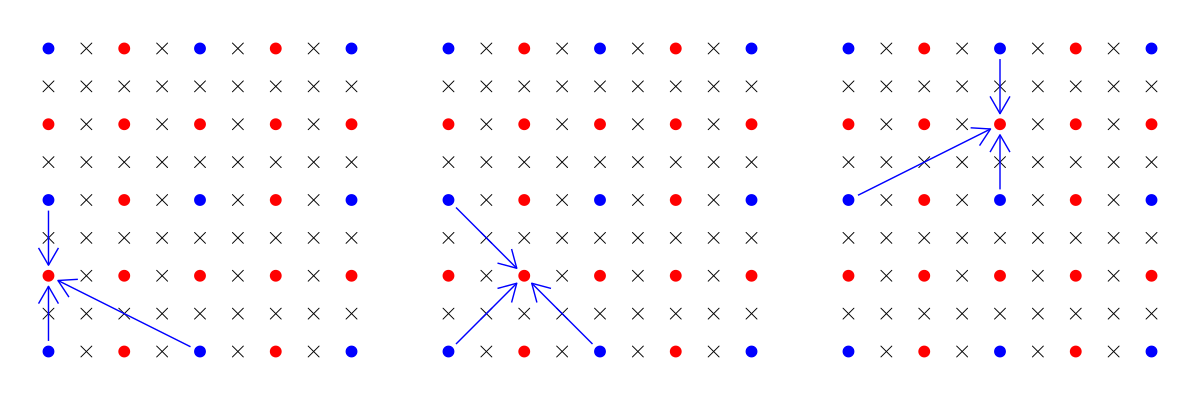}
\caption{Continuing the example in Figure \ref{fig:layers}, illustration of parent sets for $X_i\in \N_2$, with $\underline{\alpha}=1$. Red dots: current layer $\mathcal{N}_2$; Blue dots: previous layers $\mathcal{N}_0$, $\N_1$; Black crosses: all latter layers.
Blue arrows: directed edges from parent sets to children for some $X_i\in \N_2$.}\label{fig:parents} 
\end{figure}

Finally, for $\tilde{n}\ne 2^{r'}+1$, $\forall r'\in\mathbb{N}$, let $r$ be the largest integer such that $\tilde{n}> 2^r+1$. Let $\mathscr{I}\subset\{0,1,\cdots, \tilde{n}-1\}$ be an index set of cardinality $2^r+1$ such that $\{0,1,\cdots, \tilde{n}-1\}\backslash \mathscr{I}$ does not contain consecutive integers. Then we consider the subset of $\X_n$ defined by
\begin{equation}\label{eq:subgrid}
\bar{\X}_n = \big\{x=(x[1], x[2], \cdots, x[d])^T: x[h] \in \{k/(\tilde{n}-1): k\in \mathscr{I}\}, \forall 1\le h\le d \big\}.
\end{equation}
Let $\mathcal{M}$ be an order-preserving mapping from $\{k/(\tilde{n}-1): k\in \mathscr{I}\}$ to $\{k/2^r: 0\le k \le 2^r\}$. Then the set
$$\hat{\X}_n = \big\{ (\mathcal{M}(x[1]), \mathcal{M}(x[2]), \cdots, \mathcal{M}(x[d]))^T: (x[1], x[2], \cdots, x[d])^T\in \bar{X}_n\big\}$$
lies on a $d$-dimensional grid of cardinality $(2^r+1)^d$ and has a one-to-one correspondence with $\bar{X}_n$ induced by the coordinate-wise mapping $\mathcal{M}$. Thus, one can build a layered norming DAG on $\hat{\X}_n$ using the method described previously, which subsequently induces a layered norming DAG on $\bar{X}_n$. Once this is done, we define the layer $\N_{r+1} = \X_n\backslash\bar{\X}_n$. For all $X_i\in \N_{r+1}$ the parent set is defined as a corner set as in equation (\ref{eq:parent grid}). We formalize this procedure in Algorithm 1, which is deferred to Section \ref{sec:build dag alg} of the Supplementary Material. The DAG on a $2$-dimensional grid is visualized with the layer structures in Figure \ref{fig:layers} and the parent sets in Figure \ref{fig:parents}.


\section{Probabilistic properties}\label{sec:prob}

In this section, we study the probabilistic properties of the proposed Vecchia approximation of Mat\'{e}rn processes introduces in Section \ref{sec:nota}.
Standard techniques for analyzing Mat\'{e}rn GPs rely heavily on the explicit form of the covariance kernel and its Fourier transform. Unfortunately, these techniques are not applicable to Vecchia GPs. First, the covariance matrices of their finite dimensional marginals can only be obtained by inverting the precision matrices, rendering the direct analysis substantially more difficult. Moreover, the covariance kernel of the Vecchia approximated Mat\'{e}rn process is not stationary, hence techniques based on spectral density are no longer adequate. This motivates the development of new mathematical tools that exploit the finite-dimensional conditional distributions of the GP.

The section is organized as follows. In Section \ref{sec:matern theory}, we collect several results regarding the derivatives of Mat\'{e}rn covariance functions. Section \ref{sec:cond} studies the conditional distribution of $Z^{\tau,s}_X$ given the process values at $\mr{pa}(X)$. Section \ref{sec:small ball} utilizes the results from Section \ref{sec:cond} to derive small deviation bounds for the Vecchia GP $\hat{Z}^{\tau,s}_X$ associated with the Mat\'{e}rn mother GP. Finally, Section \ref{sec:decentering} investigates the Reproducing Kernel Hilbert Space (RKHS) associated with Vecchia GPs and their approximation properties, which are then used to derive lower bounds for the decentered small ball probabilities of Vecchia GPs. Throughout the section, whenever a DAG is involved, we assume that Conditions \ref{cond:layered} and \ref{cond:norming} are satisfied with some positive constants $c_d, c_L$ and $c_N$ that do not depend on the sample size $n$. 

\subsection{Derivatives of the Mat\'{e}rn covariance kernel}\label{sec:matern theory}
The differentiability properties of Mat\'{e}rn covariance functions are well-known \citep{williams2006gaussian,ghosal2017fundamentals}. However, the Lipschitz continuity result given in equation (\ref{eq:matern diff lipschitz}) below, is less commonly discussed. In particular, only the case $\alpha\not\in\mathbb{N}$ is considered in \citep{JMLR:v12:vandervaart11a}, and no proofs are provided. For clarity and completeness, we collect all relevant results here and provide the corresponding proofs in the supplement.

\begin{lemma}\label{lem:matern smooth}
Let $K(\cdot, \cdot)$ denote the Mat\'{e}rn covariance kernel with regularity parameter $\alpha>0$, see \eqref{eq:MaternCov}. Then $K(\cdot,\cdot)$ is $2\underline{\alpha}$ times differentiable on $\mathbb{R}^{2d}$, such that $\forall \;k_1,k _2\in\mathbb{N}^d$, $|k_1|+|k_2|\le 2\underline{\alpha}$ and $\forall \; x_1,x_2\in\mathbb{R}^d$, 
\begin{equation}\label{eq:matern diff}
|K^{(k_1,k_2)}(x_1,x_2)|\lesssim 1.
\end{equation}
Furthermore, for all $k\in\mathbb{N}^d$, $|k|\le\underline{\alpha}$  and $\forall \; x_1,x_2\in\mathbb{R}^d$,
\begin{equation}\label{eq:matern diff symmetry}
K^{(k,k)}(x_1,x_2) = K^{(k,k)}(x_2,x_1),\;\; K^{(k,k)}(x_1,x_1) = K^{(k,k)}(x_2,x_2).
\end{equation}
Finally, for all $|k_1|+|k_2|=2\underline{\alpha}$, $x\in \mathbb{R}^d$ and $h\in\mathbb{R}^d$ small enough,
\begin{equation}\label{eq:matern diff lipschitz}
|K^{(k_1,k_2)}(x,x)-K^{(k_1,k_2)}(x,x+h)|\lesssim \bigg\{
\begin{aligned}
&\|h\|_2^{2(\alpha-\underline{\alpha})}, \;\; &&\alpha\not\in\mathbb{N},\\
&\|h\|^2_2\ln(1/\|h\|_2), \;\; && \alpha\in\mathbb{N}.
\end{aligned}
\end{equation}
\end{lemma}

Lemma \ref{lem:matern smooth} shows that the Mat\'{e}rn covariance function $K(\cdot,\cdot)$, viewed as a map from $\mathbb{R}^{2d}\to\mathbb{R}$, is $2\underline{\alpha}$ times differentiable, with its $2\underline{\alpha}$th-order derivatives being Lipschitz continuous. The Mat\'{e}rn covariance function belongs  to $C^\alpha(\mathbb{R}^d)\times C^\alpha(\mathbb{R}^d)$ in the sense that, for integer $\alpha$, it is weakly differentiable of order $(k,k),\forall |k|=\alpha$. However, its H\"{o}lder-$(\alpha,\alpha)$ norm is infinite, since the weak derivatives of order $(k,k),\forall |k|=\alpha$ are not uniformly bounded. Moreover, the sample paths  of the process  are almost surely $\underline{\alpha}$ times differentiable with their $\underline{\alpha}$th derivatives being almost surely Lipschitz continuous of order $\alpha-\underline{\alpha}$, up to a logarithmic term.

The smoothness properties described in Lemma \ref{lem:matern smooth} also carry over to the covariance kernels of the conditional processes. More precisely, for an arbitrary finite set $A=\{w_1,w_2,\cdots,w_m\}\subset \mathbb{R}^d$, we define the process $\tilde{Z}$ as
$$\tilde{Z}_x = Z_x - \mathbb{E}(Z_x|Z_A) = Z_x - K_{x,A} K_{A,A}^{-1} Z_A.$$
By the linearity, $\tilde{Z}$ is also a centered GP with covariance function
\begin{equation}\label{eq:K tilde}
\tilde{K}(x_1,x_2) = \mathbb{E} \tilde{Z}_{x_1} \tilde{Z}_{x_2} = K(x_1,x_2) - K_{x_1,A} K_{A,A}^{-1} K_{A,x_2}.
\end{equation}
The process $\tilde{Z}$ can be regarded as the Gaussian process conditioned on the set $A$. We show below that it inherits several key properties of the mother Mat\'{e}rn process $Z$. Remarkably, these properties do not depend on the choice of the conditioning set $A$.

\begin{lemma}\label{lem:matern cond smooth}
Let $\tilde{K}(\cdot, \cdot)$ be the covariance function defined in (\ref{eq:K tilde}). Then $\tilde{K}(\cdot,\cdot)$ is $2\underline{\alpha}$ times differentiable on $\mathbb{R}^{2d}$ and for all $k\in\mathbb{N}^d$, $|k|\le\underline{\alpha}$, and $x_1,x_2\in\mathbb{R}^d$,
\begin{equation}\label{eq:matern cond diff}
|\tilde{K}^{(k,k)}(x_1,x_2)|\lesssim 1.
\end{equation}
Moreover, for all $|k|=\underline{\alpha}$, $x\in \mathbb{R}^d$ and $h\in\mathbb{R}^d$ small enough, we have
\begin{align}\label{eq:matern cond diff lipschitz}
&\big|\tilde{K}^{(k,k)}(x,x)-\tilde{K}^{(k,k)}(x,x+h)-\tilde{K}^{(k,k)}(x+h,x)+\tilde{K}^{(k,k)}(x+h,x+h)\big|\nonumber\\
&\qquad\quad\lesssim \bigg\{
\begin{aligned}
&\|h\|_2^{2(\alpha-\underline{\alpha})}, \;\; &&\alpha\not\in\mathbb{N},\\
&\|h\|_2^2\ln(1/\|h\|_2), \;\; && \alpha\in\mathbb{N}.
\end{aligned}
\end{align}
All constant multipliers in the upper bounds are independent of  $x,h \in \mathbb{R}^d$ and the choice of the set $A$.
\end{lemma}

The universality of the constants in Lemma \ref{lem:matern cond smooth} across different conditioning sets $A$ is a consequence of the fact that the variance of the conditional process $\tilde{Z}$ does not exceed the variance of its mother GP $Z$. In view of the relationship between the derivatives of GPs and the derivatives of the corresponding covariance functions, as established in Proposition I.3 of \cite{ghosal2017fundamentals}, this reasoning extends to the derivatives of the
GPs $Z$ and $\tilde{Z}$, up to order $\underline{\alpha}$. Therefore, regardless of the conditioning set $A$, the conditional process $\tilde{Z}$ inherits similar smoothness properties as the mother GP $Z$. This feature plays a crucial role in analyzing the conditional distribution in the subsequent sections. We also note that the logarithmic term in the preceding lemmas for $\alpha\in\mathbb{N}$ arises from the singularity of the modified Bessel functions of second kind at integer parameters. For notational simplicity, throughout the rest of the paper we write $[\ln(1/\|h\|_2)]^{\mathbbm{1}_{\{\alpha\in\mathbb{N}\}}}$ to indicate that the logarithmic term  appears only when $\alpha$ is an integer.

\subsection{Conditional distribution}\label{sec:cond}
The conditional distribution formulas (\ref{eq:vecchia 1}) and (\ref{eq:vecchia 2}) defining the Vecchia GPs form the starting point of our theoretical analysis. In this section, we study three aspects of these conditional distributions: the variance in equation (\ref{eq:vecchia var}), the conditional expectation in equation (\ref{eq:vecchia exp}), and the recursive application of this conditional expectation formula.

\subsubsection{Variance and conditional expectation}
Under the Norming condition \ref{cond:norming}, the distance among vertices in the set $\mr{pa}(X_i)\cup\{X_i\}$ is bounded by $\gamma^{-\eta(i)}$ up to a constant multiplier, where $\eta(i)=j$ denotes the layer allocation of $X_i\in\mathcal{N}_j$. We recall that the observations are ordered based on their layer allocations, hence the function $i\rightarrow \eta(i)$ is monotone increasing. We investigate the asymptotic behavior of the Gaussian random variable $Z^{\tau,s}_{X_{i}}-\mathbb{E}[{Z}^{\tau,s}_{X_{i}}|{Z}^{\tau,s}_{\mr{pa}(X_{i})}]$ as $i$ tends to infinity. The following theorem provides lower and upper bounds for its variance, which, by equations \eqref{eq:vecchia 2} and \eqref{eq:vecchia var}, coincides with the variance of the conditional Mat\'ern GP process.

\begin{theorem}\label{thm:var}
Let the mother Gaussian process $Z^{\tau,s}$ be the rescaled Mat\'{e}rn process with smoothness parameter $\alpha$. Suppose that Condition \ref{cond:norming} is satisfied. Then, for all $X_i\in \mathcal{N}_j$, $j=1,...,J$, we have
\begin{align}\label{eq:lem:var}
\var \big\{Z^{\tau,s}_{X_{i}}-\mathbb{E}[{Z}^{\tau,s}_{X_{i}}|{Z}^{\tau,s}_{\mr{pa}(X_{i})}]\big\} \asymp 
 s^2\tau^{2\alpha}\gamma^{-2\alpha j} (j\ln \gamma)^{\mathbbm{1}_{\{\alpha\in\mathbb{N}\}}}.
\end{align}
\end{theorem}

Theorem \ref{thm:var} plays an essential role in determining the small ball probabilities of Vecchia GPs. Specifically, for a finite training set $\X_n$, the number of layers $\eta(n)$ satisfies $\gamma^{d\eta(n)}\asymp n$. Therefore, for $X_i\in\N_{\eta(n)}$ and $\tau=s=1$, 
$ \var \big\{Z_{X_{i}}-\mathbb{E}[{Z}_{X_{i}}|{Z}_{\mr{pa}(X_{i})}]\big\} \asymp n^{- 2\alpha/d}.$ In other words,
the variance at $X_i\in\N_{\eta(n)}$ decays to zero with a polynomial rate in $n$, where the order is determined by the smoothness parameter $\alpha$ of the mother Mat\'{e}rn process. Hence smoother mother GPs exhibit faster decaying variances, as indicated in \eqref{eq:lem:var}.

We proceed to study the conditional expectation $\mathbb{E}[\hat{Z}^{\tau,s}_{X_i}|\hat{Z}^{\tau,s}_{\mr{pa}(X_i)}]$ in equation (\ref{eq:vecchia 2}), which corresponds to Gaussian  interpolation, also called kriging in the spatial statistics literature \cite{gelfand2010handbook}, of $\hat{Z}^{\tau,s}_{X_i}$ using the random vector $\hat{Z}^{\tau,s}_{\mr{pa}(X_i)}$. 
Recalling Condition \ref{cond:norming}, for all $i\ge i_0$ and $X_i\in \N_{\eta(i)}$, its parent set $\mr{pa}(X_i)$ is contained in a cube of side length $c_L \gamma^{-\eta(i)}$ and $\mr{pa}(X_i)$ forms a norming set in this cube with norming constant $c_N$. Therefore, given an arbitrary $X_i$, there exists a norming set $A=\{w_i, 1\le i\le m\}\subset [0,1]^d$ with norming constant $c_N$ and a vector $x'\in\mathbb{R}^d$ such that, for any rescaling parameter $\tau>0$,
$$\tau X_i = \tau c_L \gamma^{-\eta(i)}x^* + \tau x', \;\; \tau \mr{pa}(X_i) = \tau c_L \gamma^{-\eta(i)} A + \tau x'.$$
In view of the stationarity of the Mat\'ern process, the conditional expectation in (\ref{eq:vecchia 2}) can then be expressed as
\begin{align}
\mathbb{E}[\hat{Z}^{\tau,s}_{X_i}|\hat{Z}^{\tau,s}_{\mr{pa}(X_i)}] & = K^T_{\tau \mr{pa}(X_i), \tau X_i} K^{-1}_{\tau \mr{pa}(X_i), \tau \mr{pa}(X_i)} \hat{Z}^{\tau,s}_{\mr{pa}(X_i)} \nonumber\\
& = K^T_{\tau c_L \gamma^{-\eta(i)} A, \tau c_L\gamma^{-\eta(i)} x^*} K^{-1}_{\tau c_L\gamma^{-\eta(i)} A, \tau c_L\gamma^{-\eta(i)} A} \hat{Z}^{\tau,s}_{\mr{pa}(X_i)}\nonumber\\
&=K_{\nu A,\nu x^*}^T K_{\nu A,\nu A}^{-1} \hat{Z}^{\tau,s}_{\mr{pa}(X_i)},\label{eq:krig}
\end{align}
where $\nu = \tau c_L \gamma^{-\eta(i)} \asymp \tau \gamma^{-\eta(i)}$. Note that the conditional expectation above depends on the covariance kernel $K$, which characterizes the spatial dependencies of the process $Z$. However, we show below, that under certain regularity conditions,  for fixed $x^*$ and $A$, as $\nu$ tends to zero, the limit for the Gaussian interpolation weights $K_{\nu A,\nu A}^{-1} K_{\nu A,\nu x^*}$ becomes independent of the covariance function. This limit is often referred to as the \textit{flat limit} in the radial basis function literature \citep{driscoll2002interpolation,song2012multivariate,lee2015study}. 
For Mat\'{e}rn covariance functions, the following lemma shows that the flat limit coincides with the {polynomial interpolation weights $V_A^{-1} v_x^*$}, provided that the set $A$ is a norming set with respect to an appropriately chosen family of polynomials.

\begin{lemma}\label{lem:GP flat}
Let $A$ be a norming set on $\mathscr{P}_{\underline{\alpha}}([0,1]^d)$ with norming constant $c_N$. Suppose $x\in [0,1]^d$, then for all $\nu
\le 1$,
\begin{equation}\label{eq:GP flat}
\big\| K_{\nu A,\nu A}^{-1}K_{\nu A,\nu x} - V_A^{-1} v_{x}\big\|_1 \lesssim \nu^{2(\alpha-\underline{\alpha})} [\ln (1/\nu)]^{\mathbbm{1}_{\{\alpha\in\mathbb{N}\}}}  + \nu,
\end{equation}
where $V_A$ and $v_x$ are defined in \eqref{def:vandemonde} and \eqref{eq:v vector}. Similarly, for $\nu,\nu_0\in(0,+\infty)$ and $|\nu-\nu_0|\leq 1$,
\begin{equation}\label{eq:GP flat 2}
\big\| K_{\nu A,\nu A}^{-1}K_{\nu A,\nu x} - K_{\nu_0 A,\nu_0 A}^{-1}K_{\nu_0 A,\nu_0 x}\big\|_1 \lesssim |\nu-\nu_0|^{2(\alpha-\underline{\alpha})} [\ln (1/|\nu-\nu_0|)]^{\mathbbm{1}_{\{\alpha\in\mathbb{N}\}}}.
\end{equation}
The constant multipliers in the upper bounds are free of $x$ and $A$, but depend on $c_N$.
\end{lemma}

Lemma \ref{lem:GP flat} provides an example of flat limits in radial basis function interpolation, a topic with a long-standing literature spanning over twenty years, see for instance \cite{driscoll2002interpolation,lee2015study}. The general idea is that, for sufficiently smooth kernels (i.e. Mat\'{e}rn covariance functions in our case) and an unisolvent set $A$ with respect to a given family of polynomials, the limit of equation (\ref{eq:krig}) coincides with the limit of polynomial interpolation within that family. With this in mind, we now discuss the conditions of Lemma \ref{lem:GP flat}. 
For $A$ to serve as a norming set of a polynomial space, its cardinality has to be at least equal to the dimension of the polynomial space. The covariance function $K(\cdot,\cdot)$  of the Mat\'{e}rn process with regularity $\alpha$ belongs to $C^{2\underline{\alpha}}$, i.e. it is $\underline{\alpha}$ times differentiable in each covariate. 
Intuitively, the polynomial space describing the local behavior of Mat\'{e}rn covariance function is of order $\underline{\alpha}$, which requires the set $A$ to be a norming set for $\mathscr{P}_{\underline{\alpha}}([0,1]^d)$.

In view of Lemma \ref{lem:GP flat}, the norming constant $c_N$ appears in the upper bounds, making it important to have a tight control over it. However, as we discuss next, this is generally a highly challenging task. Recall that $A$ is a norming set if it is unisolvent with respect to a given family of polynomials. In the unidimensional case, a set $A$ with $m$ distinct elements is always unisolvent with respect to  polynomials of degree up to $(m-1)$th. In higher dimensions ($d\geq 2$), however, it is possible for the elements of $A$ to lie on a lower-dimensional polynomial manifold, in which case the corresponding Vandermonde matrix $V_A$ is singular. In such cases, the limit of Gaussian interpolation on $A$ is interpreted within a family of lower-order polynomials. We note that the assumption that $A$ is exactly a unisolvent set of $\mathscr{P}_l([0,1]^d)$ is not strictly necessary. The limit of the expression in (\ref{eq:krig}) still exists for any choice of $m$ satisfying ${l+d\choose l}\le m <{l+1+d\choose l+1}$. However, the corresponding, asymptotically equivalent formulation is more complex, as it depends on the derivatives of the covariance function at the origin, see \citep{lee2015study}. For instance, if $A$ is a strict superset of a norming set and $\alpha$ is an even integer, the limit corresponds to polyharmonic spline interpolation \cite{song2012multivariate}. We are not aware of any result addressing this limit in a more general setting. Therefore, for simplicity, we restrict our discussion to unisolvent sets from now on.

Finally, we note that our Lemma \ref{lem:GP flat} extends the results of \cite{lee2015study} by providing the explicit convergence rate of the Gaussian interpolation to the polynomial interpolation. Moreover, these rates are uniform in $x\in\mathbb{R}^d$ if $A$ is an unisolvent set with a universal norming constant $c_N$. Such explicit, tight control is needed for the small deviation bounds of the Vecchia GP, which in turn are required to obtain optimal posterior contraction rates.

\subsubsection{Recursive interpolation}\label{sec:rec:inter}
Another important aspect of Vecchia GPs is the repeated application of the conditional expectation formula \eqref{eq:krig}. Specifically, for all $X_i\in \X$, we define the \textit{double parent} set of $X_i$ as $\mr{pa}^2(X_i) \triangleq \cup_{X\in\mr{pa}(X_i)}\mr{pa}(X)$. 
Although the conditional distributions of $Z^{\tau,s}_{X_i}$ and its Vecchia approximation $\hat{Z}^{\tau,s}_{X_i}$ coincide when conditioned on the processes at the corresponding parent set $\mr{pa}(X_i)$, this is generally not the case when conditioning on the processes at the double parent set $\mr{pa}^2(X_i)$. In other words, applying the conditional expectation formula of Vecchia Gaussian processes more than once typically produces a conditional distribution that differs from that of the mother GP.

To formally study the properties of repeated Gaussian interpolation, we introduce the operator $G_j$ that maps a function $f:\cup_{j'=0}^j\N_{j'}\to\mathbb{R}$ to a function $G_j(f):\cup_{j'=0}^{j+1}\N_{j'}\to\mathbb{R}$, such that
\begin{equation}\label{def:op:G}
G_j(f)(X) = \Bigg\{ 
\begin{aligned}
&f(X), &&\forall X\in \cup_{j'=0}^j\N_{j'}, \\
&K_{\tau\mr{pa}(X),\tau X}^T K_{\tau\mr{pa}(X),\tau\mr{pa}(X)}^{-1} \bm{f}({\mr{pa}(X)}), &&\forall X\in \N_{j+1},\end{aligned}
\end{equation}
where $\bm{f}(\mr{pa}(X))$ denotes the vector $(f(x))_{ x\in \mr{pa}(X)}$, see Section \ref{sec:notation set}. The operator $G_j$ performs Gaussian interpolation at locations $X\in \N_{j+1}$ using the function values at the lower layers $ \cup_{j'=0}^j\N_{j'}$.
In view of Lemma \ref{lem:GP flat}, the limit of the Gaussian interpolation corresponds to polynomial interpolation.
Let us denote the associated polynomial interpolation operator by $P_j$, mapping a function  $f:\cup_{j'=0}^j\N_{j'}\to\mathbb{R}$ to a function $P_j(f):\cup_{j'=0}^{j+1}\N_{j'}\to\mathbb{R}$ as
$$
P_j(f)(X) = \Bigg\{ 
\begin{aligned}
&f(X), &&\forall X\in \cup_{j'=0}^j\N_{j'}, \\
&v_{\tau X}^T V_{\tau\mr{pa}(X)}^{-T} \bm{f}({\mr{pa}(X)}) = v_{X}^T V_{\mr{pa}(X)}^{-T} \bm{f}({\mr{pa}(X)}), &&\forall X\in \N_{j+1}.\end{aligned}
$$
Furthermore, let us endow the space of functions on $\cup_{j'=1}^j \N_{j'}$ with the (empirical) supremum norm $\|f\|_{\infty,j}=\sup_{X\in \cup_{j'=1}^j \N_{j'}}f(X)$. For the operator $P_j$, which maps functions from the space $L_\infty(\cup_{j'=1}^j \N_{j'})$ to $L_\infty(\cup_{j'=1}^{j+1} \N_{j'})$, we define its (empirical) $L_\infty$ operator norm as
$$
\|P_j\| \triangleq \sup_{f \in L_\infty(\cup_{j'=1}^j \N_{j'}), f\ne 0} \frac{\|P_j(f)\|_{\infty,j+1}}{\|f\|_{\infty,j}}.
$$
We define the supremem norm for $G_{j}$ and their compositions in the same way. 
Then, in view of Lemma \ref{lem:GP flat}, the difference between the operators $G_j$ and $P_j$ tends to zero in operator norm. Therefore, to obtain an upper bound for the operator norm of the recursive Gaussian interpolation operator $G_{j_2} G_{j_2-1} \cdots G_{j_1} (f)$, for some $j_1<j_2,$ it is sufficient to control the corresponding recursive interpolation using the polynomial operators.

To this end, let us introduce the notation
\begin{equation}\label{eq:poly interp cond}
\vartheta_n = \sup\left\{ \frac{\|P_{j_2} P_{j_2-1} \cdots P_{j_1} (f)\|_{\infty,j_1+1}}{\|f\|_{\infty,j}}: f \in L_\infty(\cup_{j'=1}^j \N_{j'}) , 0\le j_1\le j_2\le \eta(n)\right\}.
\end{equation}
Due to the close connection between Gaussian and polynomial interpolation, this term plays a crucial role in deriving small deviation bounds for Vecchia GPs. Note that it depends both on the data set $\X_n$ and on the DAG constructed over it. 

\begin{remark} In our analysis, we consider grid data \eqref{def:grid}. In Lemma \ref{lem:poly recursive}, deferred to Section \ref{sec:vartheta} of the Supplementary Material, we obtain tight control of $\vartheta_n$ by showing that it is bounded for grid data \eqref{def:grid} and layered norming sets. Although considerable numerical evidence suggests that $\vartheta_n=O(1)$ holds more generally, we do not have a rigorous proof of this claim. The difficulty lies in tightly controlling the multiple ($n$ dependent) recursive applications of the operators $P_j, j\in\mathbb{N}$, which is not directly implied by the uniform control of  the individual operators. Specifically, for the considered layered norming DAGs, the total number of layers satisfy $\eta(n) \asymp \ln n$. Therefore, even if we have uniform control over the operators $\|P_j\|\le c,\forall j\in\mathbb{N}$ for some constant $c> 1$, this still yields a suboptimally large upper bound for the repeated interpolation 
$$\|P_{\eta(n)-1} P_{\eta(n)-2} \cdots P_1 P_0\| \le c^{c'\ln n}=n^{c'\ln c}.$$
Therefore, it is necessary to consider the composition of these operators as a whole, which is exceedingly difficult without any additional structural or geometrical assumptions on $\X_n$. The problem becomes considerably simpler for grid data. In the proof of Lemma \ref{lem:poly recursive} we show that, in this special case the polynomial interpolation is equivalent to convolution, which can be further converted to multiplication via the Fourier transform. Semi-explicit computations then verify that $\vartheta_n=O(1)$.
\end{remark}


Then, by exploiting the close asymptotic equivalence of polynomial and Gaussian interpolations, see Lemma \ref{lem:GP flat}, we can bound the operator norm of recursive Gaussian interpolation using $\vartheta_n$ as defined in \eqref{eq:poly interp cond}.
\begin{lemma}\label{lem:GP recursive}
Under Condition \ref{cond:norming}, for all $0\le j_1\le j_2$, 
$$\|G_{j_2} \cdots G_{j_1+1} G_{j_1}\|\lesssim \tilde{\vartheta}_n,$$
where $\tilde{\vartheta}_n = \vartheta_n \exp(c\vartheta_n)$ for some absolute constant $c>0$.
\end{lemma}
Note that $\vartheta_n =O(1)$ implies $\tilde{\vartheta}_n=O(1)$, which holds for grid data. In the following, recursive Gaussian interpolation appears in several technical lemmas. For simplicity, in the rest of the remainder, whenever recursive Gaussian interpolation is involved, we explicitly include the $\tilde\vartheta_n$ term in the upper bound. This allows the results to be applied beyond the grid data case, as long as $\vartheta_n$ or, more generally $\tilde{\vartheta}_n$, can be controlled.

We conclude this section with some general remarks. Theorem \ref{thm:var} and Lemma \ref{lem:GP flat} study the variance and conditional expectations that define the Vecchia GPs, while Lemma \ref{lem:GP recursive} provides an upper bound for the operator norm of the recursive Gaussian interpolation. Together, these three results form the building blocks for deriving small deviation bounds for GPs, controlling the metric entropy of the associated RKHS and ultimately deriving posterior contraction rates for Vecchia GPs.

\subsection{Small ball probability}\label{sec:small ball}
In this section, we study the small ball probability of the Vecchia Gaussian process with respect to the supremum  on $\X_n$, defined as
$$\pr(\|\hat{Z}^{\tau,s}\|_{\infty,n} <\epsilon),\;\;\forall \epsilon\in(0,1).$$
Small ball probability is directly linked to the $\epsilon$-entropy of the associated RKHSs (Lemma I.29 and I.30 of \cite{ghosal2017fundamentals}), provides lower bounds for empirical processes (Section 7.3 of \cite{li2001gaussian}) and plays a crucial role in the convergence rate of the corresponding posterior distribution \citep{van2008rates,JMLR:v12:vandervaart11a}. For Brownian motion, the small ball probability can also be interpreted as the first exit time. For more details and applications, see for instance the survey \cite{liu2020gaussian} and references therein.

 The following lemma provides a lower  bound for the small ball probability, or equivalently, an upper bound for the small ball exponent $-\ln \pr(\|\hat{Z}^{\tau,s}\|_{\infty,n} <\epsilon)$ of the Vecchia GP.

\begin{theorem}[Small deviation bound]\label{lem:small ball}
Under Conditions \ref{cond:layered} and  \ref{cond:norming}, if the mother GP is a Mat\'{e}rn process with smoothness parameter $\alpha$, then for sufficiently small $\epsilon>0$,
$$-\ln \pr(\|\hat{Z}^{\tau,s}\|_{\infty,n} < \tilde\vartheta_n \epsilon) \lesssim \tau^d s^{d/\alpha}\epsilon^{-d/\alpha} [\ln (1/\epsilon)]^{\mathbbm{1}_{\{\alpha\in\mathbb{N}\}}},$$
where the constant in the upper bound does not depend on the sample size.
\end{theorem}

The proof of Theorem \ref{lem:small ball} follows a different route from the standard techniques used to obtain small deviation bounds of stationary Gaussian processes. Instead of deriving the small ball probability via the $\epsilon$-entropy of the corresponding Reproducing Kernel Hilbert Space, we study it by exploiting the conditional distribution results derived in the previous section. The techniques draw inspiration from Chapter 7 of \cite{ledoux2006isoperimetry} and are generalized in our analysis. We believe these techniques are of independent interest and may be applicable to a broader class of stochastic processes.

For scaling parameters $\tau=s=1$, the lower bound for the small ball probability in Theorem \ref{lem:small ball} coincides with that of the Mat\'{e}rn process, see Lemma 11.36 of \cite{ghosal2017fundamentals}. This indicates that, despite having different covariance structures and lacking stationarity, the Vecchia approximation of a Mat\'{e}rn process exhibits similar small deviation behavior as the mother process. Such property is crucial for deriving matching, optimal posterior contraction rates for the two processes.

\subsection{Decentered small ball probability}\label{sec:decentering}
In this section we extend the above small deviation bound to the decentered case, i.e. we provide lower bounds for the decentered small ball probability
$$\pr(\|\hat{Z}^{\tau,s}-f_0\|_{\infty,n} <\epsilon),\;\; \epsilon\in(0,1),$$
for $f_0\in L_\infty(\Omega)$. We begin by recalling the definition of the RKHS associated with the Vecchia GP and discuss how its approximation properties relate to decentered small ball probabilities.

Let us recall that the  RKHS of $\hat{Z}^{\tau,s}_{\X_n}$ is
$$
\mathbb{H}_n^{\tau,s} = \Big\{\sum_{i=1}^n a_i \hat{K}^{\tau,s}(X_i,\cdot): X_i\in\X_n, a_i\in\mathbb{R}\Big\}.
$$
For any function of the form $f=\sum_{i=1}^n a_i \hat{K}^{\tau,s}(X_i,\cdot)\in\mathbb{H}_n^{\tau,s}$, let $a = (a_1,a_2,\cdots, a_n)^T$. Then the RKHS norm of $f$ is given by
\begin{equation}\label{eq:f RKHS def}
\|f\|_{\mathbb{H}_n^{\tau,s}}^2 = a^T \hat{K}^{\tau,s}_{\X_n,\X_n} a.
\end{equation}
Furthermore, note that
\begin{equation*}
\bm{f}(\X_n) = \hat{K}^{\tau,s}_{\X_n,\X_n} a.
\end{equation*}
By combining the previous two equations,
\begin{align}\label{f:RKHS decom}
\|f\|_{\mathbb{H}_n^{\tau,s}}^2
= & \bm{f}(\X_n)^T (\hat{K}^{\tau,s}_{\X_n,\X_n})^{-1} \bm{f}(\X_n). 
\end{align}
We can now define the decentering function as
\begin{equation*}
\epsilon\mapsto \inf_{f\in\mathbb{H}_n^{\tau,s}:\|f-f_0\|_{\infty,n}\le \epsilon} \|f\|_{\mathbb{H}_n^{\tau,s}}^2.
\end{equation*}
The combination of the decentering function and the centered small ball exponent yields the so called \textit{concentration function}
\begin{equation}\label{def:conc:funct}
\phi_{f_0,n}^{\tau,s}(\eps)=\inf_{f\in\mathbb{H}_n^{\tau,s}:\|f-f_0\|_{\infty,n}\le \epsilon} \|f\|_{\mathbb{H}_n^{\tau,s}}^2-\ln \pr(\|\hat{Z}^{\tau,s}\|_{\infty,n} <\epsilon).
\end{equation}
Finally, in view of Proposition 11.19 of \cite{ghosal2017fundamentals}, for any $f_0$ in the closure of the RKHS with respect to the empirical supremum norm $\|\cdot\|_{\infty,n}$, the decentered small ball probability satisfies
\begin{align}
\phi_{f_0,n}^{\tau,s}(\eps)\leq -\ln \pr(\|\hat{Z}^{\tau,s}-f_0\|_{\infty,n} <\epsilon)\leq \phi_{f_0,n}^{\tau,s}(\eps/2).\label{eq:conc:ball}
\end{align}
Therefore, understanding the asymptotic behavior of the decentering function (as $\eps=\eps_n$ tends to zero with $n$) will provide us the decentered small ball probabilities.

To derive an upper bound for the decentering function, we first present several technical lemmas that are of independent interest. The first  lemma discusses the minimal RKHS norm among all functions that interpolate  a given set of points.
\begin{lemma}\label{lem:RKHS min}
Let $Z$ be a Gaussian process defined on $\X$ and $\mathbb{H}$ be the corresponding RKHS. Then, for any finite set $A=\{x_1, x_2, \cdots, x_{|A|}\}$ and any function $f_0:\X\to \mathbb{R}$,
$$\inf_{f\in\mathbb{H},\bm{f}(A)=\bm{f}_0(A)} \|f\|_{\mathbb{H}}^2 = \bm{f}_0(A)^T K_{A,A}^{-1} \bm{f}_0(A).$$
The minimum in the above equation is obtained by the function $f(\cdot)=\sum_{i=1}^{|A|} a_i K(x_i,\cdot)$, with coefficients $a = (a_1,a_2,\cdots, a_{|A|})^T = K_{A,A}^{-1} \bm{f}_0(A)$.
\end{lemma}
The next lemma quantifies the error of Gaussian process interpolation.
\begin{lemma}[Theorem 11.4 of \cite{wendland2004scattered}]\label{lem:GP interp}
For all $f\in\mathbb{H}^{\tau,s}$, we have
$$\Big|f(X_i) - \mathbb{E}[Z^{\tau,s}_{X_i}|Z^{\tau,s}_{\mr{pa}(X_i)}=\bm{f}(\mr{pa}(X_i))]  \Big| \lesssim \sqrt{\var \big\{Z^{\tau,s}_{X_{i}}-\mathbb{E}[{Z}^{\tau,s}_{X_{i}}|{Z}^{\tau,s}_{\mr{pa}(X_{i})}]\big\}  \|f\|_{\mathbb{H}^{\tau,s}}^2}.$$
\end{lemma}
Then by combining the above two lemmas with the upper bound for the variance $\var \big\{Z^{\tau,s}_{X_{i}}-\mathbb{E}[{Z}^{\tau,s}_{X_{i}}|{Z}^{\tau,s}_{\mr{pa}(X_{i})}]\big\}$, see Theorem \ref{thm:var}, we obtain a more explicit control on the approximation error of the Gaussian process interpolation.
\begin{lemma}\label{lem:Holder appro}
Let $f_0\in C^\beta_1$ with $\beta\le \alpha$, and let $Z^{\tau,s}$ be the rescaled Mat\'{e}rn process with smoothness parameter $\alpha$. Under Conditions \ref{cond:layered} and \ref{cond:norming}, for all $X_i\in\X_n$ and its parent set $\mr{pa}(X_i)$, we have
$$\left|f_0(X_{i})-\mathbb{E}[Z^{\tau,s}_{X_{i}}|Z^{\tau,s}_{\mr{pa}(X_i)}=\bm{f}_0(\mr{pa}(X_i))] \right| \lesssim \big[\tau^\alpha\gamma^{-\eta(i)\alpha}  + \gamma^{-\eta(i)\beta} \big] (\eta(i)\ln \gamma)^{\mathbbm{1}_{\{\alpha\in\mathbb{N}\}}}.$$
\end{lemma}

Building on the technical results above, we can now provide the following upper bound for the decentering function.
\begin{lemma}[Decentering]\label{lem:decentering}
Let $f_0\in C^\beta_1$ with $\beta\le \alpha$ and let $\hat{Z}^{\tau,s}$ be the Vecchia approximation of the rescaled Mat\'{e}rn process with smoothness parameter $\alpha$.
Then, under Conditions \ref{cond:layered} - \ref{cond:norming}, 
\begin{align*}
\inf_{f\in\mathbb{H}_n^{\tau,s}:\|f-f_0\|_{\infty,n}\le \tilde{\vartheta}_n\epsilon} \|f\|_{\mathbb{H}_n^{\tau,s}}^2 
\lesssim & s^{-2}\Big[\tau^d \epsilon^{-d/\alpha} + \tau^{-2\alpha}  \epsilon^{-\frac{2(\alpha-\beta) +d}{\beta}}\Big] [\ln (1/\epsilon) ]^{\mathbbm{1}_{\{\alpha\in\mathbb{N}\}}}.
\end{align*}
\end{lemma}
Finally, by combining the upper bounds derived for the (centered) small ball probability and the decentering function of the Vecchia GP in Theorem \ref{lem:small ball} and Lemma \ref{lem:decentering}, respectively, we obtain an upper bound for the concentration function \eqref{def:conc:funct}. This, in turn, via \eqref{eq:conc:ball}, provides a lower bound for the decentered small ball probability.
\begin{lemma}\label{thm:small:decent}
Let $f_0\in C^\beta_1$ with $\beta\le \alpha$ and let $\hat{Z}^{\tau,s}$ be the Vecchia approximation of the rescaled Mat\'{e}rn process with smoothness parameter $\alpha$. Then, under Conditions \ref{cond:layered} - \ref{cond:norming}, 
\begin{align*}
-\ln\pr(\|\hat{Z}^{\tau,s}-f_0\|_{\infty,n} <\tilde\vartheta_n\epsilon)
\lesssim & s^{-2}\Big[\tau^d \epsilon^{-d/\alpha} + \tau^{-2\alpha}  \epsilon^{-\frac{2(\alpha-\beta) +d}{\beta}} \Big] [\ln (1/\epsilon) ]^{\mathbbm{1}_{\{\alpha\in\mathbb{N}\}}}\\
&\qquad\qquad+  \tau^d s^{d/\alpha}\epsilon^{-d/\alpha} [\ln (1/\epsilon)]^{\mathbbm{1}_{\{\alpha\in\mathbb{N}\}}}.
\end{align*}
\end{lemma}

\section{Bayesian nonparametrics}\label{sec:BNP}
In this section, we provide posterior contraction rate guarantees for the Vecchia GP in the context of the nonparametric regression model. We note, however, that given the general theory for Gaussian processes (see \cite{van2008rates}), these results can be readily extended to other nonparametric models, including density estimation and classification. Building on the probabilistic properties developed in the previous sections, we present two main theoretical results.  First, we derive contraction rates for the Vecchia GP $\hat{Z}^{\tau,s}$ with fixed scaling hyperparameters $\tau$ and $s$, showing that optimal, oracle choices of these parameters -- depending on the regularity parameter $\beta$ of the regression function $f_0$ -- achieve minimax optimal rates. In the second part of the section, we consider a hierarchical Bayesian framework, where the scaling hyperparameters $(\tau,s)$ are endowed with a hyper-prior distribution. We show that this fully Bayesian approach can adapt to the minimax rate without requiring any knowledge of the underlying true function. 

\subsection{Posterior contraction rates}\label{sec:contract}
There is a rich literature on posterior contraction rates for Gaussian processes based on the concentration function $\phi_{f_0,n}^{\tau,s}(\epsilon_n)$,  defined in \eqref{def:conc:funct}, see e.g. \cite{van2008rates,ghosal2017fundamentals}. In fact, by relating the metric used for testing the regression function $f_0$ to the norm of the Banach space in which the GP resides, it can be shown that the posterior contracts around the truth at the rate $\eps_n$  with respect to the testing metric, provided that $\phi_{f_0,n}^{\tau,s}(\epsilon_n)\lesssim n\epsilon_n^2$. We formalize this result in the context of the nonparametric regression model by taking the empirical $L_\infty$ norm as the Banach space norm and the empirical $L_2$  norm as the testing metric, which in turn determines the posterior contraction rate.
\begin{lemma}
\label{lem:GP contraction}
Assume that $\phi_{f_0,n}^{\tau,s}(\epsilon_n)\lesssim n\epsilon_n^2$ holds for some $\epsilon_n\to 0$ for all $f_0\in C^\beta_1$. Then, for every sequence $M_n\to\infty$, the posterior distribution of the regression function $f$ satisfies
$$\Pi(\|f-f_0\|_{2,n}>M_n\epsilon_n|\D_n)\stackrel{P_{f_0}^{n}}{\to} 0.$$
\end{lemma}

In the previous section we derived lower bounds for the decentered small ball probability in Lemma  \ref{thm:small:decent}, using the upper bound on the concentration function obtained by combining Theorem \ref{lem:small ball} and Lemma \ref{lem:decentering}. By plugging this upper bound into the concentration inequality, one can obtain the posterior contraction rate in view of Lemma \ref{lem:GP contraction}.
\begin{theorem}\label{thm:contraction}
Let $f_0\in C^\beta_1$, with $\beta\le \alpha$, and let $\hat{Z}^{\tau,s}$ be the Vecchia approximation of the rescaled Mat\'{e}rn process with smoothness parameter $\alpha$. Under Conditions \ref{cond:layered} - \ref{cond:norming}, 
$$\Pi(\|f-f_0\|_{2,n}>M_n\epsilon_n|\D_n)\stackrel{P^n_{f_0}}{\to} 0$$
 for all  sequences $M_n\to\infty$ and
\begin{align}\label{eq:epsilon choice}
\epsilon_n = 
\max\Big\{  &
(ns^2\tau^{-d})^{-\frac{\alpha}{2\alpha+d}}\tilde\vartheta_n^{\frac{d}{2\alpha+d}},
n^{-\frac{\beta}{2\alpha+d}} (\tau^\alpha s)^{-\frac{2\beta}{2\alpha+d}} \tilde\vartheta_n^{\frac{2(\alpha-\beta)+d}{2\alpha+d}}, 
\nonumber  \\ 
& n^{-\frac{\alpha}{2\alpha+d}} (\tau^\alpha s )^{\frac{d}{2\alpha+d}}\tilde\vartheta_n^{\frac{d}{2\alpha+d}}
\Big\} (\ln n)^{\mathbbm{1}_{\{\alpha\in\mathbb{N}\}}}.
\end{align}
\end{theorem}
Note, that the rate obtained in Theorem \ref{thm:contraction} depends crucially on the choice of the rescaling hyperparameters $\tau$ and $s$, as well as on the operator norm upper bound of the recursive GP interpolation $\tilde\vartheta_n$. We have shown that $\tilde\vartheta_n$ is uniformly bounded in case of grid data. More generally, assuming its boundedness allows us to omit it from the rate. Furthermore, if both $\tau$ and $s$ are constants, we retrieve the well-known posterior contraction rates derived in \cite{van2008rates}. This rate is minimax optimal when the prior's regularity matches the regularity of the underlying true function, i.e. $\alpha=\beta$. At the same time, by allowing $s$ and $\tau$ to depend on $n$ and minimizing the right-hand side of equation (\ref{eq:epsilon choice}) one can obtain the minimax rate $n^{-\beta/(2\beta+d)}$ even when the regularity of the prior does not match the regularity of $f_0$, see also \cite{ghosal2017fundamentals}. These results are summarized in the following corollary.

\begin{corollary}\label{cor:rates}
Suppose the conditions of Theorem \ref{thm:contraction} hold and  $\tilde\vartheta=O(1)$. If $\tau\asymp s\asymp 1$, the posterior contracts at the rate
$\epsilon_n =n^{-\frac{\beta}{2\alpha+d}} (\ln n)^{\mathbbm{1}_{\{\alpha\in\mathbb{N}\}}}$.
On the other hand, if $\tau^\alpha s = n^{\frac{\alpha-\beta}{2\beta+d}}$ with $s\gtrsim 1$, one achieves the minimax contraction rate $\epsilon_n =n^{-\beta/(2\beta+d)} (\ln n)^{\mathbbm{1}_{\{\alpha\in\mathbb{N}\}}}$.
\end{corollary}

Theorem \ref{thm:contraction} and Corollary \ref{cor:rates} show that the Vecchia approximations of Mat\'{e}rn processes achieve the same posterior contraction rates as their mother GP  counterpart, see e.g. \cite{szabo2023adaptation}. This confirms that one can enjoy the computational scalability of Vecchia Gaussian processes without any loss in estimation accuracy. In this way, our results address a long-standing gap in the statistical guarantees for Vecchia GPs that has persisted in recent years.

\subsection{Adaptation with hierarchical Bayes}\label{sec:adapt}
While Corollary \ref{cor:rates} recovers the optimal minimax contraction rate for arbitrary $\alpha\ge \beta$, this can be achieved only for certain oracle choices of the scaling hyperparameters $\tau, s$, depending on the smoothness $\beta$ of the true regression function $f_0$. In practice, however, this information is typically unavailable. A natural solution is to introduce another layer of prior on the scaling parameters $\tau$ and $s$, allowing them to adapt automatically to the dataset. Let us consider continuous priors on these hyperparameters and denote the associated probability density function  by $p(\tau,s)$. Then, the hierarchical Vecchia GP  $\hat{Z}_x^{\pr}$ takes the form
\begin{equation}\label{def:hierVecchiaGP}
\begin{split}
(\tau,s)&\sim p(\tau,s),\\
\hat{Z}_x^{\pr}|\tau,s &\stackrel{d}{=}\hat{Z}_x^{\tau,s},
\end{split}
\end{equation}
where  $\hat{Z}^{\tau,s}$  denotes the  Vecchia approximation of the rescaled Mat\'{e}rn process with smoothness parameter $\alpha$ and scale parameters $\tau$ and $s$. To formulate our contraction rate results, we require the following condition on the hyperprior for the scale parameters.
\begin{condition}\label{cond:hyperprior}
The hyperprior on $\tau,s$ satisfies the following  equations:
\begin{align*}
&\ln \pr\left(\tau^\alpha s > n^{\frac{\alpha-\beta}{2\beta+d}} \right) \lesssim - n^{\frac{d}{2\beta+d}},\quad \ln\pr \left( \tau < n^{-\frac{\beta}{2\beta+d}}\ln^2 n\right)\lesssim-n^{\frac{d}{2\beta+d}},\\
 &\ln\pr \left( s< 1/n \right) \lesssim - n^{\frac{d}{2\beta+d}},\quad\ln \pr \left( \big\{\tau^\alpha s \in \big[n^{\frac{\alpha-\beta}{2\beta+d}}, 2n^{\frac{\alpha-\beta}{2\beta+d}}\big] \big\} \cup \{s\ge 1\}\right) \gtrsim - n^{\frac{d}{2\beta+d}}.
\end{align*}
\end{condition}
The following theorem presents our result regarding adaptation over $\tau$ and $s$.
\begin{theorem}\label{thm:adaptation}
Let $f_0\in C^\beta_1$ with $\beta\le \alpha$ and let  $\hat{Z}_x^{\pr} $ be the hierarchical Vecchia GP defined in \eqref{def:hierVecchiaGP} with hyper-prior satisfying Condition \ref{cond:hyperprior}.
If, for all $\tau,s$ in the support of $p(\tau,s)$, the process $\hat{Z}^{\tau,s}$ satisfies Conditions \ref{cond:layered} - \ref{cond:norming} and the Gaussian interpolation bound $\tilde\vartheta_n=O(1)$, then, for any sequence $M_n\to\infty$, 
$$\Pi(\|f-f_0\|_{2,n}>Mn^{-\frac{\beta}{2\beta+d}}  (\ln n)^{\mathbbm{1}_{\{\alpha\in\mathbb{N}\}}}|\D_n)\stackrel{P_{f_0}^n}{\to} 0.$$
\end{theorem}

As a special case, one can fix one of the scaling hyperparameters and endow the other with a hyper-prior. This also yields a rate-adaptive procedure. The following corollary illustrates this when the time-scaling hyperparameter $\tau$ is random and the other, $s$, is fixed.
\begin{corollary}\label{cor:tau adapt}
Fix $s=1$ and endow  $\tau$ with a hyper-prior satisfying
$$\ln p(\tau) \asymp - n^{\frac{d}{2\alpha+d}} \tau^{\frac{2\alpha d}{2\alpha+d}},\quad \forall \tau \ge 1.$$
Then Condition \ref{cond:hyperprior} is satisfied and consequently the posterior corresponding to the hierarchical Vecchia GP prior achieves the minimax contraction rate under  Conditions \ref{cond:layered} - \ref{cond:norming} and $\tilde\vartheta_n=O(1)$.
\end{corollary}

\begin{remark}
We note that adaptive contraction rate results for  Mat\'{e}rn processes have been derived in several papers. However, in most of these works, the time-scaling parameter was fixed at $\tau=1$ and adaptation was considered only with respect to the space-scaling parameter $s$. The only exceptions are the recent papers \cite{szabo2023adaptation,fang:2025}. In the former one the regularity $\alpha$ was restricted to half integers, while in the later one the time-scale parameter $\tau$ was assumed to exceed 1. The major technical challenge lies in characterizing the RKHSs corresponding to different rescaling parameters $\tau$, which do not possess a simple inclusion structure. By leveraging the probabilistic properties developed in Section \ref{sec:prob}, we are able to prove adaptation for Vecchia approximations of Mat\'{e}rn processes in arbitrary dimensions and smoothness levels, simultaneously considering both space and time rescaling hyperparameters without the restrictions imposed in previous works.
\end{remark}

\section{Discussion}\label{sec:discussion}
\subsection{Extension to non-grid data}\label{sec:non grid}
For simplicity of presentation, our paper focuses on grid data. A natural question is whether these results can be extended to more general settings.

From a methodological perspective, the main challenge is constructing DAGs that continue to satisfy the Norming condition  \ref{cond:norming}. One approach is to directly construct a DAG on non-grid training data while maintaining a tight control on the norming constant. We discuss this approach in details in Section \ref{sec:build dag general} of the Supplementary Material. There, we propose an algorithm constructing Norming sets on general datasets, though it currently lacks theoretical underpinning. A second approach is to construct a grid reference set (which may not contain any actual data point) and build a layered norming DAG on it. The predictive process is then evaluated on the  non-grid training set. A rigorous analysis of this approach is beyond the scope of the present paper and is left for future research.

From a theoretical perspective, there are two main challenges in extending the probabilistic and statistical results to non-grid data. The first, as noted above, is controlling the norming constant. The second, as discussed in Section \ref{sec:rec:inter}, is controlling the norm of the recursive polynomial interpolation. Both of these present technically challenging analytical problems. Since our focus is on the probabilistic and statistical properties of the process, we do not investigate these extensions in further details.

\subsection{Statistical guarantees of approximation methods}
Vecchia approximations of Gaussian processes have gained popularity over the past decade due to their scalability to large datasets. However, the theoretical underpinning of Vecchia GPs has lagged behind their widespread application. Previous works have primarily focused on the approximation properties of Vecchia GPs relative to the mother GP, rather than studying them as standalone processes. Therefore, the limited literature on the statistical guarantees of Vecchia GPs typically focuses on cases where the approximation error -- measured in terms of the Wasserstein distance or KL divergence -- between a Vecchia GP and its corresponding mother GP is sufficiently small. While this approach is valid, it may lead to suboptimal results by requiring unnecessarily close proximity to the original GP, which is not needed for optimal statistical inference. In fact, Vecchia GPs can exhibit fundamentally different properties from their mother GPs. For instance, the popular Mat\'{e}rn process is stationary, while this property is generally not preserved in its Vecchia approximation. In this paper, we derive fundamental probabilistic (e.g. centered and decentered small ball probabilities) and statistical properties (e.g posterior contraction rates) questions that have not yet been studied in the literature

\subsection{Polynomials and Vecchia GPs} One of the fundamental contributions of our paper is the characterization of Vecchia GPs via polynomial interpolation. In particular, Theorem \ref{thm:var} and Lemma \ref{lem:GP flat} show that, under mild conditions, Mat\'{e}rn GPs, as well as their Vecchia approximations, can be well approximated by local polynomial interpolations when conditioned on a norming set. Moreover, the cardinality of the norming set matches the dimension of the corresponding polynomial vector space. For a Mat\'{e}rn Gaussian process with regularity $\alpha$, this vector space consists of polynomials of order at most $\underline{\alpha}$, which has cardinality ${\underline{\alpha}+d \choose \underline{\alpha}}$. This implies that, regardless of the sample size $n$, a finite, well-chosen parent set is sufficient, making Vecchia GPs highly scalable to large datasets.

While the choice of norming parent sets  is relatively straightforward in one dimension, it becomes significantly more challenging in higher-dimensional spaces. As illustrated in Figure \ref{fig:norming d=2}, sets with different geometric configurations can have widely varying norming constants. It is even possible for a collection of nearby points to fail to be a norming set. This partially explains why remote locations were proposed as parent sets in the literature \citep{banerjee2015hierarchical}. In general, finding norming sets and controlling the associated norming constants is a difficult problem. Most often, norming constants are determined numerically, with only a few analytic results available, for example the Fekete points studied by \cite{bos2018fekete}.

\subsection{Further extensions} We focused on isotropic  Mat\'ern processes in the present paper. To extend our results to non-isotropic Mat\'ern processes -- where different coordinates may have different regularities or scale parameters -- one may consider interpolation based on tensor product basis. In addition, since the Mat\'ern kernel converges to the squared exponential as $\alpha$ tends to infinity, a natural direction for future work is to investigate whether our probabilistic and statistical results extend to Gaussian kernels as well. 

Recall that, outside of the reference set, the Vecchia GP is defined through conditional independence given the values of the process at the parent set. When the reference set coincides with the training data and lies on a regular grid, this does not lead to excessive volatility. However, if the data points are irregularly spaced, and particularly if they are sparse in some regions of the domain, the resulting process can become highly irregular. As discussed above, one way to avoid this issue is to decouple the reference set from the training data and instead use a fixed grid. Another option is to employ deterministic interpolation, e.g. polynomial or Gaussian, between reference points to obtain a more regular approximation. 

Finally, we note that our focus in this paper has been on establishing contraction rates. An equally important aspect of Bayesian procedures is the reliability of their uncertainty quantification, i.e. whether credible sets can be interpreted as valid confidence sets. While this question is well studied for the full GP posterior \cite{szabo:etal:2015,rousseau2020asymptotic}, considerably less is known about its Vecchia approximation. Our numerical analysis offers some preliminary insight and suggests that the approximation may provide reasonable uncertainty quantification in practice. However, a rigorous theoretical investigation of this issue lies beyond the scope of the present paper. 


\begin{acks}[Acknowledgments]
The authors would like to thank Isma\"{e}l Castillo, Yuansi Chen, Zichao Dong, Cheng Li, Andy McCormack, Galen Reeves, Bernhard Stankewitz and Dun Tang for helpful discussions, Michele Peruzzi for computational resources, as well as the Associate Editor and referees for their constructive feedback and suggestions. 
\end{acks}

\begin{funding}
Co-funded by the European Union (ERC, BigBayesUQ, project number: 101041064). Views and opinions expressed are however those of the author(s) only and do not necessarily reflect those of the European Union or the European Research Council. Neither the European Union nor the granting authority can be held responsible for them.
\end{funding}

\begin{supplement}
\stitle{Supplementary Material to ``Vecchia Gaussian Processes: Probabilistic  and statistical properties''}
\sdescription{The Supplementary Material contains numerical studies for the proposed method and background on polynomials. Furthermore, the proofs of the theorems, lemmas, and corollaries from the main paper, along with additional technical lemmas, are also included.}
\end{supplement}
\begin{supplement}
\stitle{Package ``GPDAG'': \url{https://github.com/yi-chen-zhu/GPDAG}}
\sdescription{R package for Bayesian posterior inference of Gaussian processes with general DAG structures, including the layered norming DAG proposed in this paper. The repository also contains the codes used for the numerical studies.}
\end{supplement}


\bibliographystyle{imsart-number} 
\bibliography{PVP.bib}       


\newpage
\begin{appendix}
\counterwithin{equation}{section}
\counterwithin{figure}{section}
\counterwithin{table}{section}

\counterwithin{theorem}{section}
\counterwithin{lemma}{section}
\counterwithin{corollary}{section}
\counterwithin{condition}{section}
\counterwithin{remark}{section}

\title{Supplementary Material to ``Vecchia Gaussian processes: on probabilistic and statistical properties''}
\begin{aug}
\author[A]{\fnms{Botond}~\snm{Szabo}\ead[label=e1]{botond.szabo@unibocconi.it}}
\and
\author[B]{\fnms{Yichen}~\snm{Zhu}\ead[label=e2]{yczhu@hku.hk}}
\address[A]{Bocconi Institute for Data Science and Analytics,
Bocconi University \printead[presep={ ,\ }]{e1}}
\address[B]{Department of Statistics and Actuarial Science,
the University of Hong Kong \printead[presep={,\ }]{e2}}
\end{aug}

\vspace{5mm}
The Supplementary Material contains the numerical analysis of our method for estimation and uncertainty quantification in Section \ref{sec:num}. We then provide background material on polynomial interpolation in Section \ref{sec:poly2DAG}. Section \ref{sec:algs} presents the algorithms for building layered norming DAGs. 
We discuss the posterior inference algorithm and its computational complexity in Section \ref{sec:post inf}.
Finally, sections \ref{sec:proof prob 1} - \ref{sec:proof BNP} provide the proofs of the theorems, lemmas, and corollaries of the main paper, organized in the same order in which they appear there. Auxiliary lemmas are stated and proved when necessary.
In all proofs from Section \ref{sec:alpha int |} onward, we consider only the case $\alpha\not\in\mathbb{N}$. The proofs for $\alpha\in\mathbb{N}$ follow analogously, with the inclusion of a logarithm term in the relevant formulas. 

\section{Numerical analysis}\label{sec:num}
The practical success of Vecchia Gaussian processes has been widely demonstrated in the literature. In our numerical analysis, we focus on validating the key mechanisms of Vecchia GPs that are conveyed throughout our theorems and proofs.
Since Vecchia GP methods differ only in their DAG structures, we provide a universal package, \texttt{GPDAG}, for posterior inference of Vecchia GPs to ensure a fair comparison. This package is implemented in C++ using the state-of-the-art preconditioned conjugate gradient descent method \cite{kundig2023iterative} and includes a high-level R interface. A detailed description of the posterior inference algorithms is provided in Section \ref{sec:post inf}.

\subsection{Optimal estimation without prior approximation}\label{sec:est}
In the literature, Vecchia GPs are often viewed as approximations of their corresponding mother GPs. As a result, the quality of a Vecchia GP is typically assessed by how accurately it approximates the mother GP under an appropriate metric. In this section, we numerically demonstrate one of our key arguments: Vecchia GPs can perform optimal nonparametric estimation even when they do not provide a good approximation to the corresponding mother GP. Specifically, we consider the nonparametric regression problem and investigate the performance of the posterior resulting from the Vecchia GP prior (\ref{eq:prior vecchia}). 

We implement Vecchia GPs using two different DAG structures. The first is the layered norming DAG (denoted as \textit{Norming} in the figures) described in Section \ref{sec:layered norming dag}, see also Algorithm \ref{alg:dag grid} for formalization. The second is the Nearest Neighbor Gaussian Process (NNGP) of \citep{datta2016hierarchical} with maximin ordering \citep{guinness2018permutation} (denoted as \textit{Maximin} in the figures). For simplicity, we take the interval $[0,1]$ as the domain and use an equidistant design on it as the training data $\X_n$, where $n$ ranges from $17$ to $2049$ (i.e., $n=2^k+1$, with $k$ ranging from 4 to 11). The true regression function $f_0$ is a randomly generated function from the sample paths of a Mat\'{e}rn process with smoothness $\alpha=3/2$ and scaling parameters $\tau=10, s=1$. We set the error standard deviation in equation (\ref{eq:regression})  to $\sigma=0.1$.
For both Vecchia GP methods, we use the same Mat\'{e}rn process that generated the true $f_0$. In other words, the prior smoothness matches the truth exactly and the Vecchia GP with the layered norming DAG should therefore achieve the minimax rate according to our theory. 

The DAG structures of the two Vecchia GP methods considered are as follows. For the layered norming DAG in Algorithm \ref{alg:dag grid}, since $\underline{\alpha}=1$ and $d=1$, the norming parent sets have cardinality $2$. This means for $X_i\in\N_j$, the parent set $\mr{pa}(X_i)$ consists of the two locations in $\cup_{j'=0}^{j-1}\N_{j'}$ that are closest to $X_i$. For the NNGP with maximin ordering, we choose the neighborhood size to be the integer closest to $2\ln(n)$.  While there are ongoing debate regarding the optimal size of parent sets in NNGP, several previous works propose a logarithmic rate, see for instance \cite{kang2024asymptotic,schafer2021compression}. 

\begin{figure}[h!]
    \centering
    \includegraphics[width=\linewidth]{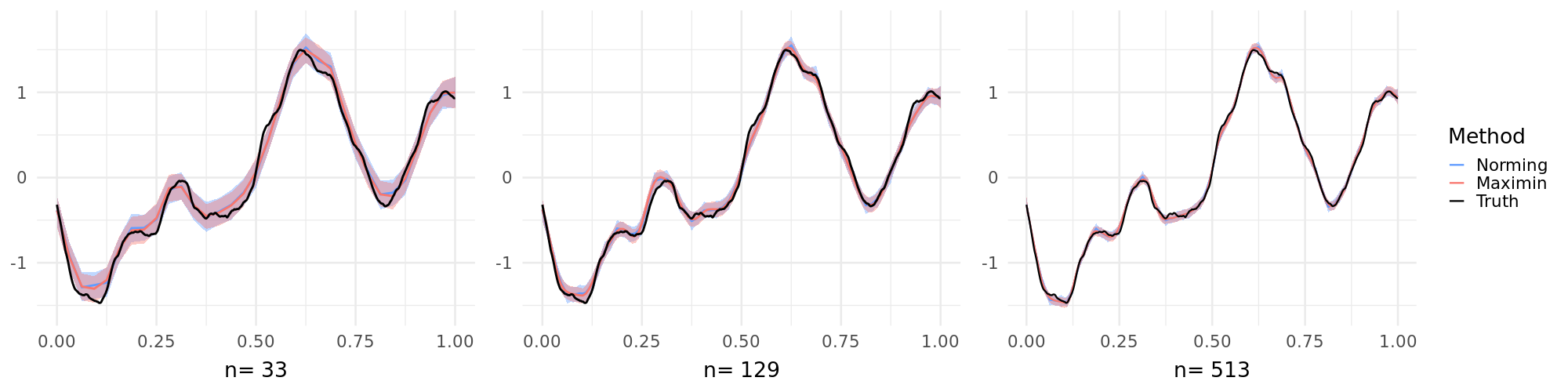}
    \caption{Illustration of nonparametric regression with Vecchia GPs. The black lines and dots represent the true regression function and the observed data, respectively. The colored lines and shaded regions depict the posterior means and the $95\%$ pointwise credible intervals obtained from the two Vecchia GP methods.}
    \label{fig:vecchia illu}
\end{figure}

\begin{figure}[h!]
    \centering
    \includegraphics[width=\linewidth]{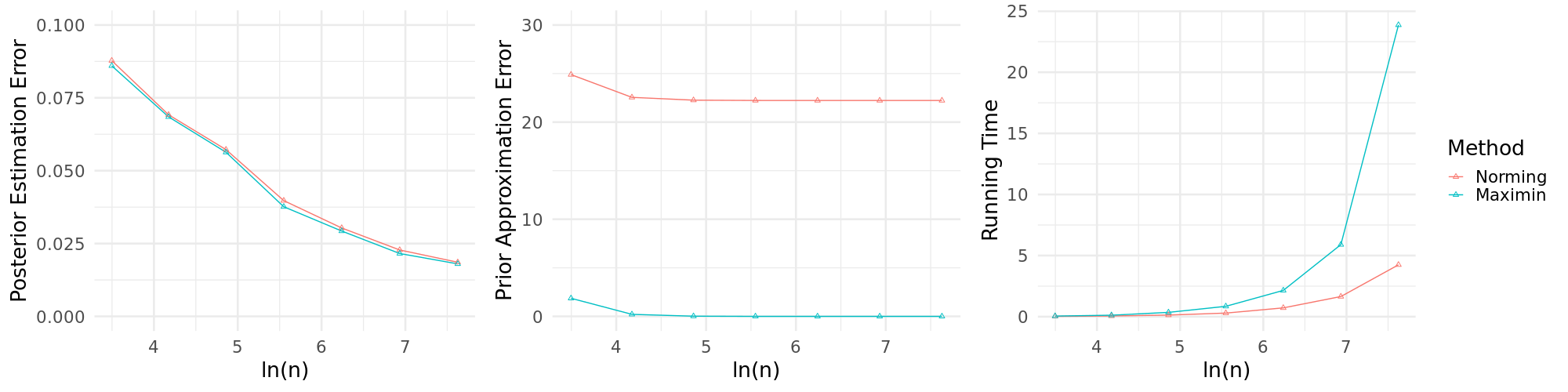}
    \caption{Qualitative results of the two Vecchia GP methods (Norming and Maximin) when the true function has H\"{o}lder regularity $\beta=1.5$. Left: posterior estimation error, measured by the $L_2$-distance between the true regression function $f_0$ and the posterior mean. Middle: prior approximation error, measured by the squared Wasserstein distance between marginals of the Vecchia GPs and their mother GPs. Right: run time of MCMC inference, measured in seconds.}
    \label{fig:est v appro}
\end{figure}

We illustrate the posterior contraction rate in the nonparametric regression model using Vecchia GPs in Figure \ref{fig:vecchia illu}.  Both Vecchia GP methods exhibit similar performances: when the sample size is small ($n=33$), they produce relatively large estimation errors and wide credible bands, but as the sample size grows ($n=513$), both methods recover the true function accurately. Figure \ref{fig:est v appro} provides a more detailed report of the two Vecchia GPs in terms of estimation accuracy, prior approximation, and run time. The estimation error is measured by the $L_2$-distance between the posterior mean and the true regression function. As the sample size $n$ increases, both Vecchia GP methods estimate the true regression function similarly well in $L_2$-norm. The prior approximation error is measured by the squared Wasserstein-2 distance between the marginal distributions of the Vecchia GP and its mother GP on $\X_n$. As shown in the figure, the Vecchia GP with layered norming DAG does not approximate the mother GP well, with its squared Wasserstein-2 distance staying above $20$ for any sample size $n$. In contrast, for the Vecchia GP with maximin ordering the squared Wasserstein-2 distance tends to zero. Taken together, the posterior estimation and prior approximation results demonstrate that Vecchia GPs do not require close approximation of their corresponding mother GPs. Instead, it is sufficient for them to inherit key regularity properties in order to recover the true function.

We also plot the computation time of the two Vecchia GPs on the right side of Figure \ref{fig:est v appro}. Since the Vecchia GP with layered norming DAG has fixed cardinality of $2$ for all parent sets, whereas the NNGP with Maximin ordering uses parent sets of size $2\ln (n)$, it is intuitively clear that the former can be computed much faster than the latter. The difference in runtime increases with the sample size. Therefore, if the objective is to recover the underlying true function, one can employ Vecchia GPs with very small parent set sizes (e.g. $2$ in this example), gaining substantial computational benefits without sacrificing statistical efficiency.

Finally, we assess the quality of the uncertainty quantification for both Norming DAG and NNGP with Maximin ordering. For each training point $X_i$, we check whether its $95\%$ posterior credible interval contains the true function value $f_0(X_i)$, returning $1$ (yes) or $0$ (no). We then average these $0-1$ values over all training locations $X_i, 1\le i\le n$ to obtain the ``marginal posterior coverage'', reported in the left side of Table \ref{tab:coverage}. For both methods, the coverage becomes reasonable close to $0.95$ when sample size exceeds $200$. This suggests that the marginal Bayesian credible intervals provide valid frequentist coverage, at least asymptotically.

\subsection{Infinite smooth truth and finite norming sets}
According to our theoretical results, choosing parent sets as norming sets with cardinality $\underline{\alpha}+d \choose \underline{\alpha}$ leads to optimal nonparametric estimation rates. While such parent sets are straightforward to implement for Mat\'{e}rn processes with fixed regularity $\alpha$, they become practically infeasible for Gaussian processes with a squared exponential covariance kernel, which can be viewed as the limit of a Matérn process as $\alpha\to\infty$. This raises the natural question: how should one choose the parent sets for Gaussian processes with very high or even analytic regularity?

A rigorous theoretical investigation of this question is beyond the scope of the present paper. Nevertheless, we conjecture that if both the truth and the prior are $\alpha$-smooth for some large value $\alpha\in (0, +\infty]$, but we select parent sets as norming sets with cardinality $\tilde{\alpha}+d \choose \tilde{\alpha}$ for some integer $\tilde{\alpha} < \underline{\alpha}$, then the posterior will contract at a slower rate, depending on $\tilde{\alpha}$.

We demonstrate our conjecture through a numerical analysis. We take a sample path from a unidimensional Gaussian process with a squared-exponential covariance kernel as the true function. Note that this provides an analytic truth (with probability one). We then compare three Vecchia GP methods, all of which use the squared exponential covariance kernel as the mother GP. From a Bayesian perspective, this means that the prior belief regarding the regularity of the truth is correct. As in to Section \ref{sec:est}, these Vecchia GP methods differ only in their DAG structures. The method denoted as \textit{Maximin} is the Nearest Neighbor Gaussian Process with maximin ordering, as in the previous subsection, with parent sets of size $2\ln(n)$. The method denoted as \textit{Norming2} and \textit{Norming4} employ the layered norming DAG with norming sets of cardinality 2 and 4, respectively. 

\begin{figure}    
\centering    \includegraphics[width=\linewidth]{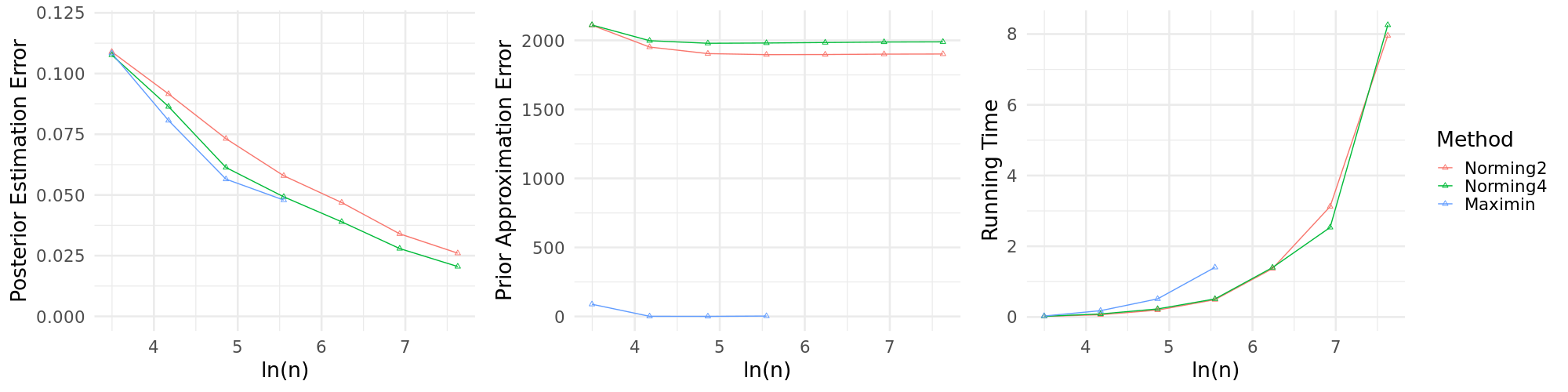}
    \caption{
    Qualitative results of three Vecchia GP methods for analytic truth. Left: posterior estimation error, measured by $L_2$-norm between the truth and the posterior mean. Middle: prior approximation error, measured by the squared Wasserstein distance between the marginals of the Vecchia GPs and their mother GPs. Right: run time of MCMC inference, measured in seconds.}
    \label{fig:inf smooth}
\end{figure}

The results of the numerical studies are shown in Figure \ref{fig:inf smooth}. Since, for the Maximin method, the covariance matrix of the parent set becomes nearly ``numerically singular'', even using double-precision floating-numbers in C++, we were only able to conduct studies up to sample size of $n=257$. As expected, the Norming2 method provides a less accurate estimation of $f_0$ compared to Maximin. Interestingly, however, Norming4 achieves similar accuracy as Maximin. This is due to two factors. First, for $\alpha=\beta=4$ and $d=1$, the optimal posterior contraction rate is $n^{-\beta/(2\beta+d)}=n^{-4/9}$, which is very close to the minimax rate for estimating analytic functions $n^{-1/2}$. As a result, the errors are visually hard to distinguish. In this case, $\underline{\alpha}=3$, so the size of the norming set is ${\underline{\alpha}+d \choose \underline{\alpha}}=4$. Therefore, a Vecchia approximation with 4 parents results in a qualitatively similar estimation to the original squared exponential covariance kernel. Second, as the size of the parent sets increases, the covariance matrix of the Vecchia GP increasingly resembles the covariance matrix of the mother GP. Consequently, the smallest eigenvalue of the Vecchia covariance matrix decays rapidly towards zero, causing numerical errors to dominate the statistical error.

Apart from the posterior estimation error, the results for prior approximation and run time behave as expected. The layered norming DAGs are not designed to approximate the mother process, hence it is not surprising that the prior approximation errors for Norming2 and Norming4 are large compared to those for the Maximin method. However, they require substantially less computation time and can be applied to larger sample sizes compared to the NNGP with maximin ordering.

Similar to Section \ref{sec:est}, we also report the marginal posterior coverage for all three methods, shown in the middle part of Table \ref{tab:coverage}. All coverage values are reasonably close to $0.95$, indicating that these Bayesian uncertainty quantifications are valid from a frequentist perspective as well.

\subsection{Estimating rough truth}\label{sec:rough}
Our theoretical analysis implies that for $\beta< 1$, i.e. when the true function is very rough, not even differentiable, a Vecchia GP with parent sets of size one can achieve minimax rates. At first glance, this may seem counter-intuitive.

To empirically validate this theoretical finding, we consider a H\"{o}lder-$1/2$ smooth true function on the unit interval $[0,1]$. We compare the behavior of three Vecchia GP methods, all using a Mat\'{e}rn process with  $\alpha=1/2$ as mother GP, differing only in their DAG structure. The Nearest Neighbor Gaussian Process with maximin ordering and parent sets of cardinality $2\log(n)$, serves as our benchmark procedure. Additionally, we consider the proposed layered norming DAG defined from Section \ref{sec:layered norming dag} with parent set cardinalities 1 and 2, denoted  ``Norming1'' and ``Norming2'', respectively.

\begin{figure}
    \centering
    \includegraphics[width=\linewidth]{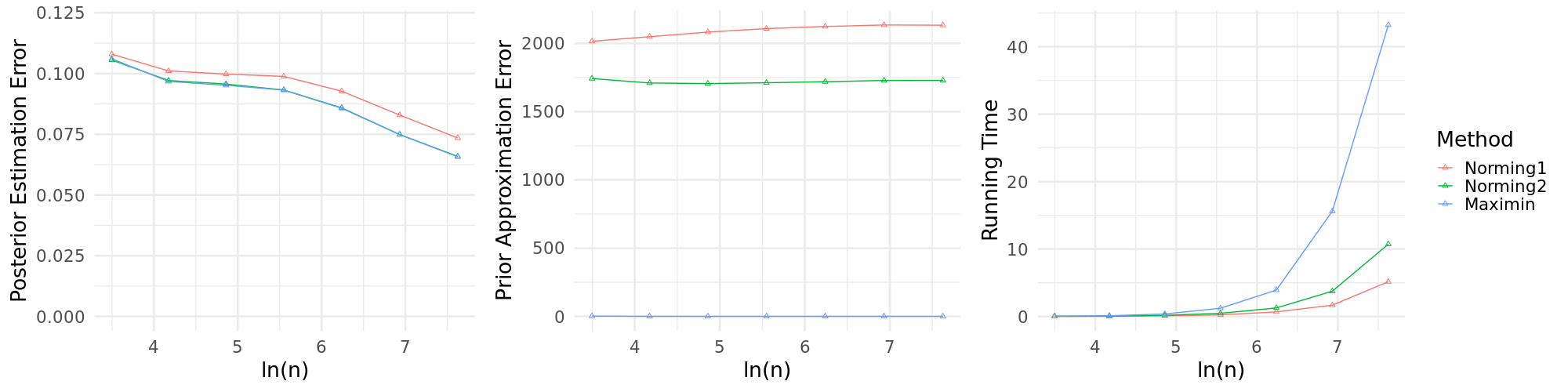}
    \caption{Qualitative results for three Vecchia GP methods for H\"{o}lder-$1/2$ truth. Left: posterior estimation error, measured by $L_2$-norm between the truth and the posterior mean. Middle: prior approximation error, measured by squared Wasserstein-2 (i.e. $W_2$) distance between the marginals of the Vecchia GPs and their mother GPs. Right: run time of MCMC inference, measured by seconds.}
    \label{fig:1d beta0.5}
\end{figure}

The numerical results are shown in Figure \ref{fig:1d beta0.5}. Compared to the other methods, ``Norming1'' with single-parent sets  performs slightly in terms of estimation accuracy, but the overall performance remains comparable and sufficiently accurate. We note that the asymptotic theory provides guarantees only for larger sample sizes, focusing on rates and ignoring constant multipliers. Hence, the application of larger parent sets, especially for small sample sizes, can improve empirical performance.
Therefore, from a practical perspective, we suggest setting the parent set size to at least $d+1$, ensuring that the resulting DAG is a connected directed graph. This corresponds to choosing $\underline{\alpha}\geq 1$. Rougher functions can still be well estimated with the inclusion of scale parameters $\tau$ and $s$.

As in previous sections, we report the marginal posterior coverage in the right portion of Table \ref{tab:coverage}. All coverage values are reasonably close to $0.95$, indicating that on average the marginal credible intervals are reliable from a frequentist perspective.

\begin{table}[h]
    \centering
    \begin{tabular}{c|cc|ccc|ccc}
    \hline
     & \multicolumn{2}{c|}{A.1} & \multicolumn{3}{c|}{A.2} & \multicolumn{3}{c}{A.3} \\ \hline
$n$ & Norming2 & Maximin & Norming2 & Norming4 & Maximin & Norming1 & Norming2 & Maximin \\ \hline
17 & 1.000 & 1.000 & 0.941 & 0.941 & 0.941  & 0.941 & 0.941 & 0.941 \\
33 & 0.939 & 0.909 & 0.939 & 0.970 & 0.970 & 0.970 & 0.939 & 0.970 \\
65 & 0.923 & 0.938 & 0.969 & 0.969 & 0.969 & 0.969 & 0.985 & 0.985 \\
129 & 0.876 & 0.876 & 0.961 & 0.984 & 0.977 & 0.984 & 0.969 & 0.977 \\
257 & 0.914 & 0.942 & 0.957 & 0.977 & 0.965 & 0.977 & 0.946 & 0.949 \\
513 & 0.918 & 0.924 & 0.906 & 0.934 & NA & 0.965 & 0.963 & 0.961 \\
1025 & 0.958 & 0.964 & 0.950 & 0.941 & NA & 0.957 & 0.950 & 0.948 \\
2049 & 0.928 & 0.939 & 0.949 & 0.956 & NA & 0.956 & 0.949 & 0.953 \\ \hline 
\end{tabular}
\vspace{1mm}
\caption{Marginal posterior coverage for $95\%$ Bayesian credible intervals, for all methods in Sections \ref{sec:est}-\ref{sec:rough}.}\label{tab:coverage}
\end{table}

\subsection{Two dimensional example}
In this subsection, we conduct numerical studies on the unit square domain $[0,1]^2$, and set the smoothness parameter $\beta$ of the truth to $0.5, 1.5, 2.5$, respectively. For each value of $\beta$, we compare two Vecchia GPs, both with Mat\'{e}rn mother GPs with smoothness $\alpha=\beta$. 
In the layered norming DAGs, we use corner parent sets of cardinality $d+\underline{\alpha} \choose \underline{\alpha}$, as suggested in Lemma \ref{lem:poly interp}. For Mat\'{e}rn process with $\alpha<1$, to avoid independent priors as described in Section \ref{sec:rough}, we use corner sets of cardinality ${d+1 \choose 1} = 3 $ as the parent sets. 

\begin{figure}[htbp]
  \centering

  \begin{subfigure}[b]{0.9\textwidth}
    \centering
    \includegraphics[width=\textwidth]{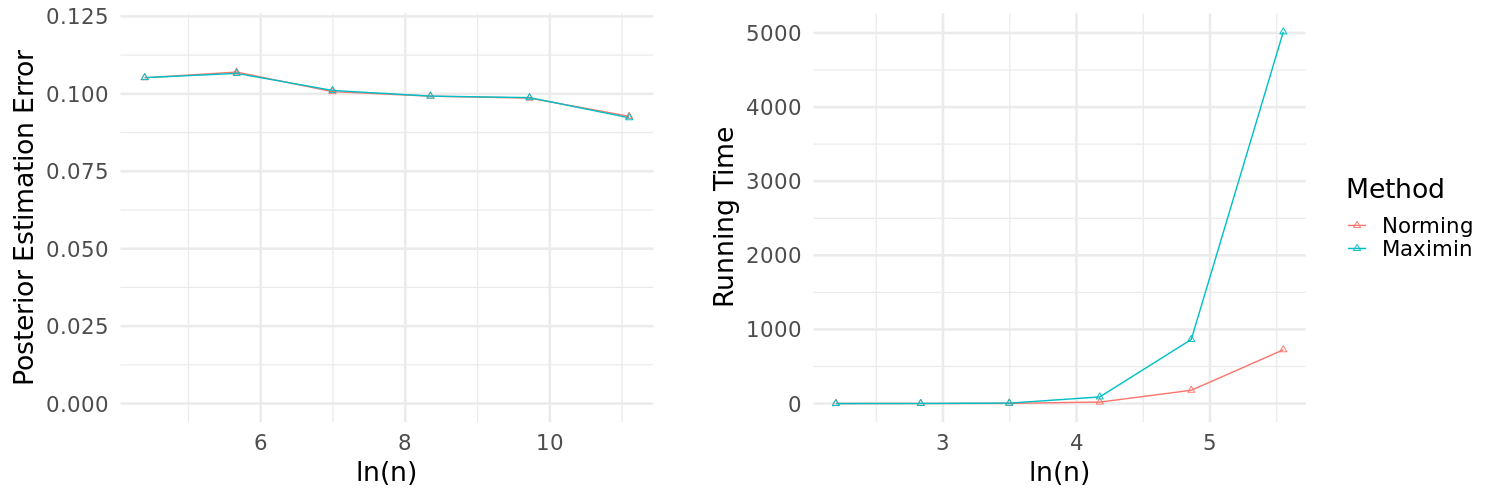}
    \caption{$\beta=0.5$}
  \end{subfigure}
  \hfill
  \begin{subfigure}[b]{0.9\textwidth}
    \centering
    \includegraphics[width=\textwidth]{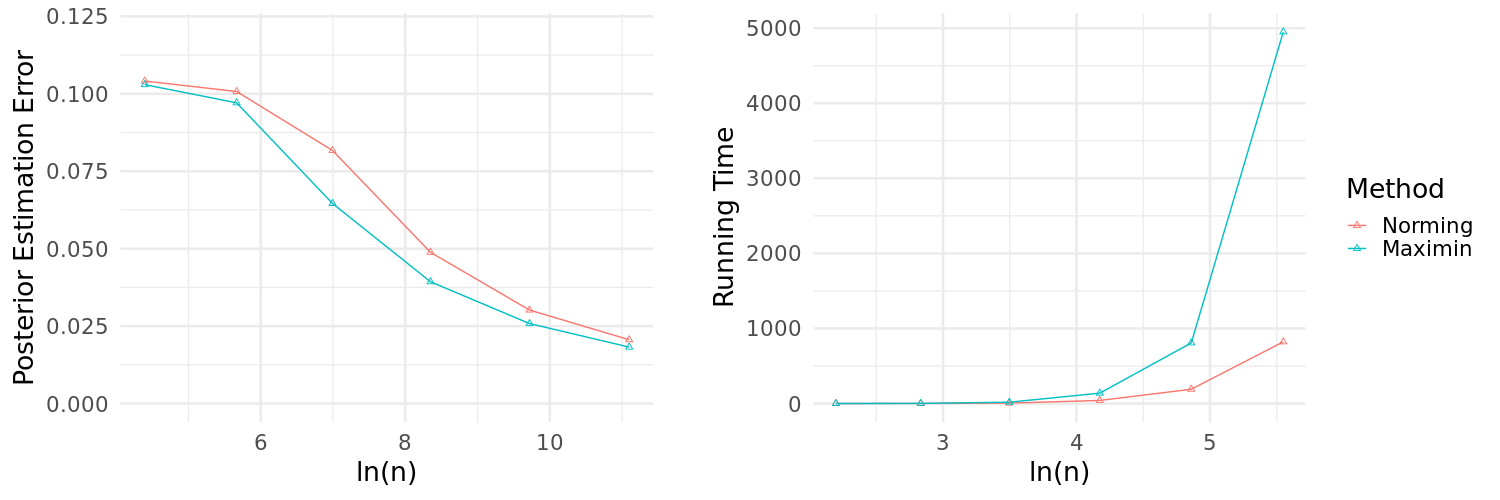}
    \caption{$\beta=1.5$}
  \end{subfigure}
  \hfill
  \begin{subfigure}[b]{0.9\textwidth}
    \centering
    \includegraphics[width=\textwidth]{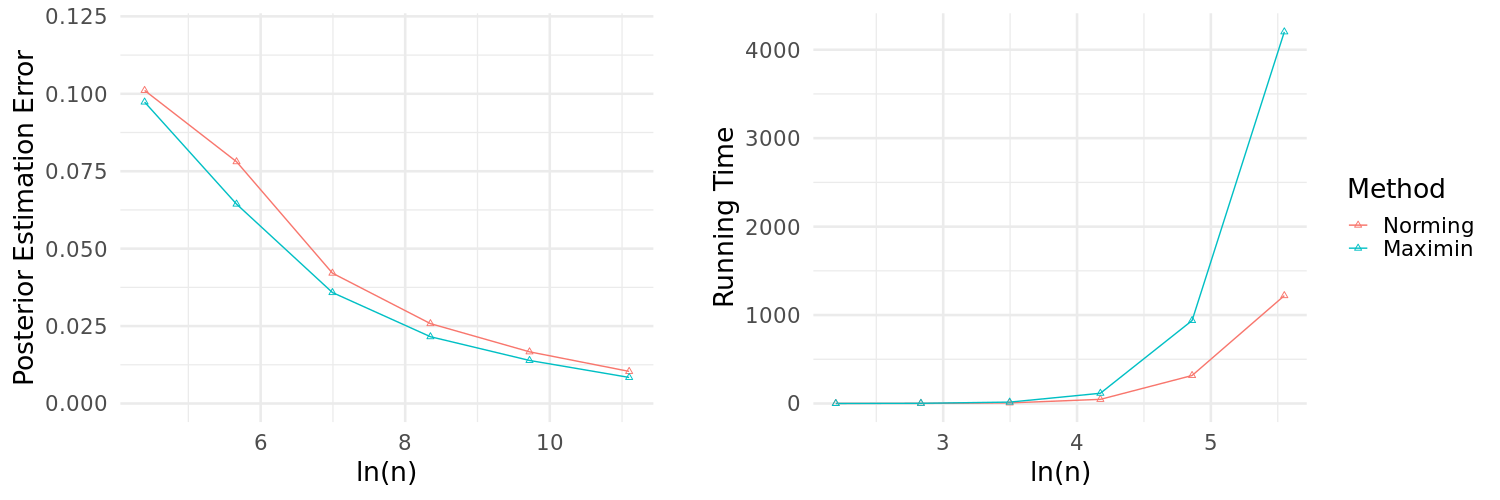}
    \caption{$\beta=2.5$}
  \end{subfigure}

  \caption{Vecchia GPs in two dimensions. Subfigures (a), (b) and (c) display the results when the truth has H\"{o}lder regularity $0.5$, $1.5$ and $2.5$, respectively. In each subfigure, the left panel shows the posterior estimation error, measured by $L_2$ norm between the truth and the posterior mean, while the right panel reports the MCMC run time in seconds.}
  \label{fig:2d}
\end{figure}

Figure \ref{fig:2d} presents the results on posterior estimation error and computation time. We omit the prior approximation accuracy, since it has already been demonstrated in the univariate setting that Vecchia GPs with layered norming DAGs are not intended for this purpose. Moreover, evaluating the $W_2$ distance in two dimensions is substantially more computationally demanding. 

We highlight two phenomena in the figure that support our methodology and theoretical findings. First, the two methods perform very similarly in terms of estimation accuracy, whereas Vecchia GPs with layered norming DAGs require only a fraction of the computational time needed for the NNGP with maximin ordering. Second, as the smoothness of the truth and the prior increases, the posterior estimation error decreases significantly. Note that for $\beta=0.5$, the $L_2$ posterior contraction rate is only $n^{-1/6}$, leading to very slow convergence. For $\beta=2.5$, the contraction rate improves to $n^{-5/14}$, and this substantially faster rate is clearly reflected in the numerical results.

\subsection{Adaptation with hierarchical Bayes}
In this subsection, we examine the adaptation phenomena discussed in Section \ref{sec:adapt} of the main paper. We consider a H\"{o}lder-$\beta$ smooth true function $f_0$ with $\beta=1.5$. In practice, the regularity parameter $\beta$ is typically unknown. Therefore, we employ a Vecchia approximation to a Mat\'{e}rn process with regularity parameter $\alpha$, chosen to be sufficiently large, and endow the scale parameter $\tau$ with a hyper-prior. In this example, we set $\alpha=2.5$ and use the following hyperprior on $\tau$
\begin{equation}\label{eq:tau prior exp}
\log p(\tau) = -\frac{1}{2} n^{\frac{d}{2\alpha+d}} \tau^{\frac{2\alpha d}{2\alpha+d}}+c,
\end{equation}
where $c$ is the normalizing constant.

We compare the posterior estimation error and runtime of four methods. (1) Baseline Norming: the Norming method with the oracle regularity $\alpha=\beta=1.5$ and scale parameter $\tau=1$, identical to the Norming method in Figure \ref{fig:est v appro}. (2) Oversmoothed Norming: the Norming method with $\alpha=2.5$ and $\tau=1$. (3) Hierarchical (data-driven) procedure: the method with $\alpha=2.5$ and $\tau$ equipped with the hyper-prior \eqref{eq:tau prior exp}. In our experiments, the rescaling parameter $\tau$ stabilizes after a few hundred MCMC iterations, therefore, we sample $\tau$ via Gibbs steps only once every ten MCMC iterations. (4) Baseline NNGP (Maximin): the NNGP method with maximin ordering and the oracle choice of hyperparameters $\alpha=\beta=1.5$ and $\tau=1$, corresponding to the Maximin method in Figure \ref{fig:est v appro}.

\begin{figure}
    \centering
    \includegraphics[width=0.8\linewidth]{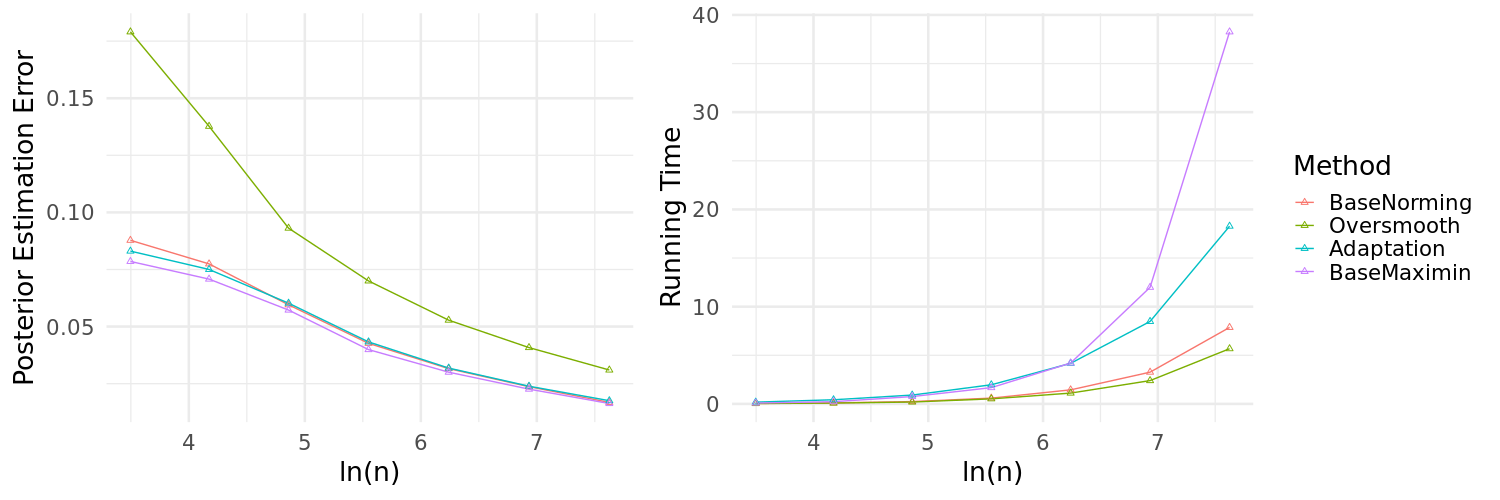}
    \caption{Illustration of Adaptation of Vecchia GP using hierarchical Bayes. The truth is H\"{o}lder 1.5 smooth. ``BaseNorming'' stands for layered norming DAGs with $\alpha=1.5$ prior, ``Oversmooth''  stands for layered norming DAGs with $\alpha=2.5$ prior, ``Adaptation''  stands for layered norming DAGs with $\alpha=2.5$ and hyper prior on rescaling parameter $\tau$, whereas ``BaseMaximin''  stands for NNGP with Maximin ordering and $\alpha=1.5$}
    \label{fig:adaptation}
\end{figure}

The results for posterior estimation error and runtime are presented in Figure \ref{fig:adaptation}. As expected, the behavior of two baseline methods are very similar in terms of posterior estimation error. In contrast, the oversmoothed method has poor performance, its contraction rate $n^{-\frac{\beta}{2\alpha+d}}\gg n^{-\frac{\beta}{2\beta+d}}$ is highly suboptimal. The hierarchical procedure with $\tau$ endowed with a hyper-prior, however, achieve comparable posterior estimation error to the two baseline oracle methods.

In terms of runtime, adaptation does result in additional computation costs. This is expected, as hierarchical Bayesian inference requires sampling extra variables. Specifically, each time a new value of $\tau$ is sampled, a new sparse Cholesky decomposition of the precision matrix must be computed. In summary, adaptation to unknown smoothness is entirely feasible in practice, but it comes with increased computational cost (although it still remains faster than the NNGP with a $2\log n$ parent size).

\subsection{Norming Sets}
\begin{figure}[h!]
    \centering
    \includegraphics[width=\linewidth]{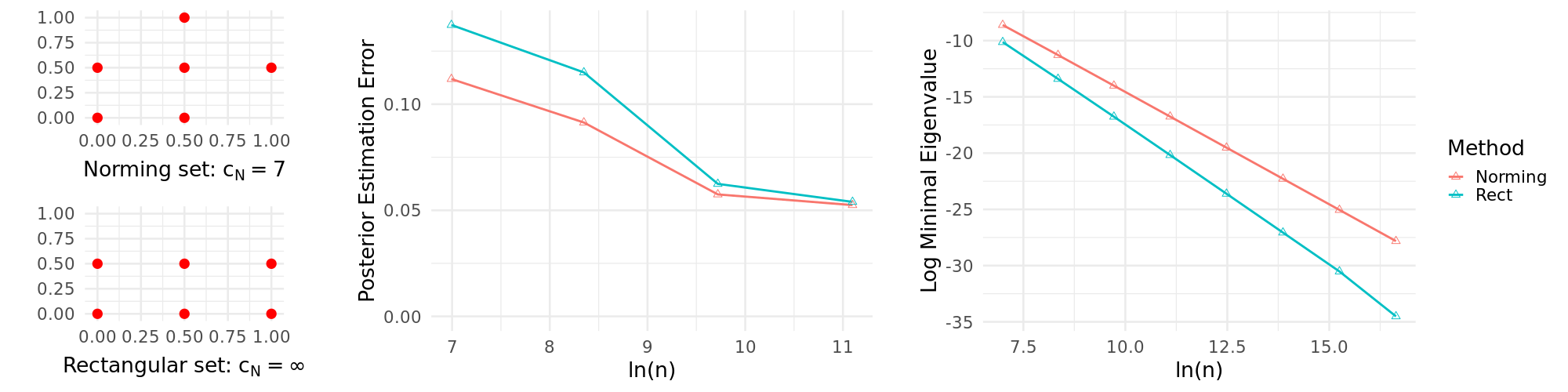}
    \caption{Comparison between the layered norming DAGs and Rectangle (non-norming) DAGs across different sample sizes. Left: geometric configurations of the norming and rectangular (non-norming) sets. Middle: $L_2$-distance between the truth and the posterior mean. Tight: infimum over $X_i$ of the logarithm of the minimal eigenvalues of the matrix $K_{\mr{pa}(X_i),\mr{pa}(X_i)}$ among all $X_i\in\X_n$.}
    \label{fig:norming}
\end{figure}

In view of our theoretical results, norming parent sets are sufficient to guarantee optimal estimation rates. In this section, we investigate the necessity of choosing norming sets. We keep all other aspects of the model fixed and compare DAGs with norming sets and non-norming sets as parents. We demonstrate numerically that violating the norming sets condition negatively affects the empirical performance of Vecchia GP methods. Since in a one-dimensional space, all sets with distinct elements are norming sets, we consider the unit square $[0,1]^2$ as the domain and set $\X_n$ to be a grid on it. We set the true regression function generating the data as $f_0(x) = 2(x[1]-0.5)^2 + 4(x[2]-0.5)^2 + x[1]x[2]^2.$
The mother Gaussian process is a Mat\'{e}rn process with smoothness parameter $\alpha=5/2$ and scaling parameters $\tau=8,s=1$. We consider Vecchia GPs with two different DAG structures in the comparison. The first is the layered norming DAG with $\underline{\alpha}=2$ and norming sets of cardinality ${\underline{\alpha}+d \choose \underline{\alpha}} = 6$. The second DAG is nearly identical to the first, except that we intentionally choose all parent sets so that they do not form norming sets. All these parent sets have identical geometric shapes, as shown in the left plot of Figure \ref{fig:norming}. Since the elements form a rectangle, we name these DAGs as Rectangle DAGs (denoted ``Rect'' in the figures).

We compare Vecchia GPs with layered norming DAGs and Rectangle DAGs from two perspectives: posterior estimation error and computational stability at large sample sizes, see Figure \ref{fig:norming}. Since they both approximate their mother GPs poorly and have almost the same computational time, we omit those quantities. 

From the figure, one can observe that Vecchia GPs with Norming DAGs generally perform better in terms of posterior estimation error than Rectangle DAGs. 
In light of our theoretical results, this difference in approximation accuracy stems from using all quadratic polynomials (Norming DAGs) versus only a subspace of them (Rectangle DAGs). Specifically, the norming set in the left panel of Figure \ref{fig:norming} uniquely determines all functions in the linear space 
$$\mathcal{P}_2(\Omega) = \mr{span}\{1,x[1],x[2],x[1]^2,x[1]x[2],x[2]^2\},$$
whereas the rectangle set in the same panel uniquely determines only functions in the space
$$\tilde{\mathcal{P}}_2(\Omega) = \mr{span}\{1,x[1],x[2],x[1]^2,x[1]x[2]\}.$$
The latter space excludes quadratic terms in $x[2]$, resulting in large estimation error for the true regression function, which is twice differentiable in $x[2]$.

However, one can also observe that the difference in estimation error  gradually diminishes as $n$ increases. Although we do not have a definitive explanation for this phenomenon, we conjecture that it is caused by numerical error in the posterior sampling algorithms. As discussed in Section \ref{sec:post inf}, a conjugate gradient step with a fixed numerical tolerance is employed to solve the sparse linear inverse problem, introducing a fixed numerical error into all posterior estimates. As sample size $n$ increases, the statistical error eventually decreases to the point where this numerical error dominates, causing the estimation errors for the Norming DAGs and Rectangle DAGs to become indistinguishable.

Another issue with non-norming sets is the numerical instability for large sample sizes, or in other words, the instability in the flat limit. In view of equations (\ref{eq:vecchia exp}) and (\ref{eq:vecchia var}), computing the conditional distributions of Vecchia GPs requires inverting the Mat\'{e}rn covariance matrix $K_{\mr{pa}(X_i),\mr{pa}(X_i)}$ for every $X_i\in\X_n$. As $n$ increases, however, the eigenvalues of these matrices may approach zero, rendering the inversion numerically unstable. To illustrate this, we plot the minimal eigenvalues of the matrices $K_{\mr{pa}(X_i),\mr{pa}(X_i)}$, taken over all locations $X_i\in\X_n$, in the right panel of Figure \ref{fig:norming}. One can observe that the difference between the two methods is substantial. For example, when $n=4,198,401$, the minimal eigenvalue of the Rectangle DAG is only $1/238$ of that for the Norming DAG. In fact, the Rectangle DAG exhibits significantly smaller minimal eigenvalues across all sample sizes and reaches numerical singularity much sooner than the Norming DAG. Recall that in Section \ref{sec:cond}, all results concerning the flat limit of Mat\'{e}rn processes require the parent sets to form norming sets. When the norming condition is violated, it is unclear whether the flat limits still exist, and if they do,what their structure might be. Our numerical results, however, serve as a warning in this regard.

\section{Background on polynomial interpolation}\label{sec:poly2DAG}
In this section we provide a brief introduction to polynomial interpolation. Some of the lemmas have already appeared in the literature, often without proofs. For completeness, we collect these results, cite the monographs where they are stated, and provide proofs for them as well.

\subsection{Vandermonde matrices and polynomial interpolation}\label{sec:vandermonde}
We first introduce the multidimensional version of the Vandermonde matrix. For a finite subset $A=\{w_1,w_2,\cdots,w_m\}\subset\mathbb{R}^d$ and ordered multi-indices $k_{(1)},...,k_{(m)}\in\mathbb{N}^d$ (with respect to lexicographical ordering), the $d$-dimensional Vandermonde matrix is defined as
\begin{equation}\label{def:vandemonde}
V_{A} = \begin{bmatrix}
1 & 1 & \cdots & 1 \\
w_1^{k_{(2)}} & w_2^{k_{(2)}} & \cdots & w_m^{k_{(2)}} \\
w_1^{k_{(3)}} & w_2^{k_{(3)}} & \cdots & w_m^{k_{(3)}} \\
\vdots & \vdots & & \vdots \\
w_1^{k_{(m)}} & w_2^{k_{(m)}} & \cdots & w_m^{k_{(m)}}
\end{bmatrix},
\end{equation}
where $w_i^{k_{(j)}}=\prod_{\ell=1}^d w_i[\ell]^{k_{(j)}[\ell]}$ and recall that $k_{(1)}=(0,0,\cdots,0)^T$. Note that the standard, one-dimensional Vandermonde matrix is a special case of the above $d$-dimensional version. Furthermore, for $x\in\mathbb{R}^d$, we define the vector $v_{x}$ as
\begin{equation}\label{eq:v vector}
v_{x} = (1,{x}^{k_{(2)}}, {x}^{k_{(3)}}, \cdots , {x}^{k_{(m)}})^T. 
\end{equation}

Next, we discuss polynomial interpolation. Let $\mathscr{P}_l(\Omega)$ denote the space of polynomials of order (not exceeding) $l$ on $\Omega\subseteq\mathbb{R}^d$. Given $m$ interpolation points
$$\{(w_i,y_i):1\le i\le m, w_i\in\mathbb{R}^d, y_i\in\mathbb{R}\}$$
and an order $l\in\mathbb{N}$, the polynomial interpolation problem consists of finding all the polynomials $P\in\mathscr{P}_l(\mathbb{R}^d)$ such that 
$$P(w_i)=y_i, \forall 1\le i\le m.$$
Here, we focus on \textit{unisolvency,} when the above problem has a unique solution. This topic has a history spanning over 100 years, see for instance the surveys \cite{gasca2000polynomial} and \cite{gasca2001history}. The connection between unisolvency and Vandermonde matrices has been stated in several sources, for example, in Proposition 3.5 of \cite{wendland2004scattered}. 
\begin{lemma}\label{lem:poly interp}
Let $m={l+d \choose l}$, $A=\{w_1,w_2,\cdots,w_m\}\subset\mathbb{R}^d$ and $Y=(y_1,y_2,\cdots,y_m)^T\in\mathbb{R}^m$. Then there exists a unique polynomial $P\in\mathscr{P}_l(\mathbb{R}^d)$ satisfying
$P(w_i)=y_i$ if and only if the Vandermonde matrix $V_A$ is invertible. Moreover, this polynomial can be expressed as
\begin{equation}\label{eq:poly i form}
P(x) = Y^T V_A^{-1} v_x.
\end{equation}
\end{lemma}
\begin{proof}
We prove the equivalence statement by contradiction. Suppose the Vandermonde matrix $V_A$ is not invertible. Since $V_A$ is a square matrix, there exists a nonzero vector $b$, such that $V_A b =0$. Assume furthermore, that there exists a polynomial $P(x)$ satisfying ${P}(w_i)=y_i,\forall 1\le i\le m$; otherwise, the statement holds. Then the polynomial $\tilde{P}(x) = P(x) + \sum_{i=1}^m b[i] x^{k_{(i)}}$ also satisfies $\tilde{P}(w_i)=y_i,\forall 1\le i\le m$, and thus, in view of $b\ne 0$, there exist at least two polynomails satisfying the interpolation condition. This contradicts uniqueness, proving the equivalence.

Next, assume that the Vandermonde matrix $V_A$ is invertible. Then, equation (\ref{eq:poly i form}) provides a solution. We prove by contradiction that this solution is unique. Suppose that there are at least two distinct polynomials $P_1\ne P_2$ satisfying $P_j(w_i)=y_i$ for $j=1,2$. Then $P_1-P_2\ne 0$ and $(P_1-P_2)(w_i)=0,\forall 1\le i\le m$. Writing the difference in the form $(P_1-P_2)(x) = \sum_{i=1}^m b[i] x^{k_{(i)}}$ and introducing the notation $b = (b[1],b[2],\cdots ,b[m])^T\in\mathbb{R}^m$, $b\neq 0$, we obtain  $V_A b = 0$, which contradicts with the assumption that $V_A$ is invertible. Hence the solution is unique.
\end{proof}

\subsection{Results regarding Norming sets}
In view of Lemma \ref{lem:poly interp}, the unisolvency of the polynomial interpolation problem is equivalent to the invertibility of the corresponding Vandermonde matrix.  Building on this, and recalling the notion of norming sets introduced in Section \ref{sec:norming}, we aim to establish qualitative relationships among the properties of the polynomials, the minimal singular value of the Vandermonde matrices, and the spatial dispersion of the set $A$, beginning with the following lemma.

\begin{lemma}\label{lem:norming def}
Let $\Omega= [0,1]^d$. For $l\in\mathbb{N}$, and a set $A=\{w_1,w_2,\cdots, w_m\}$ with cardinality $m={l+d \choose l}$, the following two statements are equivalent:
\begin{enumerate}[(1)]
\item The set $A$ is a norming set for $\mathscr{P}_l(\Omega)$ with norming constant $c_N>0$;
\item The minimal singular value of the Vandermonde matrix $V_A$ is bounded from below by $c_S>0$.
\end{enumerate}
Moreover, there exists a universal constant $C_0\geq1$, independent of $A$, such that
$$C_0^{-1} c_N\leq c_S^{-1}\leq C_0 c_N. $$
\end{lemma}
\begin{proof}
We first prove (2) $\Rightarrow$ (1). Since the smallest singular value of the Vandermonde matrix $V_A$  is positive, $V_A$ is invertible. In view of Lemma \ref{lem:poly interp}, the set $A$ is unisolvent with respect to $\mathscr{P}_l(\mathbb{R}^d)$ and, therefore, also unisolvent for $\mathscr{P}_l(\Omega)$. Consequently, every polynomial $P\in\mathscr{P}_l(\mathbb{R}^d)$ can be expressed in the form
\begin{equation}\label{eq:unisolvent P}
P(x) = P(A)^T V_A^{-1} v_x.
\end{equation}
This allows to control $P(x)$, $\forall x\in\Omega$ as
\begin{align*}
|P(x)| \le & \|P(A)\|_2 \|V_A^{-1}\|_2 \|v_x\|_2 \le 
\sqrt{ m } \max_{x'\in A} |P(x')| \cdot c_S^{-1} \cdot \sqrt{m} \le 
{l+d \choose l} c_S^{-1} \max_{x'\in A}| P(x')|,
\end{align*}
where we utilized that $\|v_x\|_2 \le \sqrt{m}$ for $x\in\Omega = [0,1]^d$. Thus, $A$ is a norming set with norming constant 
\begin{equation*}
c_N \le {l+d \choose l} c_S^{-1}.
\end{equation*}

We proceed to prove (1) $\Rightarrow$ (2). 
First we verify by contradiction that the Vandermonde matrix $V_A$ is invertible. Suppose it is not. Then there exists a non-zero vector $b\in\mathbb{R}^m$ such that $V_A^T b = 0$. The corresponding polynomial $P(x) = b^T v_x\in\mathcal{P}_l(\Omega)$ satisfies $P(x)=0$,  $\forall x\in A$, i.e. $\|P\|_A=0$.
However, since $b\neq 0$, $\exists x'\in\Omega$ such that $P(x')\ne 0$, i.e. $\sup_x|P(x)|>0$. This contradicts the assumption that $A$ is a norming set.

To establish $c_S^{-1}<C_0 c_N$, we introduce a ``standard'' norming set $A^*$ to serve as a ``bridge'' between $c_S$ and $c_N$. Define
$$A^* = \{(x\in\mathbb{R}^d: x[i]\in\{0/l,1/l,\cdots,l/l\},\forall 1\le i\le d; \|x\|_1 \le 1\}.$$
The set $A^*$ is a corner set, as in Section 3 of \citep{neidinger2019multivariate} with cardinality $m={l+d \choose l}$, and is unisolvent. 
Denote the minimal singular value of $V_{A^*}$ by $c^*$. Then for all $u\in \mathbb{R}^m$ with $\|u\|_2=1$, 
\begin{equation*}
c^* \le \|u^T V_{A^*}\|_2 = \sqrt{\sum_{x'\in A^*}(v_{x'}^T u)^2} \le \sqrt{m} \sup_{x'\in A^*}|v_{x'}^T u| \le \sqrt{m} \sup_{x\in \Omega}|v_{x}^T u|.
\end{equation*}
Therefore, for all $u\in \mathbb{R}^m$ with $\|u\|_2=1$,
\begin{equation}\label{eq:vx spread}
\sup_{x\in\Omega} |v_x^T u| \ge c^*/\sqrt{m}.
\end{equation}

Next, for the Vandermonde matrix $V_A$ associated with the set $A$, consider its singular value decomposition (SVD) $V_A = U_L D U_R^T$, where both $U_L$, $U_R$ are unitary matrices, and $D$ is a diagonal matrix with positive diagonal entries. Let  $u_{l,m}$ and $u_{r,m}$ denote the left and right singular vectors corresponding to the minimal singular value of $V_A$, respectively.
Let $P^*$ be the interpolation polynomial such that $P^*(A) = u_{r,m}$. Then $P^*(x) = P^*(A)^T V_A^{-1} v_x$. By the lower bound in (\ref{eq:vx spread}), 
\begin{align}\label{eq:cs lower}
\sup_{x\in\Omega} |P^*(x)| &= \sup_{x\in\Omega} |P^*(A)^T V_A^{-1} v_x| = \sup_{x\in\Omega} |u_{r,m}^T U_R D^{-1} U_L^T v_x|\nonumber\\
& = c_S^{-1}\sup_{x\in\Omega} |u_{l,m}^T v_x| \ge c_S^{-1}c^*/\sqrt{m}.
\end{align}
Also, by the definition of norming set, we have
\begin{equation}\label{eq:cs upper}
\sup_{x\in\Omega} |P^*(x)| \le c_N \sup_{x'\in A} |P^*(x')| = c_N \|u_{r,m}\|_\infty \le c_N.
\end{equation}
Combining equations (\ref{eq:cs lower}) and (\ref{eq:cs upper}), we have
$$c_S^{-1} \le c_N\sqrt{m} / c^*.$$
Since the constant $c^*$ is independent of the set $A$, we finished the proof.
\end{proof}

While the equivalence between norming sets and unisolvent sets is well known in the literature, the above lemma characterizes the qualitative relationship between the constants $c_N$ and $c_S$. This explicit dependence plays a crucial role in computing the small deviation bounds for Vecchia GPs in Section \ref{sec:prob}. We note, that an important difference between Lemma \ref{lem:poly interp} and Lemma \ref{lem:norming def} lies in the support of the polynomials. In view of the unboundedness of polynomials  on the whole domain $\mathbb{R}^d$, finite norming sets do not exist when $\Omega=\mathbb{R}^d$.

The norming constant $c_N$ has  an explicit expression and is  ``scale-invariant'' as stated in the following lemma, which is also discussed in \cite{bos2018fekete}.
\begin{lemma}\label{lem:norming invariant}
Let $\Omega$ be a compact subset of $\mathbb{R}^d$ and let $A\subset\Omega$ with cardinality $m={l+d \choose l}$ be a norming set for $\mathscr{P}_l(\Omega)$ with norming constant $c_N$. Then we have
\begin{equation}\label{eq:compute cn}
c_N = \sup_{x\in\Omega} \|V_A^{-1} v_x\|_1.
\end{equation}
Moreover, for all $\tau>0$, the set $\tau A$ is a norming set for $\mathscr{P}_l(\tau\Omega)$ with the same norming constant $c_N$.
\end{lemma}
\begin{proof}
In view of Lemma \ref{lem:norming def}, the Vandermonde matrix $V_A$ is invertible. Thus following from Lemma \ref{lem:poly interp}, every polynomial $P\in\mathscr{P}_l(\Omega)$ can be expressed as $P(x) = P(A)^T V_A^{-1} v_x$, $\forall x\in\Omega$, which implies
$$|P(x)| \le \|P(A)\|_\infty \sup_{x\in\Omega}\|V_A^{-1} v_x\|_1.$$
Hence, we obtain the bound $c_N \le \sup_{x\in\Omega}\|V_A^{-1} v_x\|_1.$
On the other hand, since $\Omega$ is compact, $\exists x_0\in\Omega$ satisfies $\|V_A^{-1} v_{x_0}\|_1 = \sup_{x\in\Omega}\|V_A^{-1} v_x\|_1$. Let us take the unique polynomial $P\in \mathscr{P}_l(\Omega)$ satisfying $P(A)[i] = (\mr{sign}(V_A^{-1} v_{x_0}))[i]$. Then
$$P(x_0) = P(A)^T V_A^{-1} v_{x_0} = \|V_A^{-1} v_{x_0}\|_1 = \|V_A^{-1} v_{x_0}\|_1 \|P(A)\|_\infty.$$
This provides the lower bound $c_N \ge \sup_{x\in\Omega}\|V_A^{-1} v_x\|_1$ and hence verifies \eqref{eq:compute cn}.

Next, consider a scale parameter $\tau>0$. Define the diagonal matrix $\Gamma\in\mathbb{R}^{m\times m}$ with diagonal entries $\Gamma[j,j]=\tau^{|k_{(j)}|}$, where $k_{(j)}$ is the $j$th multi-index in $\mathbb{N}^d$ with respect to lexicographical ordering. Note that $V_{\tau A} = \Gamma V_A$ and $v_{\tau x} = \Gamma v_x$. Thus the Vandermonde matrix $V_{\tau A}$ is also invertible, and $\forall\tau x\in \tau \Omega$,
$V_{\tau A}^{-1} v_{\tau x} = V_A^{-1} v_x$.
Therefore, the scaled set $\tau A$ is a norming set for $\mathscr{P}_l(\tau\Omega)$ with the same norming constant $c_N$.
\end{proof}

While the above lemma provides the explicit formula (\ref{eq:compute cn}) for the norming constant $c_N$, computing it can still be challenging, since it requires inverting the corresponding Vandermonde matrix and taking a supremum over the entire compact set $\Omega$.  In the unidimensional case, $\Omega =[0,1]$, we can provide an alternative formulation. Let $A=\{w_1,w_2,\cdots, w_{l+1}\}\subset[0,1]$ consist of $l+1$ distinct points. Then $c_N$ can be expressed using Lagrange polynomials as
\begin{equation}\label{eq:c_N d=1}
c_N =\sup_{x\in [0,1]} \sum_{j=1}^{l+1} \bigg| \frac{\prod_{i\ne j}(x-w_i)}{\prod_{i\ne j}(w_j-w_i)} \bigg|.
\end{equation}
 In higher dimensions, however, without additional technical assumptions on the set $A$, equation (\ref{eq:compute cn}) can not be further simplified.

\begin{figure}[h!]
    \centering
    \includegraphics[width=\linewidth]{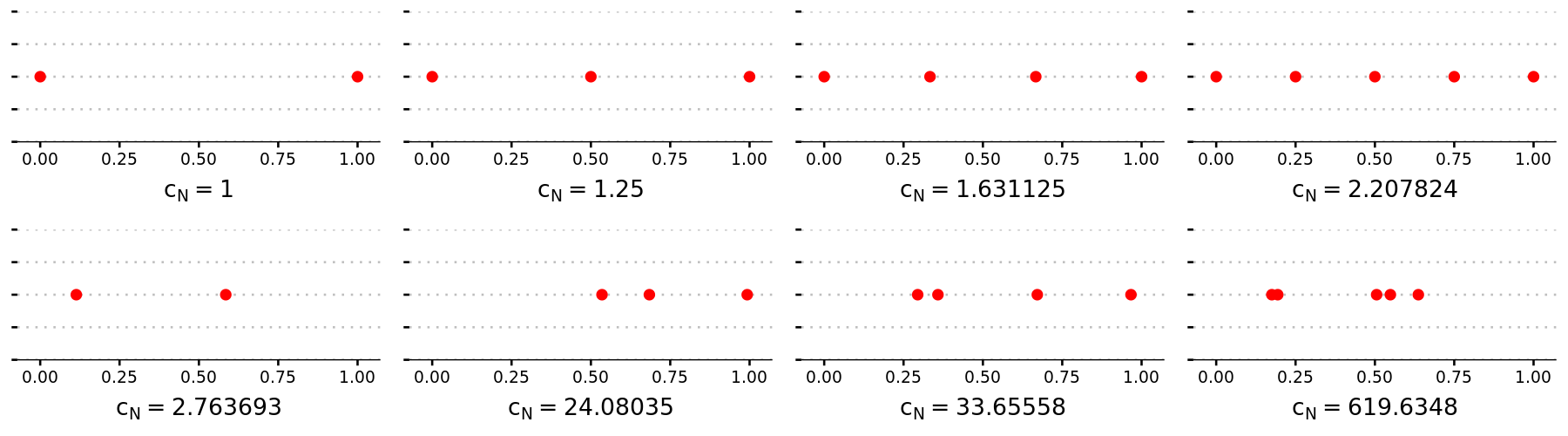}
    \caption{Norming sets for a one-dimensional space for polynomials of order $1,2,3$ and $4$, respectively. Top row: equally spaced points on $[0,1]$; bottom row: i.i.d. samples from the uniform distribution on $[0,1]$.}
    \label{fig:norming d=1}
\end{figure}

To better understand the properties of norming sets and their norming constants, we graphically illustrate several natural candidates. We begin with the one-dimensional case, where norming sets  for $\ \mathscr{P}_l([0,1])$ are exactly the finite sets with cardinality $l+1$. In Figure \ref{fig:norming d=1}, we plot several norming sets with their associated norming constants.

First, note that as the polynomial order $l$, and consequently the cardinality of the corresponding norming sets, increases, the norming constants also increase. This is intuitive: as the dimension of the polynomial space grows with $l$, it becomes more difficult to characterize all polynomials using a finite set of functionals. Second, one can observe that for randomly chosen points $w_j$, the norming constant is substantially larger than for the equidistant design.

\begin{figure}[h!]
    \centering
    \includegraphics[width=\linewidth]{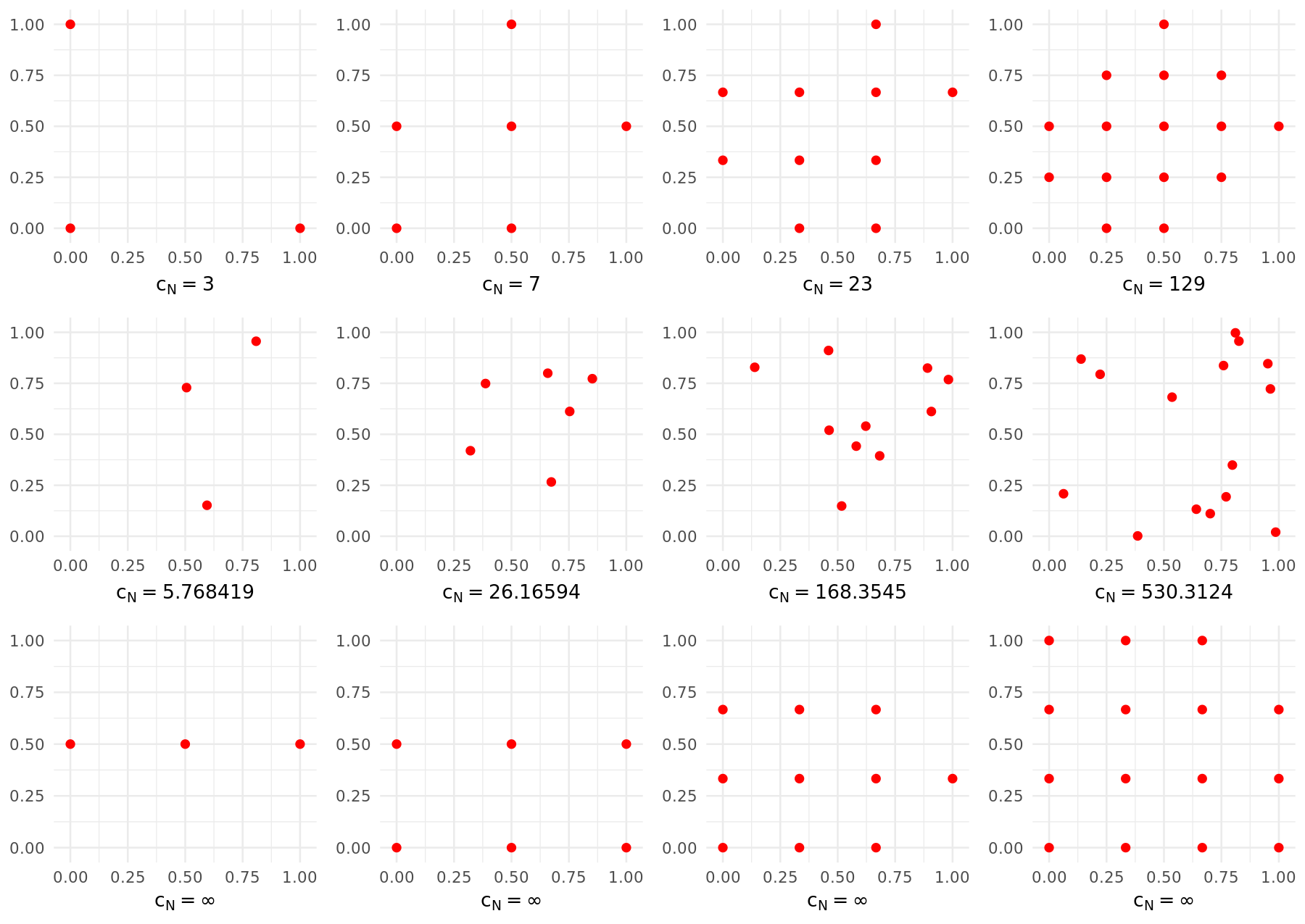}
    \caption{Norming sets in $[0,1]^2$. The cardinality of the sets in the $l$th column  ($1\le l\le 4$) is ${l+2 \choose l}$, which equals the dimension of the polynomial space $\mathscr{P}_l(\Omega)$.  Top row: corner sets; middle row: i.i.d. samples from the uniform distribution; bottom row: non-norming sets for $\mathscr{P}_l(\Omega)$.}
    \label{fig:norming d=2}
\end{figure}

Next, we consider the case $d=2$, see Figure \ref{fig:norming d=2}. Similarly to the one-dimensional case,the norming constants increase with the order of the polynomials. Furthermore, the norming constants are smaller for the so called \textit{corner sets} than for  i.i.d. samples generated from a uniform distribution. A corner set is a  carefully designed unisolvent set that captures the local polynomial structure, see Section \ref{sec:build dag} for its definition and further discussion. Moreover, for $d\ge 2$, unlike in the one-dimensional case, not every set with cardinality ${l+d \choose l}$ is a norming set for $\mathscr{P}_l(\Omega)$. Examples of such non-norming sets are shown in the bottom row of Figure \ref{fig:norming d=2}.

\section{Algorithms for building layered norming DAGs}\label{sec:algs}

\subsection{Algorithms for building layered norming DAGs}\label{sec:build dag alg}

We provide Algorithm \ref{alg:dag grid} used to construct the layered norming DAGs discussed in Section \ref{sec:build dag}. 
\begin{algorithm}[h!]
\SetAlgoLined
Input: grid data $\X_n$ with $n=\tilde{n}^d$ and order $l$. Set $m={l+d \choose l}$ and let $r$ be the largest integer such that $\tilde n\ge 2^r+1$.\\
Define the set $\bar{\X}_n$ as in (\ref{eq:subgrid}). Choose layers $\N_j, 1\le j\le r$ as in (\ref{eq:layer grid}), satisfying $\cup_{j=0}^r \N_j = \bar{X}_n$.\\
\eIf{$\tilde{n}>2^r+1$}{
Let $J=r+1$. Set $\N_{r+1}=\X_n\backslash\bar{\X}_n$.
}{
Let $J=r$.
}
Let $j_0=\arg\max_{j}\{2^{j-1}+1<l\}$. Set $\mr{pa}(X_1)=\emptyset$.\\
\For{$1\le j\le j_0$}{
For all $X_i\in\N_j$, set $\mr{pa}(X_i)=\cup_{j'=0}^{j-1}\N_{j'}$. 
}
\For{$j_0+1\le j \le J$}{
\For{$X_i\in\N_j$}{
\For{$1\le h\le d$}{
Order the projection of $\cup_{0\le j'\le j-1}\N_{j'}$ onto the $h$th coordinate by increasing distance to $X_i[h]$ as $B_h=(x_{s(1),h},x_{s(2),h},\cdots,x_{s(n'),h})$.
}
Set $\mr{pa}(X_i)$ as in equation (\ref{eq:parent grid}).
}
}
\KwResult{Return $\mr{pa}(X_i),\forall 1\le i\le n$.}
\caption{Construct layered norming DAGs on Grid.}\label{alg:dag grid}
\end{algorithm}

We briefly discuss the computational complexity of Algorithm \ref{alg:dag grid}. 
We first consider the case where the dataset $\D_n$ is stored in a structured format, so that retrieving data at any grid point takes $O(1)$ time. The grid structure ensures a simple, explicit form for the layers $\N_j$, without requiring costly numerical methods. Thus the computational complexity of building the layers is $O(n)$.
 Allocating parent sets requires sorting with respect to each coordinate. However, due to the grid structure, the set $B_h, \forall 1\le h\le d$ and all $X_i$ has a simple, closed-form solution. Therefore, for each $X_i$, assigning its parent set of cardinality $m$ requires only $O(m)$ steps. Overall, the computational complexity of Algorithm \ref{alg:dag grid} is $O(nm)$, with $m = {\underline{\alpha}+d \choose \underline{\alpha} }$. If the dataset $\D_n$ is not stored in a structured format, one has to start with appropriately sorting it, increasing the overall computational complexity to $O(dn\ln n + nm)$.

\subsection{Building DAGs satisfying the Norming Condition \ref{cond:norming} for general $\X_n$}\label{sec:build dag general}
In this section, we discuss how to build DAGs satisfying Conditions \ref{cond:layered} and \ref{cond:norming} for general design points $\X_n$. Without loss of generality, we assume $\X_n\subset[0,1]^d$.

The DAG construction follows a three-steps procedure. First, we sort all elements in $\X_n$ according to maximin ordering, defined sequentially. The first element is selected arbitrarily. Once the first $i$ elements are chosen, the $(i+1)$th element is selected as the point in $\X_n$ that maximizes the minimal distance to the first $i$ elements, with respect to a metric specified by the user. The procedure is repeated until all $n$ elements are sorted.
With a slight abuse of notation, we denote the sorted elements as $X_1, X_2, \cdots, X_n$. 

Second, for all $i\ge 2$, let $\mr{dist}_i$ denote the distance between $X_i$ and the set $\{X_{i'},1\le i'< i\}$ as . Let $\N_0 = \{X_1\}$. Then, for all $j\ge 1$, define the layers $\N_j$ as
\begin{equation}\label{eq:layered}
\N_j = \{X_i\in\X_n\backslash\bigcup_{j'=0}^{j-1}\N_{j'}: \mr{dist}_i\in (c_d\gamma^{-(j+1)}, c_d\gamma^{-j}]\},
\end{equation}
for a given constant $\gamma>1$.

Finally, once the layers are specified,  for all $j$ large enough such that $\cup_{0\le j'\le j-1}\N_{j'}$ contains at least $l^d$ elements, for each $X_i\in\N_j$, we choose the corresponding parent set $\mr{pa}(X_i)$ via the following iterative procedure. 
Let $c_S^*>0$ be a given constant and set $c_S = c_S^*$. Initialize $\mr{pa}(X_i)$ with the element in $\cup_{j'=0}^{j-1} \N_{j'}$ that is closest to $X_i$. Suppose that after $m'\in\{1,...,m-1\}$ iterations, the parent set is $\mr{pa}(X_i)=\{w_1,w_2,\cdots,w_{m'}\}$. Recall that $k_{(i)}, 1\le i\le m$, denotes the $i$th multi-index in lexicographical ordering ``$\prec$'' and define the diagonal matrix $\Gamma=\mr{diag}\{1, \gamma^{j|k_{(2)}|}, \cdots, \gamma^{j|k_{(m)}|} \}\in\mathbb{R}^{m\times m}$. Then, we add a new element $X'$ to $\mr{pa}(X_i)$ if the matrix
\begin{equation}\label{eq:nonsquare van}
M=\Gamma\;[v_{w_1-X_i} \;v_{w_2-X_i} \cdots v_{w_{m'}-X_i} \;v_{X'-X_i}]
\end{equation}
satisfies $\sigma_{\min}(M) \ge c_S$, where $\sigma_{\min}(M)$ denotes the minimum singular value of $M$, and $v_x$ is defined in (\ref{eq:v vector}). If all $m$ elements are selected this manner, the process terminates. Otherwise, we set $c_S=c_S/2$ and repeat the procedure with this updated value. The process terminates either when a parent set of size $m$ is constructed, or when $c_S$ drops below a prescribed, numerical tolerance level. In the latter case, we simply fill up the parent set with the closest data points to $X_i$ until the prescribed cardinality $m$ is reached. We formalize this procedure in Algorithm \ref{alg:dag general}.

\begin{algorithm}[h]
\SetAlgoLined
Input: data set $\X_n$, initial singular value $c_S^*>0$, constant $\gamma>1$ for the layers, order $l$ for polynomials and numerical tolerance $\epsilon_{\mr{tol}}$. Set $m={l+d \choose l}$.\\
Randomly choose an initial location and sort all elements in $\X_n$ by maximin ordering. With a slight abuse of notation, we denote the sorted elements by $X_1, X_2, \cdots, X_n$.\\
Set layers $\N_0=\{X_1\}$ and $\N_1, \N_2,\cdots \N_J$ according to (\ref{eq:layered}).\\
Let $j_0=\arg\max_{j}\{|\cup_{j'=0}^{j-1}\N_{j'}|<m\}$. Set $\mr{pa}(X_1)=\emptyset$.\\
\For{$1\le j\le j_0$}{
For all $X_i\in\N_j$, let $\mr{pa}(X_i)=\cup_{j'=0}^{j-1}\N_{j'}$. 
}
\For{$j_0+1\le j \le J$}{
\For{$X_i\in\N_j$}{
Sort the elements of $\cup_{j'=0}^{j-1}\N_{j'}$ in order of increasing distance to $X_i$, denoted as $B=(X_{i_1},X_{i_2},\cdots, X_{i_{n'}})$, where $n'=|\cup_{j'=0}^{j-1}\N_{j'}|$.\\
Set $c_S = c_S^*$, $m'=0$ and $\mr{pa}(X_i)=\emptyset$.\\
\While{$c_S>\epsilon_{\mr{tol}}$}{
\For{$X'\in B$}{
\If{
$\lambda_{\min}(M)\geq c_S$, for matrix $M$ given in (\ref{eq:nonsquare van}) 
\\
}{
Let $m'=m'+1$, $w_{m'}=X'$,  $\mr{pa}(X_i) = \mr{pa}(X_i)\cup\{w_{m'}\}$ and $B=B\backslash \{X'\}$.
}
\If{$m'=m$}{
\textbf{Break}
}
}
$c_S = c_S/2$.
}
\If{$m'<m$}{
Add the first $m-m'$ elements of the set $B$ to $\mr{pa}(X_i)$.
}
}
}
\KwResult{Return $\mr{pa}(X_i),\forall 1\le i\le n$.}
\caption{Constructing layered norming DAGs on general positions.}\label{alg:dag general}
\end{algorithm}

We note that, although Algorithm \ref{alg:dag general} always produces a parent set of cardinality $m$, the resulting norming constant may be arbitrarily small depending on the properties of the dataset $\X_n$. This limitation is, however, unavoidable. For certain irregular datasets, it is possible that none of the candidate parent sets for $X_i$ are unisolvent, in which case the norming constant $c_N$ becomes effectively infinity. To achieve a tight control over the norming constant, additional structural assumptions on the data set $\X_n$ are required. For example, grid data (or its slightly perturbed version), as considered in the previous section, provides such strong control.

\section{Computational aspects}\label{sec:post inf}
In this section, we discuss the computational aspects of the posterior. We start by recalling the Bayesian hierarchical model introduced in Section \ref{sec:GP regression}. The observed data are $(X_i,Y_i)_{i=1,...,n}$, where $X_i$ are fixed design points in $[0,1]^d$ and
\begin{align*}
Y_i = f(X_i) + \epsilon_i, \quad \epsilon_i \overset{\mr{i.i.d}}{\sim} N(0, \sigma^2).
\end{align*}
The model parameters $f$ and $\sigma^2$ are endowed with a Vecchia Mat\'{e}rn GP and inverse-gamma prior, respectively, i.e.
\begin{align*}
\sigma^2 &\sim \mr{IG}(a_0, b_0),\\
f &\sim (\hat{Z}_X^{\tau,s}, X\in\X).
\end{align*}
Additionally, one can also introduce a hyper-prior on the scale parameters $\tau$ and $s$ (and/or the regularity parameter $\alpha$) of the Vecchia GP.

In the considered regression model with nuggets, i.e. observational errors $\eps_i$, the Cholesky factor of the posterior precision matrix is no longer sparse. As a result, we use a Gibbs sampling framework for inference. More concretely, we alternate sampling the variance parameter $\sigma^2$, the Vecchia GP $\hat{Z}_{\X_n}^{\tau,s}$, and any hyperparameters (if included) from their respective full conditional distributions.

Before providing further details, we introduce some additional notation. Let  $\hat{\Phi}_n$ denote the precision matrix of the Vecchia GP $\hat{Z}_X^{\tau,s}$ evaluated on the training data $\X_n$, and let $\mathcal{Y}_n = [Y_1, Y_2, \cdots, Y_n]^T$. 
For a Vecchia GP defined in Section \ref{sec:vecchia GP}, its DAG structure (or equivalent, the parent sets $\mr{pa}(X_i)$ for all $i$) uniquely determines a Cholesky decomposition of the precision matrix 
$$\hat{\Phi} = \hat{L} \hat{L}^T.$$ 
Moreover, if the DAG is sparse, the Cholesky factor $\hat{L}$ is also sparse.
For a detailed derivation, see Lemma S3 and its proof in \cite{zhu2024radial}. 

Next, we discuss sampling $\sigma^2$, $\hat{Z}_{\X_n}^{\tau,s}$ and any hyperparameters from their full conditional distributions. By conjugacy the full conditional distribution of $\sigma^2$ is also inverse-gamma, allowing for direct sampling. The most challenging and computationally intensive part is sampling $\hat{Z}_{\X_n}^{\tau,s}|\X_n,\Y_n,\sigma^2$, which follows a multivariate $n$-dimensional Gaussian distribution
\begin{align*}
\hat{Z}_{\X_n}^{\tau,s}|\X_n,\Y_n,\sigma^2  \sim N\Big(  (\hat{\Phi}+I_n/\sigma^2)^{-1}\Y_n /\sigma^2, (\hat{\Phi}+I_n/\sigma^2)^{-1}\Big).
\end{align*}
Note that the sparse Cholesky decomposition of $\hat{\Phi}$ does not directly induce a sparse Cholesky decomposition of $\hat{\Phi} + I_n/\sigma^2$. Instead, one can sample from this distribution by solving a linear inverse problem using conjugate gradient solver \citep{nishimura2023prior}. This approach only requires matrix-vector multiplications and can fully utilize the sparse structures of the Cholesky factor $\hat{L}$.

If the model includes hyper-priors on the GP parameters, e.g. scale or regularity, the corresponding hyper-posteriors are sampled using a Metropolis-Hasting step within Gibbs sampler. This requires evaluating both the prior density and the likelihood. The likelihood $\Y_n|\hat{Z}_{\X_n}^{\tau,s},  \sigma^2$ is a product of independent Gaussian densities. The GP prior density for $\hat{Z}_{\X_n}^{\tau,s}$, in view of the sparse Choleksy decomposition of $\hat{\Phi}$, can be efficiently computed as
\begin{equation}\label{eq:prior chol}
p(\hat{Z}_{\X_n}^{\tau,s}) = (2\pi)^{-n/2} \det(\hat{L})^2 \exp\big[-(\hat{L} \hat{Z}_{\X_n}^{\tau,s})^2 \big],
\end{equation}
where the determinant reduces to the product of the diagonal entries of $\hat{L}$ due to its triangular structure.

We summarize the full sampling scheme in Algorithm \ref{alg cg}. Note that the method is fairly general and applies to any Vecchia GPs with a well-defined DAG structure -- not only the layered norming DAGs and Mat\'{e}rn mother GPs proposed in this work.

\begin{algorithm}[ht]
\SetAlgoLined
Input: training data $\X_n, \Y_n$; regularity $\alpha$ and scale parameters $\tau,s$ of the GP (or their hyper-priors); number of MCMC iterations $n_{mc}$.\\
Build a DAG on the training locations $\X_n$. \\
\For{$1\le i_{mc} \le n_{mc}$}{
Sample $\sigma^2$ from its full conditional distribution;\\

Generate random vector ${W}\sim N(\Y_n/\sigma^2, \hat{\Phi}+{I}_n/\sigma^2)$ via 
${W} = \Y_n/\sigma^2 + \hat{{L}}{W}_1 +  {W}_2/\sigma,$   
where ${W}_1, {W}_2\stackrel{iid}{\sim}N(0,{I}_n)$ and $\hat{L}$ is the Cholesky factor of the precision matrix $\hat{\Phi}$; \\
Obtain a sample of $\hat Z_{\X_n}^{\tau,s}$ by solving the linear system $(\hat{\Phi}+{I}_n/\sigma^2)\hat Z_{\X_n}^{\tau,s}  = {W}$ using the conjugate gradient method;\\
If hyper-priors are present, sample the hyperparameters, where the prior density of $\hat{Z}_{\X_n}^{\tau,s}$ is computed via \eqref{eq:prior chol}.
}
\KwResult{Output the posterior samples of $\sigma^2$, $\hat{Z}_{\X_n}^{\tau,s}$, and hyperparameters (if applicable).}
\caption{Posterior sampling for Vecchia GP}\label{alg cg}
\end{algorithm}
\vspace{0.2cm}

\textbf{Computational complexity} 
Considering a sparse DAG structure with parent set size bounded by $m$ -- e.g., for the layered norming DAG introduced in Section \ref{sec:layered norming dag}, we have $m={\underline{\alpha}+d \choose \underline{\alpha}}$ -- the computational complexity of Algorithm \ref{alg cg} is $O(n m n_{cg} n_{mc} )$, where $n_{mc}$ denotes the number of MCMC iterations and $n_{cg}$ is the average number of conjugate gradient steps required per iteration. The term $nm$ corresponds to the cost of a single conjugate gradient step, as each matrix-vector multiplication involving the sparse Cholesky factor $\hat{L}$ can be carried out in $O(nm)$ time.

\section{Proof of the lemmas and theorems in Section \ref{sec:matern theory} and \ref{sec:cond}}\label{sec:proof prob 1}
\subsection{Proof of Lemma \ref{lem:matern smooth}}
\begin{proof}
By the definition of the Mat\'{e}rn covariance kernel in (\ref{eq:MaternCov}), for all multi-indices $k_1,k_2$ satisfying $|k_1|+|k_2|\le 2 \underline{\alpha}$,
\begin{align*}
|K^{(k_1,k_2)}(x_1,x_2)| \propto & \left|\int_{\mathbb{R}^d} e^{-\iota\langle x_1-x_2,\xi\rangle} (-\iota)^{|k_1|}\xi^{k_1} \iota^{|k_2|} \xi^{k_2} (1+\|\xi\|_2^2)^{-(\alpha+d/2)} d\xi \right|\\
\le & \int_{\mathbb{R}^d} (1+\|\xi\|_2^2)^{-(\alpha-\underline{\alpha}+d/2)} d\xi \lesssim 1,
\end{align*}
proving assertion (\ref{eq:matern diff}).

Next, we consider the derivative $K^{(k,k)}(x_1,x_2)$ for $|k|\le\underline{\alpha}$, which can be directly computed as
\begin{equation*}
K^{(k,k)}(x_1,x_2) \propto \int_{\mathbb{R}^d} e^{-\iota\langle x_1-x_2,\xi\rangle} \xi^{2k} (1+\|\xi\|_2^2)^{-(\alpha+d/2)} d\xi.
\end{equation*}
It is straightforward to see that $K^{(k,k)}(x_1,x_1)=K^{(k,k)}(x_2,x_2),\;\forall x_1,x_2\in\mathbb{R}^d$. 
Since the function $\xi\rightarrow \xi^{2k} (1+\|\xi\|_2^2)^{-(\alpha+d/2)}$ is even, the equation $K^{(k,k)}(x_1,x_2)=K^{(k,k)}(x_2,x_1)$ also holds, providing (\ref{eq:matern diff symmetry}).

Finally, we consider the derivative $K^{(k_1,k_2)}(x_1,x_2)$ for $|k_1|+|k_2|=2\underline{\alpha}$.
Since $K(\cdot,\cdot)$ and its derivatives are real valued functions,
\begin{align*}
K^{(k_1,k_2)}(x,x)&-K^{(k_1,k_2)}(x,x+h)\\
\propto & \;\mr{Re}\bigg(\int_{\mathbb{R}^d}(1-\expinner{-h}) \xi^{k_1+k_2} (-\iota)^{|k_1|} \iota^{|k_2|} (1+\|\xi\|_2^2)^{-(\alpha+d/2)} d\xi \bigg) \\
= & (-1)^{|k_1|+\underline{\alpha}}\int_{\mathbb{R}^d}\mr{Re}(1-e^{\iota\inner{h,\xi}}) \xi^{k_1+k_2}  (1+\|\xi\|_2^2)^{-(\alpha+d/2)} d\xi\\
= & \frac{(-1)^{|k_1|+\underline{\alpha}}}{2} \int_{\mathbb{R}^d} (1-e^{-\iota\langle h,\xi\rangle}) \overline{(1-e^{-\iota\langle h,\xi\rangle})} \xi^{k_1+k_2} (1+\|\xi\|_2^2)^{-(\alpha+d/2)} d\xi.
\end{align*}
Therefore, the absolute value of $K^{(k_1,k_2)}(x,x)-K^{(k_1,k_2)}(x,x+h)$, is bounded as
\begin{align}\label{eq:K derive drop}
&\big|K^{(k_1,k_2)}(x,x)-K^{(k_1,k_2)}(x,x+h)\big|\\
&\qquad\quad\lesssim  \int_{\mathbb{R}^d} (1-e^{-\iota\langle h,\xi\rangle}) \overline{(1-e^{-\iota\langle h,\xi\rangle})}  (1+\|\xi\|_2^2)^{-(\alpha-\underline{\alpha}+d/2)} d\xi.\nonumber
\end{align}
Noting that the right hand side (RHS) of (\ref{eq:K derive drop}) (up to a constant multiplier) coincides with $K(x,x)-K(x,x+h)$ for regularity parameter $\alpha-\underline{\alpha}$, it is sufficient to prove the upper bound (\ref{eq:matern diff lipschitz}) for $\alpha\in(0,1]$. From now on we will apply the representation (\ref{eq:MaternCov o}) of the covariance kernel. For $\alpha\in(0,1)$, combining equations (10.25.2) and (10.27.4) of \cite{DLMF} yields the following expansion for the Mat\'{e}rn covariance kernel near zero
\begin{align*}
K(x,x+h) = &\frac{2^{1-\alpha}}{\Gamma(\alpha)} \|h\|_2^\alpha \mathcal{K}_\alpha(\|h\|_2) \\
= & \frac{2^{-\alpha}\pi}{\Gamma(\alpha)\sin(\alpha \pi)} \|h\|_2^\alpha [\mathcal{I}_{-\alpha}(\|h\|_2) - \mathcal{I}_{\alpha}(\|h\|_2)] \\
= & \frac{\pi}{\Gamma(\alpha)\sin(\alpha \pi)} 
\left[ \sum_{i=0}^\infty \frac{4^{-i}\|h\|_2^{2i}}{i!\Gamma(-\alpha+i+1)} - 2^{-2\alpha}\|h\|_2^{2\alpha} \sum_{i=0}^\infty \frac{4^{-i}\|h\|_2^{2i}}{i!\Gamma(\alpha+i+1)} \right] \\
= & 1 - \frac{\pi}{2^{2\alpha}\Gamma(\alpha)\Gamma(\alpha+1)\sin(\alpha\pi)} \|h\|_2^{2\alpha} + O(\|h\|_2^2),
\end{align*}
where $\mathcal{I}_\alpha(\cdot)$ is the modified Bessel function of the first kind. Therefore,
$$K(x,x)-K(x,x+h) = \frac{\pi}{2^{2\alpha}\Gamma(\alpha)\Gamma(\alpha+1)\sin(\alpha\pi)} \|h\|_2^{2\alpha} + O(\|h\|_2^2),$$
verifying the upper bound (\ref{eq:matern diff lipschitz}) for $\alpha\not\in\mathbb{N}$.
For $\alpha\in\mathbb{N}$, combining equations (10.25.2) and (10.31.1) of \cite{DLMF} yields
\begin{align*}
K(x,x+h) = & \|h\|_2 \mathcal{K}_1(\|h\|_2) \\
= & \|h\|_2 \big[\|h\|_2^{-1} + \ln(\|h\|_2/2) \mathcal{I}_1(\|h\|_2) \big] + O(\|h\|_2^2) \\
= & 1 + \frac{1}{2} \|h\|_2^2 \ln(\|h\|_2/2) + O(\|h\|_2^2).
\end{align*}
Therefore, 
$$K(x,x)-K(x,x+h) = -\frac{1}{2}\|h\|_2^2 \ln(\|h\|_2) + O(\|h\|_2^2),$$
concluding the proof of the lemma.
\end{proof}

\subsection{Proof of Lemma \ref{lem:matern cond smooth}}\label{sec:alpha int |}
\begin{proof}
For all $|k_1|+|k_2|\le 2\underline{\alpha}$, the $(k_1,k_2)$th derivative of $\tilde{K}$ is
$$\tilde{K}^{(k_1,k_2)}(x_1,x_2) = K^{(k_1,k_2)}(x_1,x_2) - \frac{\partial^{k_1}}{\partial x_1^{k_1}}K_{x_1,A} K_{A,A}^{-1} \frac{\partial^{k_2}}{\partial x_2^{k_2}}K_{A,x_2}.$$
Thus $\tilde{K}(\cdot,\cdot)$ is $2\underline{\alpha}$ times differentiable on $\mathbb{R}^{2d}$.

Next, using the notation  $\inner{A,\xi} = (\inner{w_1,\xi},\inner{w_2,\xi},\ldots,\inner{w_m,\xi})^T$ and $A-x=(w_1-x,...,w_m-x)$, one can note that
\begin{align}
\tilde{K}_{x_1,x_2}&= K_{x_1,x_2} - K_{x_1,A} K_{A,A}^{-1} K_{A,x_2}\nonumber\\
&= K_{x_1,x_2} - K_{x_1,A} K_{A,A}^{-1} K_{A,x_2} - K_{x_2,A} K_{A,A}^{-1} K_{x_1,A}^T + K_{x_1,A} K_{A,A}^{-1} K_{A,A} K_{A,A}^{-1} K_{A,x_2} \nonumber\\
&= \int_{\mathbb{R}^d} \Big[
\expinner{x_1-x_2} - K_{x_1,A} K_{A,A}^{-1}\expinner{A-x_2}- K_{x_2,A}K_{A,A}^{-1}\expinner{x_1-A}  \nonumber\\
& \qquad\qquad + K_{x_1,A}K_{A,A}^{-1}\expinner{A} (\expinner{-A})^T K_{A,A}^{-1}K_{A,x_2}
\Big] (1+\|\xi\|^2_2)^{-(\alpha+d/2)} d\xi \nonumber\\
&=\int_{\mathbb{R}^d} \big(e^{-\iota\inner{x_1,\xi}}-K_{x_1,A}K_{A,A}^{-1}e^{-\iota\inner{A,\xi}} \big) \nonumber\\
&\qquad\qquad\times \overline{ \big(e^{-\iota\inner{x_2,\xi}}-K_{x_2,A}K_{A,A}^{-1}e^{-\iota\inner{A,\xi}} \big) } (1+\|\xi\|^2_2)^{-(\alpha+d/2)} d\xi.\label{eq:cond alter}
\end{align}

Define the function $x\mapsto g(x;\xi)$ as
\begin{equation}\label{eq:g}
g(x;\xi) = e^{-\iota \inner{x,\xi}}-K_{x,A}K_{A,A}^{-1}e^{-\iota\inner{A,\xi}}.
\end{equation}
Noting that $K_{x,A}K_{A,A}^{-1}e^{-\iota\inner{A,\xi}}$ is the Gaussian interpolation passing through $\{(x, e^{-\iota \inner{x,\xi}}), \forall x\in A\}$, the function $g(x;\xi)$ is $2\underline{\alpha}$ times differentiable and $g(x;\xi) = 0,\forall x\in A$. Then, for each multi-index $k$, we denote its $k$th derivative as
\begin{equation*}
g^{(k)}(x;\xi) = (-\iota)^{|k|}\xi^k\expinner{x} - \frac{\partial^k}{\partial x^{k}}K_{x,A} K_{A,A}^{-1} \expinner{A}.
\end{equation*}
The derivatives of $\tilde{K}(x_1,x_2)$ are then given as
\begin{equation}\label{eq:K tilde deriv}
\tilde{K}^{(k_1,k_2)}(x_1,x_2) = \int_{\mathbb{R}^d} g^{(k_1)}(x_1;\xi) \overline{g^{(k_2)}(x_2;\xi)} (1+\|\xi\|^2_2)^{-(\alpha+d/2)} d\xi.
\end{equation}
When $x_1=x_2=x$ and $k_1=k_2=k$, equation (\ref{eq:K tilde deriv}) yields
$$\tilde{K}^{(k,k)}(x,x) =  \int_{\mathbb{R}^d} |g^{(k)}(x;\xi)|^2(1+\|\xi\|^2_2)^{-(\alpha+d/2)}  d\xi\ge 0.$$
Furthermore, since $K_{A,A}^{-1}$ is a positive definite matrix, we also have $\tilde{K}^{(k,k)}(x,x)\le K^{(k,k)}(x,x)$, implying
$$0\leq \tilde{K}^{(k,k)}(x,x)\leq K^{(k,k)}(x,x)\lesssim 1.$$
In addition, following from (\ref{eq:K tilde deriv}),
\begin{align*}
|\tilde{K}^{(k,k)}(x_1,x_2)| 
\le & \int_{\mathbb{R}^d} \big(|g^{(k)}(x_1;\xi)|^2+ |g^{(k)}(x_2;\xi)|^2\big)(1+\|\xi\|^2_2)^{-(\alpha+d/2)} d\xi \\
\le &\tilde{K}^{(k,k)}(x_1,x_1)+ \tilde{K}^{(k,k)}(x_2,x_2) \lesssim 1,
\end{align*}
providing (\ref{eq:matern cond diff}).

Next, we deal with (\ref{eq:matern cond diff lipschitz}). For all $|k|=\underline{\alpha}$ and $h\in\mathbb{R}^d$ in a neighborhood of zero, let us define an auxiliary process
$$Y_x^h = \frac{\partial^k}{\partial x^k}(\tilde Z_{x+h} - \tilde Z_x). $$
By Proposition I.3 of \cite{ghosal2017fundamentals}, $Y_x^h$ is still a Gaussian process with covariance function
\begin{align}\label{eq:cov Y}
\cov(Y_{x_1}^h, Y_{x_2}^h) = & \frac{\partial^{k}}{\partial x_1^k} \frac{\partial^{k}}{\partial x_2^k}\cov(\tilde Z_{x_1+h} - \tilde Z_{x_1}, \tilde Z_{x_2+h} - \tilde Z_{x_2}) \nonumber\\
= & \tilde{K}^{(k,k)}(x_1,x_2)-\tilde{K}^{(k,k)}(x_1,x_2+h) \nonumber\\
&\qquad\qquad-\tilde{K}^{(k,k)}(x_1+h,x_2)+\tilde{K}^{(k,k)}(x_1+h,x_2+h).
\end{align}
Thus, it is sufficient to study the covariance function of the process $Y_x^h$. For $x_1=x_2=x$, we have
\begin{align*}
\var(Y_x^h) = & 
\big[K^{(k,k)}(x+h,x+h)-K^{(k,k)}(x,x+h)-K^{(k,k)}(x+h,x)+K^{(k,k)}(x,x)\big] \\
& - \frac{\partial^k}{\partial x^k} (K_{x+h,A}-K_{x,A}) K_{A,A}^{-1} \frac{\partial^k}{\partial x^k} (K_{A,x+h}-K_{A,x}).
\end{align*}
Since $K_{A,A}^{-1}$ is a positive definite matrix, we have
\begin{align*}
0\le \var(Y_x^h) \le & \big[K^{(k,k)}(x+h,x+h)-K^{(k,k)}(x,x+h)-K^{(k,k)}(x+h,x)+K^{(k,k)}(x,x)\big] \\
= & 2 \big[K^{(k,k)}(x,x)-K^{(k,k)}(x,x+h)\big],
\end{align*}
where the second inequality utilizes equation (\ref{eq:matern diff symmetry}). Since $Y_x^h$ is a Gaussian process, we have
\begin{equation}\label{eq:cov Y 2}
\cov(Y_{x_1}^h, Y_{x_2}^h)\le \max\{\var(Y_{x_1}^h), \var(Y_{x_2}^h) \} \le 2 \sup_{x\in\mathbb{R}^d}\big[K^{(k,k)}(x,x)-K^{(k,k)}(x,x+h)\big].
\end{equation}
The combination of (\ref{eq:matern diff lipschitz}), (\ref{eq:cov Y}) and (\ref{eq:cov Y 2}) provides our statement (\ref{eq:matern cond diff lipschitz}).
\end{proof}

\subsection{Proof of Theorem \ref{thm:var}}
\begin{proof}
The proof is organized in three steps. In the first and second steps, we fix the rescaling parameters $(\tau,s)=(1,1)$ and derive matching lower and upper bounds. In the third step, we apply the rescaling parameters to complete the proof.\\

\textbf{Direction $\gtrsim$:}\mbox{}
The case $i=1$ is straightforward. For all $i\ge 2$, we obtain
$$\var \big\{Z_{X_{i}}-\mathbb{E}[{Z}_{X_{i}}|{Z}_{\mr{pa}(X_{i})}]\big\} =  K_{X_i,X_i} - K_{\mr{pa}(X_i),X_i}^T K_{\mr{pa}(X_i),\mr{pa}(X_i)}^{-1} K_{\mr{pa}(X_i),X_i}.$$
Define the set $A_i:= \{X_i\}\cup \mr{pa}(X_i)$. Then, applying Theorem 12.3 of \cite{wendland2004scattered} (with $\Phi$ being the Mat\'{e}rn covariance function, $\varphi_0(M) \asymp \inf_{\|\xi\|_2\le 2M} (1+\|\xi\|_2^2)^{-\alpha+d/2} = M^{-(2\alpha+d)}$, $q_X$ being $\gamma^{-j}$ and $M \asymp d/q_X$),
the minimal eigenvalue of the matrix $K_{A_i,A_i}$ satisfies that
\begin{equation*}
\lambda_{\min}(K_{A_i,A_i}) \gtrsim M^{-(2\alpha+d)} M^d = M^{-2\alpha} \asymp \gamma^{-2\alpha j}.
\end{equation*}
Therefore, 
\begin{align*}
& 1-K_{\mr{pa}(X_i),X_i}^T K_{\mr{pa}(X_i),\mr{pa}(X_i)}^{-1} K_{\mr{pa}(X_i),X_i} \\
= & \begin{bmatrix} 1 \\ - K_{\mr{pa}(X_i),\mr{pa}(X_i)}^{-1} K_{\mr{pa}(X_i),X_i} \end{bmatrix}^T 
\begin{bmatrix}
1 & K_{\mr{pa}(X_i),X_i}^T \\ K_{\mr{pa}(X_i),X_i} & K_{\mr{pa}(X_i),\mr{pa}(X_i)}
\end{bmatrix}
\begin{bmatrix} 1 \\ - K_{\mr{pa}(X_i),\mr{pa}(X_i)}^{-1} K_{\mr{pa}(X_i),X_i} \end{bmatrix} \\
\ge & \left\| \begin{bmatrix} 1 \\ - K_{\mr{pa}(X_i),\mr{pa}(X_i)}^{-1} K_{\mr{pa}(X_i),X_i} \end{bmatrix} \right\|^2 \lambda_{\min}\left( K_{A_i,A_i}\right) \\
\ge & 1 \cdot \lambda_{\min}\left( K_{A_i,A_i}\right) \gtrsim \gamma^{-2\alpha j}.
\end{align*}

\textbf{Direction $\lesssim$:}\mbox{}
Take $A=\mr{pa}(X_i)$ in Lemma \ref{lem:matern cond smooth} and recall the formula for the derivative $\tilde{K}^{(k,k)}(x_1,x_2)$ in equation (\ref{eq:K tilde deriv}) and the definition of $g(x;\xi)$  in (\ref{eq:g}). Since $g$ is $2\underline{\alpha}$ times differentiable, recalling the definition of the cube $\C$ in Condition \ref{cond:norming}, the $\underline{\alpha}$ order Taylor expansion of $g(x;\xi)$ around $x_0\in \C$  can be expressed as
\begin{equation*}
g(x;\xi) = Q(x;\xi) + R(x;\xi),    
\end{equation*}
where $Q(x;\xi)$ is a polynomial of order (up to) $\underline{\alpha}$ and $x\rightarrow R(x;\xi)$ is a $\underline{\alpha}$ times differentiable function with
$R^{(k)}(x_0;\xi) = 0, \;\forall |k|\le\underline{\alpha}.$
Thus, for any multi-index $k$ satisfying $|k|=\underline{\alpha}$, the $k$th derivative of $Q(x;\xi)$ is constant. Denoting this constant by $Q_k$, we get
\begin{equation}\label{eq:h deriv}
g^{(k)}(x;\xi) = Q_k + R^{(k)}(x;\xi).
\end{equation}
By plugging in equation (\ref{eq:h deriv}) into equation (\ref{eq:K tilde deriv}) and the latter one further into the inequality (\ref{eq:matern cond diff lipschitz}), we get for all $\alpha\not\in\mathbb{N}$ and $|k|=\underline{\alpha}$,
\begin{align*}
\|h\|_2^{2(\alpha-\underline{\alpha})} 
\gtrsim & \tilde{K}^{(k,k)}(x_0,x_0) - \tilde{K}^{(k,k)}(x_0,x_0+h) - \tilde{K}^{(k,k)}(x_0+h,x_0) + \tilde{K}^{(k,k)}(x_0+h,x_0+h) \\
= & \int_{\mathbb{R}^d} \big[ g^{(k)}(x_0;\xi) - g^{(k)}(x_0+h;\xi) \big] \overline{\big[ g^{(k)}(x_0;\xi) - g^{(k)}(x_0+h;\xi) \big]} (1+\|\xi\|_2)^{-(\alpha+d/2)} d\xi \\
= & \int_{\mathbb{R}^d} \big[ R^{(k)}(x_0;\xi) - R^{(k)}(x_0+h;\xi) \big] \overline{\big[ R^{(k)}(x_0;\xi) - R^{(k)}(x_0+h;\xi) \big]} (1+\|\xi\|_2)^{-(\alpha+d/2)} d\xi \\
\overset{R^{(k)}(x_0;\xi)=0}{=} & \int_{\mathbb{R}^d} |R^{(k)}(x_0+h;\xi)|^2 (1+\|\xi\|_2)^{-(\alpha+d/2)} d\xi.
\end{align*}
For all $x\in \C$, let $x=x_0+h$, then $\|h\|\lesssim\gamma^{-j}$ and in view of the previous display,
\begin{equation}\label{eq:R x}
\int_{\mathbb{R}^d} |R^{(k)}(x;\xi)|^2 (1+\|\xi\|_2)^{-(\alpha+d/2)} d\xi \lesssim \gamma^{-2(\alpha-\underline{\alpha})j}.
\end{equation}
We prove below that the above inequality can be extended for multi-indices satisfying $|k|\le\underline{\alpha}$, i.e.
\begin{equation}\label{eq:R k x}
\int_{\mathbb{R}^d} |R^{(k)}(x;\xi)|^2 (1+\|\xi\|_2)^{-(\alpha+d/2)} d\xi \lesssim \gamma^{-2(\alpha-|k|)j}.
\end{equation}

Applying the above inequality for $|k|=0$, we have
\begin{equation}\label{eq:R int bound}
\int_{\mathbb{R}^d} |R(x;\xi)|^2(1+\|\xi\|_2^2)^{-(\alpha+d/2)} d\xi \lesssim \gamma^{-2\alpha j}, \;\forall x\in \C.
\end{equation}
Recalling the definition of $g$ in \eqref{eq:g} and the fact that $g(x;\xi)=0$ for all $x\in A \supset\mr{pa}(X_i)$, we have
$|Q(x;\xi)|=|R(x;\xi)| $. Furthermore, since  $\mr{pa}(X_i)$ is a norming set on the domain $\C$ for the collection of polynomials  up to degree $\underline{\alpha}$,
\begin{equation*}
\sup_{x\in \C} |Q(x;\xi)|\lesssim \sup_{x\in\mr{pa}(X_i)} |Q(x;\xi)|= \sup_{x\in\mr{pa}(X_i)} |R(x;\xi)|.
\end{equation*}
This in turn implies that for all $ x\in \C$,
\begin{equation*}
|g(x;\xi)| \lesssim \sup_{t\in\mr{pa}(X_i)} |R(t;\xi)| + |R(x;\xi)|.
\end{equation*}
Then by plugging in the preceding upper bound into  (\ref{eq:K tilde deriv})  and using \eqref{eq:R int bound}, we arrive at
\begin{align*}
\var \big\{Z_{X_{i}}-\mathbb{E}[{Z}_{X_{i}}|{Z}_{\mr{pa}(X_{i})}]\big\} 
= & \int |g(X_i;\xi)|^2 (1+\|\xi\|_2^2)^{-(\alpha+d/2)} d\xi \\
\lesssim & \int \big[\sup_{t\in\mr{pa}(X_i)} |R(t;\xi)| + |R(X_i;\xi)|\big]^2 (1+\|\xi\|_2^2)^{-(\alpha+d/2)} d\xi \\
\le & 2 \int \big[\sum_{t\in\mr{pa}(X_i)} |R(t;\xi)|^2 + |R(X_i;\xi)|^2\big] (1+\|\xi\|_2^2)^{-(\alpha+d/2)} d\xi \\
\lesssim & \gamma^{-2\alpha j}, 
\end{align*}
proving the upper bound in the statement for $s=\tau=1$.

Therefore, it only remains to prove assertion (\ref{eq:R k x}). We proceed by induction on the decreasing order of $|k|$. The case $|k|=\underline{\alpha}$ is covered in (\ref{eq:R x}). Now suppose (\ref{eq:R k x}) holds for all $k$ satisfying $|k|>l\in\mathbb{N}$. Denoting by  $e_j\in\mathbb{N}^d$,  $e_j[i] = \mathbbm{1}_{\{j=i\}}$ the elementary basis vector, we have for $|k|=l$,
\begin{small}
\begin{align*}
& \int_{\mathbb{R}^d}|R^{(k)}(x;\xi)|^2 (1+\|\xi\|_2)^{-(\alpha+d/2)} d\xi \\
= & \int_{\mathbb{R}^d} |R^{(k)}(x;\xi)-R^{(k)}(x_0;\xi)|^2 (1+\|\xi\|_2)^{-(\alpha+d/2)} d\xi \\
= & \int_{\mathbb{R}^d} \bigg|\sum_{j=1}^d \big[R^{(k)}((x_0[1],\ldots,x_0[j-1],x[j],\ldots,x[d])^T;\xi)\\
& \quad -R^{(k)}((x_0[1],\ldots,x_0[j],x[j+1],\ldots,x[d])^T;\xi)\big] \bigg|^2 (1+\|\xi\|_2^2)^{-(\alpha+d/2)} d\xi \\
= & \int_{\mathbb{R}^d} \bigg|\sum_{j=1}^d \int_{x_0[j]}^{x[j]} R^{(k+e_j)}((x_0[1],\ldots,x_0[j-1],t, x[j+1],\ldots,x[d])^T;\xi) dt \bigg|^2 (1+\|\xi\|_2^2)^{-(\alpha+d/2)} d\xi \\
\le & d \sum_{j=1}^d \int_{\mathbb{R}^d} \big|x[j]-x_0[j]\big|  \bigg[\int_{x_0[j]}^{x[j]} \big| R^{(k+e_j)}((x_0[1],\ldots,x_0[j-1],t, x[j+1],\ldots,x[d])^T;\xi) \big|^2 dt \bigg]\\
&\qquad\qquad \times (1+\|\xi\|_2^2)^{-(\alpha+d/2)} d\xi \\
\lesssim & d\gamma^{-j} \sum_{j=1}^d \int_{x_0[j]}^{x[j]} \gamma^{-2(\alpha-(l+1))j} dt
\lesssim  \gamma^{-2(\alpha-l)j},
\end{align*}
\end{small}
\noindent where in the first inequality we used Cauchy-Schwarz and in the second inequality we used Fubini's theorem together with the induction assumption. This concludes the proof of  (\ref{eq:R k x}) and the upper bound for the statement for $s=\tau=1$.\\

\textbf{Rescaling:}
For general scaling parameters $\tau$ and $s$, compared to the case $\tau=s=1$, the domain is mapped to $\tau\X$ and the covariance matrix is scaled by a factor of $s^2$. Therefore, by replacing $\gamma^{-j}$ with $\tau \gamma^{-j}$ and multiplying the final term with $s^2$, we obtain
$$\var \big\{Z^{\tau,s}_{X_{i}}-\mathbb{E}[{Z}^{\tau,s}_{X_{i}}|{Z}^{\tau,s}_{\mr{pa}(X_{i})}]\big\} \asymp s^2 (\tau\gamma^{-j})^{2\alpha } = s^2 \tau^{2\alpha} \gamma^{-2\alpha j},
$$
concluding the proof.
\end{proof}

\subsection{Proof of Lemma \ref{lem:GP flat}}
\begin{proof}
\textbf{Proof of (\ref{eq:GP flat 2}):} For $\alpha\le 1$, in view of (\ref{eq:matern diff lipschitz}), noticing $\|\nu x' - \nu_0 x'\|_2 \le \sqrt{d} |\nu-\nu_0|$ for all $x'\in A\cup\{x_0\}$,
\begin{equation}\label{eq:tau waaa2}
\|K_{\nu A, \nu A} - K_{\nu_0 A, \nu_0 A}\|_1 \lesssim |\nu - \nu_0|^{2\alpha}, \quad \|K_{\nu A, \nu x} - K_{\nu_0 A, \nu_0 x}\|_1 \lesssim |\nu - \nu_0|^{2\alpha}.
\end{equation}
Since $\nu_0$ is a non-zero constant, the minimal and maximal eigenvalues of the matrix $K_{\nu_0 A, \nu_0 A}$ are bounded away from zero and infinity. For $\tau$ sufficient close to $\nu_0$,  the minimal and maximal eigenvalues of the matrix $K_{\nu A, \nu A}$ are also bounded away from zero and infinity. Therefore, we have 
$$\|K_{\nu_0 A, \nu_0 A}^{-1}\|_1\lesssim 1, \quad \|K_{\nu A, \nu A}^{-1}\|_1 \lesssim 1$$
and subsequently
\begin{equation}\label{eq:tau waaa1}
\|K_{\nu A, \nu A}^{-1} - K_{\nu_0 A, \nu_0 A}^{-1}\|_1 \le  
\|K_{\nu_0 A, \nu_0 A}^{-1}\|_1 \|K_{\nu A, \nu A} - K_{\nu_0 A, \nu_0 A}\|_1 \|K_{\nu A, \nu A}^{-1}\|_1 \lesssim |\nu-\nu_0|^{2\alpha}.
\end{equation}
Combining the inequalities (\ref{eq:tau waaa2}) and (\ref{eq:tau waaa1}) results in
\begin{align*}
& \|K_{\nu A,\nu A}^{-1}K_{\nu A,\nu x} - K_{\nu_0 A,\nu_0 A}^{-1}K_{\nu_0 A,\nu_0 x}\|_1 \\
= & \|(K_{\nu A,\nu A}^{-1}-K_{\nu_0 A,\nu_0 A}^{-1})K_{\nu A,\nu x} +  K_{\nu_0 A,\nu_0 A}^{-1}(K_{\nu A,\nu x}-K_{\nu_0 A,\nu_0 x})\|_1 \\
\lesssim & |\nu-\nu_0|^{2\alpha} (\|K_{\nu A,\nu x}\|_1 + \|K_{\nu_0 A,\nu_0 A}^{-1}\|_1) \lesssim |\nu-\nu_0|^{2\alpha},
\end{align*}
 providing (\ref{eq:GP flat 2}) for the case $\alpha\le 1$. 

For $\alpha>1$, in view of Lemma \ref{lem:matern smooth}, $K(\cdot,\cdot)$ is at least twice differentiable with uniformly bounded derivatives, hence,
$$\|K_{\nu A, \nu A} - K_{\nu_0 A, \nu_0 A}\|_1 \lesssim |\nu - \nu_0|^{2}, \quad \|K_{\nu A, \nu x} - K_{\nu_0 A, \nu_0 x}\|_1 \lesssim |\nu - \nu_0|^{2}.
$$
With exactly the same derivation as the case $\alpha\le 1$, we obtain
$$\|K_{\nu A,\nu A}^{-1}K_{\nu A,\nu x} - K_{\nu_0 A,\nu_0 A}^{-1}K_{\nu_0 A,\nu_0 x}\|_1 \lesssim |\nu-\nu_0|^2 ,$$
concluding the proof of (\ref{eq:GP flat 2}).\\

\textbf{Proof of (\ref{eq:GP flat}):}
For $\delta\geq0$, let us write $M_\nu=O(\nu^{\delta})$ for a matrix $M_\nu$ indexed by a parameter $\nu$, if $\|M_\nu\|_1\lesssim \nu^{\delta}$ as $\tau$ tends to 0. We recall that for finite dimensional matrices, the matrix $L_1$ norm is equivalent to various other matrix norms, including $L_2$, $L_\infty$ and trace norms.

Then, for all $1\le j \le m$, define the set $A_{j,x}$ as
$$A_{j,x} = \{w_1,w_2,\cdots, w_{j-1},x,w_{j+1},\cdots, w_m\},$$
i.e. replacing the element $w_j$ with $x$. Then, since the matrix $K_{\nu A,\nu A}$ is  positive definite for all $\nu>0,$ we get by Cramer's rule that
\begin{align}
 \label{eq:cramer 1}
\big(K_{\nu A,\nu A}^{-1} K_{\nu A,\nu x}\big)[j] = & \frac{\det(K_{\nu A,\nu A_{j,x}})}{\det(K_{\nu A,\nu A})}.
\end{align}
Furthermore, let us introduce the matrix
\begin{align}
 B_{j,x}:=V_A^{-1} V_{A_{j,x}}.\label{def:mtx:b}
\end{align}
Note that except for the $j$th column, the matrix $B_{j,x}$ coincides with the identity matrix, while its $j$th column is $V_A^{-1} v_x$. Therefore,
$$\det(B_{j,x}) = B_{j,x}[j,j] = (V_A^{-1} v_x)[j].$$
Furthermore, we show below that 
\begin{align}
\frac{\det(K_{\nu A,\nu A_{j,x}})}{\det(K_{\nu A,\nu A})} 
= \det(B_{j,x}) + O(\nu + \nu^{2(\alpha-\underline{\alpha})}).\label{eq:help:lem7}
\end{align}
The combination of the previous four displays implies that
\begin{align*}
\big(K_{\nu A,\nu A}^{-1} K_{\nu A,\nu x}\big)[j]= (V_A^{-1} v_x)[j] + O(\nu + \nu^{2(\alpha-\underline{\alpha})}),\;\forall 1\le j\le m,
\end{align*}
providing the statement of the lemma.

Hence, it only remains to show that \eqref{eq:help:lem7} holds. We proceed by reformulating the matrices $K_{\tau A,\tau A_{j,x}}$ and $K_{\tau A,\tau A}$ into a more suitable, product form, by taking their Taylor expansions. We also introduce some additional notation to facilitate the arguments.

First, we consider the case when the origin $0$ belongs to the convex hull of $A\cup\{x\}$. In view of Lemma \ref{lem:matern smooth}, the $2\underline\alpha$ order Taylor expansion of $K(\cdot,\cdot)$ around $0$ can be written as
\begin{equation}\label{eq:K taylor}
K(x_1,x_2) = \sum_{|k_1|+|k_2|\le 2\underline{\alpha}} q_{k_1,k_2}x_1^{k_1}x_2^{k_2} + O\big((\|x_1\|_2+\|x_2\|_2)^{2\alpha}\big),
\end{equation}
where $q_{k_1,k_2} = K^{(k_1,k_2)}(0,0)/(k_1!k_2!)$ is the multivariate Taylor expansion coefficient. 
Let $m' = {2\underline{\alpha}+d \choose 2\underline{\alpha}}-{\underline{\alpha}+d \choose \underline{\alpha}}$ and
define the matrices $Q_{11}\in\mathbb{R}^{m\times m}$, $Q_{12}\in\mathbb{R}^{m\times m'}$ as
$$
Q_{11}[i,j] = q_{k_{(i)},k_{(j)}},\;\; Q_{12}[i,j] = \bigg\{\begin{aligned}
&q_{k_{(i)},k_{(j+m)}},   && |k_{(i)}|+|k_{(j+m)}|\le 2\underline{\alpha},\\
& 0, && |k_{(i)}|+|k_{(j+m)}|> 2\underline{\alpha},
\end{aligned}
$$
where $k_{(i)}\in\mathbb{N}^d$ denotes the $i$th multi-index in lexicographical ordering. Furthermore, let us define the matrix $Q\in\mathbb{R}^{(m+m')\times(m+m')}$ in blocks as
$$Q = \begin{bmatrix}
Q_{11}  & Q_{12} \\ Q_{12}^T & 0
\end{bmatrix}.$$
By the symmetry of $K(\cdot,\cdot)$, note $Q[i,j]=Q[j,i]$ and hence the matrix $Q$ satisfies
$$Q[i,j] = \bigg\{\begin{aligned}
&q_{k_{(i)},k_{(j)}},   && |k_{(i)}|+|k_{(j)}|\le 2\underline{\alpha},\\
& 0, && |k_{(i)}|+|k_{(j)}|> 2\underline{\alpha}.
\end{aligned}$$
Finally,  define the matrix $\tilde{V}_A\in\mathbb{R}^{m'\times m}$ as $\tilde V_A[i,j] = w_j^{k_{(m+i)}}$, where $w_j$s are the elements of $A$, and the
diagonal matrices $\Gamma_1\in\mathbb{R}^{m\times m}$, $\Gamma_2\in\mathbb{R}^{m'\times m'}$ and $\Gamma\in\mathbb{R}^{(m+m')\times(m+m')}$ as
$$
\Gamma_1[i,i] = \nu^{|k_{(i)}|}, \;\; \Gamma_2[i,i] =\nu^{|k_{(m+i)}|},\;\; 
\Gamma = \begin{bmatrix}
 \Gamma_1 & 0 \\  0 & \Gamma_2 
\end{bmatrix},
$$
respectively. Since $k_{(i)}$ is the $i$th multi-index in lexicographical ordering and $m = {\underline{\alpha}+d \choose \underline{\alpha}}$, we have
$$\|\Gamma_1^{-1}\|_2 \lesssim \nu^{-\underline{\alpha}}, \;\; \|\Gamma_2\|_2 \lesssim \nu^{\underline{\alpha}+1}.$$
Taking the Taylor expansion (\ref{eq:K taylor}) for each element of the matrix $K_{\tau A,\tau A}$ yields
\begin{align}\label{eq:Ktau decom}
K_{\nu A,\nu A} 
= & [V_A^T \; {\tilde{V}_A}^T] 
\begin{bmatrix}
\Gamma_1 & 0 \\  0 & \Gamma_2   
\end{bmatrix} 
\begin{bmatrix}
Q_{11}  & Q_{12} \\ Q_{12}^T & 0
\end{bmatrix}
\begin{bmatrix}
\Gamma_1 & 0 \\  0 & \Gamma_2   
\end{bmatrix}
\begin{bmatrix}
V_A \\ \tilde{V}_A
\end{bmatrix} +O(\nu^{2\alpha})\nonumber\\
= & V_A^T \Gamma_1 \big[ Q_{11} + Q_{12}\Gamma_2 \tilde{V}_A V_A^{-1}\Gamma_1^{-1} + \Gamma_1^{-1} V_A^{-T}  {\tilde{V}_A}^T \Gamma_2 Q_{12}^T \big]  \Gamma_1 V_A+  O(\nu^{2\alpha}) \nonumber\\
= & V_A^T \Gamma_1 \tilde{Q}_{11} \Gamma_1 V_A \big[ I + O(\nu^{2(\alpha-\underline{\alpha})})\big] ,
\end{align}
where $\tilde{Q}_{11} = Q_{11} + Q_{12}\Gamma_2 \tilde{V}_A V_A^{-1}\Gamma_1^{-1} + \Gamma_1^{-1} V_A^{-T}  {\tilde{V}_A}^T \Gamma_2 Q_{12}^T = Q_{11} + O(\nu).$

Next we proceed with reformulating the covariance matrix $K_{\tau A,\tau A_{j,x}}$ in a similar manner. First note, that in view of Lemma \ref{lem:norming def} and its proof thereafter, for a norming set $A$  with norming constant $c_N$, for all $P\in\mathscr{P}_l([0,1]^d)$, 
$$P(x) = P(A)^T V_A^{-1} v_x \le c_N \|P(A)\|_\infty,\qquad \forall x\in[0,1]^d.$$
Hence, for any $x\in[0,1]^d$, by choosing $P(A)[i] = \mr{sign}( (V_A^{-1} v_x) [i]), \forall 1\le i\le d$, we have
\begin{equation}
\|B_{j,x}[\cdot,j]\|_1=\|V_A^{-1} v_x\|_1 \le c_N.\label{eq: UB:Bjx}
\end{equation}
Next, let us define the matrix $\tilde V_{A_{j,x}}$ similarly as $\tilde V_A$. Then we can write
$$\begin{bmatrix} V_{A_{j,x}} \\ \tilde{V}_{A_{j,x}}\end{bmatrix} = \begin{bmatrix} V_{A} \\ \tilde{V}_{A}\end{bmatrix} B_{j,x} + 
\begin{bmatrix}
0 \\ R_{j,x}
\end{bmatrix},$$
where $R_{j,x} = \tilde{V}_{A_{j,x}} - \tilde{V}_{A} B_{j,x}$. Since $A\cup\{x\}\subset [0,1]^d$, in view of \eqref{eq: UB:Bjx}, the error term $R_{j,x}$ has bounded $L_1$ norm independent of the set $A$ and point $x$.  Then, the covariance matrix $K_{\tau A, \tau A_{j,x}}$ can be written as
\begin{align}\label{eq:Ktau x decom}
K_{\tau A, \tau A_{j,x}} 
= & [V_A^T \; {\tilde{V}_A}^T] 
\begin{bmatrix}
\Gamma_1 & 0 \\  0 & \Gamma_2   
\end{bmatrix}
\begin{bmatrix}
Q_{11}  & Q_{12} \\ Q_{12}^T & 0
\end{bmatrix}
\begin{bmatrix}
\Gamma_1 & 0 \\  0 & \Gamma_2   
\end{bmatrix}
\left(\begin{bmatrix}
V_{A} \\ \tilde{V}_{A}
\end{bmatrix} B_{j,x} + 
\begin{bmatrix}
0 \\ R_{j,x}
\end{bmatrix}\right) + O(\tau^{2\alpha})\nonumber\\
= & V_A^T \Gamma_1 \tilde{Q}_{11} \Gamma_1 V_A B_{j,x} +
V_A^T \Gamma_1 {Q}_{12} \Gamma_2 R_{j,x} + O(\tau^{2\alpha}) \nonumber\\
= & V_A^T \Gamma_1 \tilde{Q}_{11} \Gamma_1 V_A \big[ B_{j,x} + V_A^{-1} \Gamma_1^{-1} \tilde{Q}_{11}^{-1} Q_{12} \Gamma_2 R_{j,x} \big] + O(\tau^{2\alpha}) \nonumber\\
= & V_A^T \Gamma_1 \tilde{Q}_{11} \Gamma_1 V_A \big[ B_{j,x} + O(\tau) + O(\tau^{2(\alpha-\underline{\alpha})}) \big].
\end{align}
Since the determinant is a multiplicative map, we obtain, in view of the previous display and \eqref{eq:Ktau decom} that
\begin{align*}
\frac{\det(K_{\tau A,\tau A_{j,x}})}{\det(K_{\tau A,\tau A})} =& \frac{\det(B_{j,x} + O(\tau) + O(\tau^{2(\alpha-\underline{\alpha})}) )}{\det(I+O(\tau^{2(\alpha-\underline{\alpha})}))}
=  \det(B_{j,x}) + O(\tau + \tau^{2(\alpha-\underline{\alpha})}),    
\end{align*}
concluding the proof of \eqref{eq:help:lem7}.

Finally, in case the origin doesn't belong to the convex hull of $A\cup\{x\}$, take an arbitrary  $x_0\in A\cup\{x\}$, and note that $0\in A\cup\{x\}-x_0$. Hence, by working with the shifted process $(Z_{x-x_0}^{\tau,s}|Z_{A-x_0}^{\tau,s})$, the conclusion holds.
\end{proof}

\subsection{Bounding $\vartheta_n$ defined in \eqref{eq:poly interp cond}}\label{sec:vartheta}
The quantity $\vartheta_n$ plays a crucial role in controlling the recursive polynomial and Gaussian interpolations, and hence in deriving small deviation bounds for Vecchia GPs. We focus on grid data and verify that, for $d\leq 2$ and moderate values of the regularity parameter $\alpha$, $\vartheta_n=O(1)$. Controlling this term in higher dimensions  or for more general design settings is highly challenging and beyond the scope of the present paper. Nevertheless, we conjecture that the results extend beyond the restricted setting considered here. To simplify the analysis, we consider a slightly modified version of Algorithm \ref{alg:dag arg} to handle data points near the boundary. We discuss this technical extension in detail the next subsection.

\begin{lemma}\label{lem:poly recursive}
Suppose Condition \ref{cond:norming} holds, $\mathcal{X}_n$ is on a grid and the DAG of the Vecchia GP is constructed according to Algorithm \ref{alg:dag arg}. If one of the following conditions is satisfied:
\begin{itemize}
    \item $d=1$ and   $\alpha\le 100$ 
    \item $d = 2$ and $\alpha\le 17$,
\end{itemize}
then $\vartheta_n=O(1)$.
\end{lemma}
We first discuss the conditions of Lemma \ref{lem:poly recursive}. Although the specified ranges for $d$ and $\alpha$ may appear somewhat restricted, they actually cover the cases commonly encountered in geostatistical applications, where Vecchia GPs are routinely employed. In such applications, the domain typically corresponds to a geological space with $d=1$ or $2$, and the regularity parameter $\alpha$ generally does not exceed $5/2$.


The main idea of the proof stems from the observation that, on an infinite grid, the geometric shapes of $\mr{pa}(X_i)$ are identical. Therefore, the polynomial interpolation operators $P_j,\forall j$ can be viewed as convolution operators, which can be further converted to pointwise multiplication via the Fourier transform. This approach greatly simplifies our analysis. However, we consider a finite grid, where the behavior near the boundary differs from that of an infinite grid.
Polynomial interpolation on a finite grid exhibits ``boundary effects'', which has to be treated separately. To address this, we use data augmentation. We expand the dataset $\X_n$ with artificially placed locations, so that any $X_i\in\X_n$ and its ancestors are sufficiently far from the boundary. Specifically, we first embed the finite grid $\X_n$ into a larger grid and then construct a layered norming DAG on the expanded grid, following a procedure analogous to Algorithm \ref{alg:dag grid}. We make this general strategy precise in the remainder of the section.

Our proof is divided into three parts. The first part describes how to extend a grid to mitigate boundary effects. The second part defines polynomial interpolation on an infinite grid and shows that it is equivalent to finite polynomial interpolation on the extended grid. In the final part, we derive
upper bounds for the polynomial interpolation on infinite grid using Fourier analysis.

\subsubsection{Extension of grid}
We first introduce some notation.
Without loss of generality, let the dataset $\X_n$ be a $n=(2^r+1)^d$ grid for some $r\in\mathbb{N}$ and choose $\gamma = 2$, i.e.
$$\X_n = \{(t_1,t_2, \cdots, t_{d})^T/2^r: t_h\in\mathbb{N}, 0\le t_h\le 2^r, \forall 1\le h \le d\}.$$
For all $j\in\mathbb{N}$, we also define the infinite grid with $2^{-j}$ separation distance as
$$\mathcal{G}_j = \{(t_1,t_2, \cdots, t_{d})^T/2^j: t_h\in\mathbb{Z}, \forall 1\le h \le d\}.$$
Then layer $j\in\{0,...,\eta(n)\}$ from Algorithm \ref{alg:dag grid} takes the form
$$\N_j = \X_n \cap (\mathcal{G}_j\backslash \cup_{j'=0}^{j-1}\mathcal{G}_{j'}) \subset \mathcal{G}_j.$$
Observe, that the above DAG structure satisfies Condition \ref{cond:layered} with $c_d=1$ and Condition \ref{cond:norming} with a constant $c_N>0$.
Next, we augment each $\N_j$ with elements from $\mathcal{G}_j$ to form new sets $\tilde{\N}_j$, ensuring that none of the elements $x\in\mathcal{X}_n$ are affected by boundary effects.  Specifically, for all $X_i\in\X_n$, let $\mathcal{A}(X_i)$ be the collection of all ancestors of $X_i$, i.e., the collection of all elements of $\X_n$ such that there exists a directed path from that element to $X_i$. Then  for all $X\in \mathcal{A}(X_i)\cap\N_j$, $j<\eta(i)$, under Condition \ref{cond:norming} and recalling that $\gamma=2$, we have
\begin{equation}\label{eq:ancestor dist}
\|X_i-X\|_\infty \le \sum_{j'=j}^{\eta(i)} c_L\gamma^{-j'} \le 2c_L\cdot 2^{-j}.
\end{equation}

Next, we describe the iterative data augmentation procedure. Recall that in Algorithm \ref{alg:dag grid} the parent sets for elements in layers $j\leq j_0$, with $j_0$ only depending on $\alpha$ and $d$, were initialized as $\mr{pa}(X_i)=\cup_{j'=0}^{j-1}\N_{j'}$. In this case we do not use data augmentation and set  $\tilde{\N}_j = \N_j$, $\forall\;, 0\le j \le j_0$. For all $j_0< j \le \eta(n)-1$, we iteratively augment each layer $\N_{j}$ to $\tilde{\N}_{j}$ as
\begin{equation}\label{eq:arg}
\tilde{\N}_{j} = [-4c_L2^{-j},1+4c_L2^{-j}]^d \cap (\mathcal{G}_j \backslash  \cup_{j'=0}^{j-1} \tilde{\N}_{j'}). 
\end{equation}
Note that by construction $ \N_j\subseteq\tilde{\N}_{j}$. Finally, for the last layer we take $\tilde{\N}_{\eta(n)} = \N_{\eta(n)}$.
The DAG on the new dataset 
$$\tilde{\X}_n = \cup_{j=0}^{\eta(n)} \tilde{\N}_j$$
is defined in the same way as in Algorithm \ref{alg:dag grid} and satisfies Condition \ref{cond:norming}. This also means that relation \eqref{eq:ancestor dist} holds for all $X_i\in \tilde{\X}$ as well.

For all $X_i\in \tilde{\X}_n$ and $j_0\le j\le\eta(i)-1$, the new layer $\tilde{\mathcal{N}}_j$ includes all grid points in $\mathcal{G}_j$ within $4c_L\cdot 2^{-j}-4c_L\cdot 2^{-\eta(i)}\geq 2c_L\cdot 2^{-j}$ distance of $X_i$. In view of (\ref{eq:ancestor dist}), all potential ancestors of $X_i$ in the $j$th layer lie within a distance of $2c_L\cdot 2^{-j}$ from $X_i$, hence the DAG structure of Algorithm \ref{alg:dag arg} eliminates all boundary effects for $X_i\in\cup_{j=j_0}^{\eta(n)}\tilde{\N}_j\subset\tilde{\X}_n$.
{
\makeatletter
\renewcommand{\algocf@name}{Algorithm}   
\renewcommand{\thealgocf}{1a}            
\begin{algorithm}[h]
\SetAlgoLined
Input grid data $\X_n$ with $n=\tilde{n}^d$ and order $l$. Set $m={l+d \choose l}$ and let $r$ be the largest integer such that $\tilde n\ge 2^r+1$.\\
Define the set $\bar{\X}_n$ as in (\ref{eq:subgrid}). Choose layers $\N_j, 1\le j\le r$ as in (\ref{eq:layer grid}), satisfying $\cup_{j=0}^r \N_j = \bar{X}_n$.\\
\eIf{$\tilde{n}>2^r+1$}{
Let $J=r+1$. Set $\N_{r+1}=\X_n\backslash\bar{\X}_n$.
}{
Let $J=r$.
}
Let $j_0=\arg\max_{j}\{2^{j-1}+1<l\}$. Set $\mr{pa}(X_1)=\emptyset$.\\
\For{$j_0\le j \le r-1$}{
Augment $\N_{j}$ to $\tilde{\N}_{j}$ according to \eqref{eq:arg}. \\
}

\For{$1\le j\le j_0$}{
For all $X\in\tilde{\N}_j$, set $\mr{pa}(X)=\cup_{j'=0}^{j-1}\tilde{\N}_{j'}$. 
}
\For{$j_0+1\le j \le r$}{
\For{$X\in\tilde{\N}_j$}{
\For{$1\le h\le d$}{
Order the projection of $\cup_{0\le j'\le j-1}\tilde{\N}_{j'}$ onto the $h$th coordinate by increasing distance to $X[h]$, denoting the resulting sequence as $B_h=(x_{s(1),h},x_{s(2),h},\cdots,x_{s(n'),h})$.
}
Set $\mr{pa}(X)$ as in equation (\ref{eq:parent grid}).
}
}
\KwResult{Return $\mr{pa}(X),\forall X\in\cup_{j=0}^r \tilde{\N}_j$.}
\caption{Construct layered norming DAGs on Grid with data Augmentation.}\label{alg:dag arg}
\end{algorithm}
}

We note that, instead of data augmentation, one could also trim the dataset to eliminate boundary effects, sacrificing only a negligible fraction of the data (vanishing as $n\rightarrow\infty$). This can be done in a manner analogous  to data augmentation and is left as an exercise for the reader.

\subsubsection{Polynomial interpolation on infinite grid}
With  boundary effects eliminated via data augmentation, we can define the polynomial interpolation operator on the infinite grid and align it with the finite-grid operators $P_j,j\in\mathbb{N}$.
\paragraph{Definition}
To begin with, define a multi-index set
\begin{equation}\label{eq:K set}
\mathbb{K} = \{k=(k[1],k[2],\cdots,k[d)]^T\in\mathbb{Z}^d: |k|\ge 1, 0\le k_j\le 1\}.
\end{equation}
For all multi-index $k=(k[1],k[2],\cdots,k[d])^T\in\mathbb{K}$, let $k/2$ be the element $$(k[1]/2,k[2]/2,\cdots,k[d]/2)^T\in\mathbb{R}^d.$$ Furthermore, let us denote by $\mr{pa}(k/2)$ the corner parent set defined in (\ref{eq:parent grid}) for $k/2\in \N_j$ and $\cup_{0\le j'\le j-1}\N_{j'}$ a sufficiently large subset of $\mathbb{Z}^d$.
Recalling the definition of the Vandermonde matrix $V_A$ from \eqref{def:vandemonde}, we define a linear operator $\mr{In}_k$ (``$\mr{In}$'' taking the name from interpolation) on functions $f:\mathbb{Z}^d\rightarrow \mathbb{R}$ as 
\begin{equation}
\mr{In}_k(f)(z) =  v_{k/2}^T V_{\mr{pa}(k/2)}^{-1} \bm{f}(z+\mr{pa}(k/2)),\qquad \forall z\in\mathbb{Z}^d,\label{def:Lk}
\end{equation}
where $ \bm{f}(z+\mr{pa}(k/2))\in\mathbb{R}^{|\mr{pa}(k/2)|}$. The linear operator $P$ is then defined as
\begin{equation}\label{eq:fourier combine}
P(f)(z) = \Bigg\{ 
\begin{aligned}
& f(z/2), && \text{if}\, z[j] \text{ is even},\;\forall 1\le j\le d, \\
& \mr{In}_k(f)(\lfloor z/2\rfloor),\;&& \text{else, with $k= z - 2\lfloor z/2\rfloor$, }
\end{aligned}
\end{equation}
where $\lfloor z/2 \rfloor$ denotes the lower integer part of each coordinate of $z/2$. We note that in the second case, $k[j]=0$ for even $z[j]$ and $k[j]=1$ for odd $z[j]$, $j=1,...,d$. The operator $P$ can be thought of as the polynomial interpolation operator at $k/2$ given the function values $ \bm{f}(z+\mr{pa}(k/2))$ at the corresponding parent set $\mr{pa}(k/2)$.

\paragraph{Alignment with infinite grid}
This part aligns an arbitrarily function $f$ defined on the finite grid $\tilde{\X}_n$ with an appropriate function $f'$ defined on the infinite grid $\mathbb{Z}^d$, such that, under proper scaling, both functions yield the same polynomial interpolation over the grid.

Specifically, for all $0\le j_1\le j_2\le\eta(n)-1$ and $f\in L_(\tilde{\N}_{j_1})$, the objective of this paragraph is to find an integer $j_3\in\mathbb{N}$ satisfying $j_1\le j_3\le j_2$ and a function $f'\in L_1(\mathbb{Z}^d)$, such that 
\begin{equation}\label{eq:align 1}
P_{j_2} P_{j_2-1} \cdots P_{j_1}(f)(x) = P^{j_2-j_3}(f')(2^{j_2+1} x),  \;\; \forall x\in \tilde{\X}_n\cap\mathcal{G}_{j_2+1},
\end{equation}
and
\begin{equation}\label{eq:align 2}
\|f'\|_\infty \lesssim \|f\|_\infty.
\end{equation}

 First, let us consider the case  $j_0\leq j_1$ and define
 $$
f'(2^{j_1}x) = \left\{ \begin{aligned}
& f(x), \;\; &&\text{if\,\,} x\in \cup_{j=0}^{j_1}\tilde{\N}_j, \\
& 0, \;\; &&\text{if\,\,} x\not\in \cup_{j=0}^{j_1}\tilde{\N}_j.
\end{aligned} \right.
$$
Let $j_3=j_1$, then (\ref{eq:align 2}) is immediate. Furthermore, in view of Algorithm \ref{alg:dag arg}, for all $X_i\in\cup_{j=j_1}^{\eta(n)}\tilde{\N}_j\subset \cup_{j=j_0}^{\eta(n)}\tilde{\N}_j$, its ancestors are in $\tilde\X_n$. Therefore, $X_i$ has no boundary effect, implying directly  (\ref{eq:align 1}).

For $j_1<j_0$, let us define an intermediate function $f_{j_0}\in L_1(\cup_{j=0}^{j_0}\tilde{\N}_j)$ as
$$f_{j_0} = P_{j_0-1} P_{j_0-2} \cdots P_{j_1} f.$$
Then define the function $f'$ as
$$
f'(2^{j_0}x) = \left\{ \begin{aligned}
& f_{j_0}(x), \;\; && x\in \cup_{j=0}^{j_0}\tilde{\N}_j, \\
& 0, \;\; && x\not\in \cup_{j=0}^{j_0}\tilde{\N}_j.
\end{aligned} \right.
$$
Let $j_3=j_0$, then $\|f'\|_\infty = \|f_{j_0}\|_\infty$ holds. Furthermore, since $j_0$ only depends on $d,\alpha$ and $P_j,\forall j$ has bounded $ L_\infty$ norm, we also have $\|f_{j_0}\|_\infty \lesssim \|f\|_\infty$. Thus, the upper bound (\ref{eq:align 2}) is satisfied. Furthermore, since for arbitrary $X_i\in\cup_{j=j_0}^{\eta(n)}\tilde{\N}_j$, there is no boundary effect, equation (\ref{eq:align 1}) also holds.

\paragraph{Connecting operator norms}
Assertions (\ref{eq:align 1}) and (\ref{eq:align 2})
align a function defined on the finite grid $\tilde\X_n$ with an appropriate function defined on the infinite grid $\mathbb{Z}^d$. 
The polynomial interpolation operator $P$ possesses a localization property in the sense that the value of $P(f)(x)$ depends only on $f(x')$ for a limited number of points $x'$ that are near $x$. Building on this idea, by introducing a surrogate function on $\mathbb{Z}^d$, we can derive sharp bounds for the operator norms, which are suitable analysis with Fourier methods.
Details are provided below.

For $j_3$ and $f'\in L_1(\mathbb{Z}^d)$ satisfying (\ref{eq:align 1}) and (\ref{eq:align 2}), and $X_i\in\tilde{\X}_n$, define the function $f_i'\in L_1(\mathbb{Z}^d)$ as
$$
f_i'(2^{j_3}x) = \left\{ \begin{aligned}
& f'(2^{j_3}x), \;\; && \|x-X_i\|_\infty \le 2c_L\cdot 2^{-j_3}, \\
& 0, \;\; && \mr{otherwise}.
\end{aligned} \right.
$$
In view of assertion (\ref{eq:ancestor dist}) , the set $\{x:\|x-X_i\|_\infty \le 2c_L\cdot 2^{-j_3}\}$  contains all the ancestors of $X_i$ that belong to the layers $\cup_{j=j_3}^{j_2-1}\tilde{\N}_j$.

Therefore, for the above $X_i$, 
\begin{equation}\label{eq:align local 1}
P_{j_2} P_{j_2-1} \cdots P_{j_1}(f)(X_i) = P^{j_2-j_3}(f')(2^{j_2+1} X_i) = P^{j_2-j_3}(f_i')(2^{j_2+1} X_i),  
\end{equation}
and
\begin{equation}\label{eq:align local 2}
\|f_i'\|_1 \le(4c_L+1)^d \|f'\|_\infty \lesssim \|f\|_\infty.
\end{equation}
Furthermore, by recalling (\ref{eq:align local 1}), (\ref{eq:align local 2}) and noting that  in view of the finite co-domain $ L_1(\cup_{j=0}^{j_1}\tilde{\N}_j)= L_\infty(\cup_{j=0}^{j_1}\tilde{\N}_j)$, we have
\begin{align}\label{eq:poly interp infinite}
\|P_{j_2} P_{j_2-1} \cdots P_{j_1}\|
= & \sup_{f\in L_1(\cup_{j=0}^{j_1}\tilde{\N}_j)}\frac{\|P_{j_2} P_{j_2-1} \cdots P_{j_1} (f)\|_\infty}{\|f\|_\infty} \nonumber\\
= & \sup_{f\in L_1(\cup_{j=0}^{j_1}\tilde{\N}_j)}\frac{ \sup_{X_i\in\tilde{\X}_n} |P_{j_2} P_{j_2-1} \cdots P_{j_1} (f)(X_i)|}{\|f\|_\infty} \nonumber\\
\lesssim &\sup_{f_i'\in L_1(\mathbb{Z}^d)}\sup_{X_i\in\tilde{\X}_n} \frac{ |P^{j_2-j_3}(f_i')(2^{j_2+1} X_i)|}{\|f_i'\|_1} \nonumber\\
\le & \sup_{j\in\mathbb{N}}\sup_{f_i'\in L_1(\mathbb{Z}^d)} \frac{\|P^{j}(f_i')\|_\infty}{\|f'_i\|_1}.
\end{align}
Assertion (\ref{eq:poly interp infinite}) basically bounds the $ L_\infty$ norm of the recursive polynomial interpolation operator on a finite grid with the supremum norm of the recursive polynomial interpolation operator on the infinite grid. From now on we will focus on bounding the latter term.

\subsubsection{Bounding the interpolation operator via Fourier transform}
The polynomial interpolation operator $P$ on the infinite grid is also a convolution operator, which motivates us to use Fourier analysis. For all $f\in L_1(\mathbb{Z}^d)$, its (discrete time) Fourier transform is denoted by
$\mathcal{F}: L_1(\mathbb{Z}^d) \to  L_1(\mathbb{T}([0,1]^d))$, with $\mathbb{T}([0,1]^d)$ denoting the $d$-dimensional torus, such that for all $\xi\in \mathbb{T}([0,1]^d)$, 
\begin{equation}\label{eq:DTFT}
\mathcal{F}(f)(\xi) = \sum_{z\in\mathbb{Z}^d} e^{-2\pi\iota \langle z, \xi\rangle} f(z).
\end{equation}
 Furthermore, for $\rho\in L_1(\mathbb{T}([0,1]^d))$ and $z\in\mathbb{Z}^d$, we also define the inverse (discrete time) Fourier transform as
\begin{equation}\label{eq:DTFT inverse}
\mathcal{F}^{-1}(\rho)(z) = \int_{\xi\in\mathbb{T}([0,1]^d)} e^{2\pi\iota \langle z, \xi\rangle} \rho(\xi) d\xi.
\end{equation}
For all $f\in L_1(\mathbb{Z}^d)$, we have the following identity
$$\mathcal{F}^{-1}(\mathcal{F}(f)) = f.$$

We now proceed to the proof of Lemma \ref{lem:poly recursive}.

\begin{proof}[Proof of Lemma \ref{lem:poly recursive}]\mbox{}
In view of the analysis above, it suffices to bound the right-hand side term in (\ref{eq:poly interp infinite}). The proof is divided into two parts. First, we consider the special case where $d=1$ and $m=\underline{\alpha}+1$ is an even integer. This simpler case illustrates the main ideas of the proof. We show that the polynomial interpolation operator $P$, acting on the function $f$ in the Fourier domain, can be expressed as a Markov operator acting on $\mathcal{F}(f)$. Building on this observation, we derive the stated upper bound for the recursive application of $P$. In the second part, we address the more general cases. Here, the previous technique cannot be applied directly, and additional technical steps are required.\\

\textbf{Case 1:} We first consider the case $d=1$ and $m=\underline{\alpha}+1$ an even integer. Then the multi-index set $\mathbb{K}$, see \eqref{eq:K set}, has only a single element $k=1$ and the parent set of $k/2$ is $A=\{-m/2+1,-m/2+2,\cdots,m/2-1,m/2\}$. The operator $\mr{In}_1$, given in \eqref{def:Lk} for $k=1$, is
$$\mr{In}_1(f)(z) = v_{1/2}^T V_{A}^{-1} \bm{f}(z+A).$$
By symmetry, the $i$th and $(m-i)$th elements of the vector $V_A^{-1}v_{1/2}$ are the same. Thus, we can write $V_A^{-1}v_{1/2} = (b_{m/2},b_{m/2-1},\cdots, b_1,b_1,b_2,\cdots,b_{m/2})^T$ and
$$\mr{In}_1(f)(z) = \sum_{j=1}^{m/2} b_j[f(z+j) + f(z-j+1)].$$
Because $\|f\|_1<\infty$, the Fourier transform of $f$ exists and takes the form
$$
\mathcal{F}(\mr{In}_1(f))(\xi) = \sum_{j=1}^{m/2} b_j (e^{2j\pi\iota \xi} + e^{2(-j+1)\pi\iota \xi} ) \mathcal{F}(f)(\xi), \;\;\forall \xi\in\mathbb{T}([0,1]).
$$
Then, in view of the definition of the operator $P$ in (\ref{eq:fourier combine}), the Fourier transform of $P(f)$ can be computed as
\begin{align*}
\mathcal{F}(P(f))(\xi) = &\sum_{z=-\infty}^\infty P(f)(z) e^{-2\pi\iota z\xi} \nonumber\\
= & \sum_{z=-\infty}^\infty P(f)(2z) e^{-2\pi\iota \xi\cdot 2z} + \sum_{z=-\infty}^\infty P(f)(2z+1) e^{-2\pi\iota \xi\cdot (2z+1)} \nonumber\\
= & \sum_{z=-\infty}^\infty f(z) e^{-2\pi\iota \cdot 2\xi\cdot z} + \sum_{z=-\infty}^\infty L_1(f)(z) e^{-2\pi\iota \cdot 2\xi\cdot z} e^{-2\pi\iota \xi} \nonumber\\
= & \mathcal{F}(f)(2\xi) + e^{-2\pi\iota\xi}\mathcal{F}(\mr{In}_1(f))(2\xi) \nonumber\\
= & \mathcal{F}(f)(2\xi)\left[ 1 + \sum_{j=1}^{m/2} b_j (e^{-2\pi\iota(2j-1)\xi} + e^{-2\pi\iota(-2j+1)\xi} ) \right] \nonumber\\
= & \mathcal{F}(f)(2\xi)\bigg[ 1 + \sum_{j=1}^{m/2} b_j 2\cos(2\pi(2j-1)\xi) \bigg].
\end{align*}
Let us denote the above linear operator by $M:  L_1(\mathbb{T}([0,1]))\to  L_1(\mathbb{T}([0,1]))$, i.e.
\begin{equation}\label{eq:F-M}
\mathcal{F}(P(f)) = M (\mathcal{F}(f)).
\end{equation}

Next, we show that $M$ is a Markov operator on densities. That is, for all density $\rho(\xi)$ on $\mathbb{T}([0,1])$, $M\rho$ remains  a probability density and hence
\begin{equation}\label{eq:markov 2}
\|M \rho\|_1 = \|\rho\|_1, \forall \rho \in  L_1(\mathbb{T}).
\end{equation}
By definition, $\mathcal{F}(P(f))$ is the Fourier transform of the polynomial interpolation of the function $f$. Thus, $\sum_{j=1}^{m/2} b_j 2\cos(2\pi(2j-1)\xi)$ corresponds to the frequency response of polynomial interpolation in the literature of signal processing, whose complex modulus does not exceed one (see page 42 of 
\cite{laakso1996split}). Thus $M\rho\geq 0$ whenever $\rho\geq 0$.  It remains to show that $M\rho$ integrates to one. We directly compute the integral as
\begin{align*}
\int_0^1 M\rho(\xi) d\xi 
= & \int_0^{1} \rho(2\xi) \bigg[ 1 + \sum_{j=1}^{m/2} b_j 2\cos(2\pi(2j-1)\xi) \bigg] d\xi \\
= & \int_0^{1/2} \rho(2\xi) \bigg[ 1 + \sum_{j=1}^{m/2} b_j 2\cos(2\pi(2j-1)\xi) +1 + \sum_{j=1}^{m/2} b_j 2\cos(2\pi(2j-1)(\xi+1/2)) \bigg] d\xi \\
= & \int_0^{1/2} 2 \rho(2\xi) d\xi = \int_0^1 \rho(\xi) d\xi = 1.
\end{align*}

Then, for all $f$, define $\mathcal{F}(f)_+(\xi) = \max\{\mathcal{F}(f)(\xi),0\}$ and $\mathcal{F}(f)_-(\xi) = \mathcal{F}(f)(\xi) - \mathcal{F}(f)_+(\xi)$. Furthermore, define $f_+ = \mathcal{F}^{-1}(\mathcal{F}(f)_+)$ and $f_- = f - f_+$.
Combining equations (\ref{eq:F-M}) and (\ref{eq:markov 2}), we have
\begin{align*}
\|P^j(f)\|_\infty \le & \|\mathcal{F}(P^j(f))\|_1 = 
\|\mathcal{F}(P^j(f_+)) + \mathcal{F}(P^j(f_-))\|_1 \\
\le & \|\mathcal{F}(P^j(f_+))\|_1 + \|\mathcal{F}(P^j(f_-))\|_1 \\ 
\overset{(\ref{eq:F-M})}{=} &\|M^j(\mathcal{F}(f_+))\|_1 + \|M^j(\mathcal{F}(f_-))\|_1 \\
\overset{(\ref{eq:markov 2})}{=} & \|\mathcal{F}(f)_+\|_1 + \|\mathcal{F}(f)_-\|_1 \\
= & \|\mathcal{F}(f)\|_1 \le \|\mathcal{F}(f)\|_\infty \le \|f\|_1,
\end{align*}
where in the last line, $\|\mathcal{F}(f)\|_1 \le \|\mathcal{F}(f)\|_\infty$ follows from the fact that $\mathcal{F}(f)$ is defined on the torus $\mathbb{T}([0,1]^d)$. The inequalities $\|\mathcal{F}(f)\|_\infty \le \|f\|_1$ and $\|P^j(f)\|_\infty \le  \|\mathcal{F}(P^j(f))\|_1$ follow from the definition of Fourier transform \eqref{eq:DTFT} and \eqref{eq:DTFT inverse}, respectively. Therefore, the right-hand side of \eqref{eq:poly interp infinite} is bounded from above by one, which in turn implies that $\vartheta_n=O(1)$.\\

\textbf{Case 2:}
Next we extend the results to the more general case of $d$ and $\alpha$. The proof is divided into two steps. In the first step, analogous to Case 1, we study the operator $P$ in the Fourier domain, showing that it corresponds to a linear operator $M$ and identifying the associated kernel $\mathcal{K}_M$. However, in this general setting, $M$ is not necessarily a Markov operator and therefore cannot be directly used to derive the stated upper bounds. In the second step we address this issue by constructing a modified operator $\tilde{M}$, which is then employed in the analysis.\\

{\it \underline{Step 1}}
Recall the multi-index set $\mathbb{K}$ defined in (\ref{eq:K set}), and for $k\in\mathbb{K}$, introduce the vector 
$$(b_{k,1},b_{k,2},\cdots,b_{k,m})^T=V_{\mr{pa}(k/2)}^{-1}v_{k/2}.$$
Denote the parent set of $k/2$ as $\mr{pa}(k/2) = \{z_{k,1},z_{k,2},\cdots,z_{k,m}\}$, which is a corner set as in Section \ref{sec:build dag} with $m ={\underline{\alpha}+d \choose \underline{\alpha}}$.
Then the discrete time Fourier transform of $\mr{In}_k(f)$, given in \eqref{def:Lk}, is
$$\mathcal{F}(\mr{In}_k(f))(\xi) = \sum_{k=1}^m b_{k,m} e^{2\pi \iota \langle z_{k,m},\xi\rangle}\mathcal{F}(f)(\xi).$$
In view of the definition of the operator $P$ in (\ref{eq:fourier combine}),
\begin{align}\label{eq:freq p multi}
\mathcal{F}(P(f))(\xi) = & \sum_{z\in\mathbb{Z}^d} P(f)(z) e^{-2\pi\iota \langle z,\xi\rangle} \nonumber\\
= & \sum_{k: 0\le k[j]\le 1} \sum_{z\in\mathbb{Z}^d} P(f)(2z+k) e^{-2\pi\iota \langle 2z+k,\xi\rangle} \nonumber\\
= & \sum_{z\in\mathbb{Z}^d}f(z)e^{-2\pi\iota \langle z,2\xi\rangle} + \sum_{k: |k|\ge 1, 0\le k[j]\le 1} \sum_{z\in\mathbb{Z}^d} \mr{In}_k(f)(z) e^{-2\pi\iota \langle z+k/2, 2\xi\rangle} \nonumber\\
= & \mathcal{F}(f)(2\xi) + \sum_{k: |k|\ge 1, 0\le k[j]\le 1} e^{-2\pi\iota \langle k,\xi\rangle} \mathcal{F}(\mr{In}_k(f))(2\xi) \nonumber\\
= & \mathcal{F}(f)(2\xi) \left( 1 + \sum_{k: |k|\ge 1, 0\le k[j]\le 1} e^{-2\pi\iota \langle k,\xi\rangle} \sum_{k=1}^m b_{k,m} e^{4 \pi \iota \langle z_{k,m},\xi\rangle}  \right).
\end{align}

Next, let us define the function $\psi:\mathbb{T}([0,1]^d)\to\mathbb{C}$ as
$$\psi(\xi) = 1 + \sum_{k: |k|\ge 1, 0\le k[j]\le 1} e^{-2\pi\iota \langle k,\xi\rangle} \sum_{i=1}^m b_{k,m} e^{4 \pi \iota \langle z_{k,m},\xi\rangle}$$
and a linear operator $M: L_1(\mathbb{T}([0,1]^d)) \to  L_1(\mathbb{T}([0,1]^d)$ as
$$M(\rho)(\xi) = \rho(2\xi) \psi(\xi).$$
Then equation (\ref{eq:freq p multi}) can be reformulated as
\begin{equation}\label{eq:M high-d}
\mathcal{F}(P(f)) = M (\mathcal{F}(f)).
\end{equation}

Similarly to the first case, we continue with investigating $M$ in the Fourier domain. Let us define the kernel $\mathcal{K}_M: \mathbb{T}([0,1]^d)\times \mathscr{B}(\mathbb{T}([0,1]^d)) \to \mathbb{R}$, such that for all $A \in  \mathscr{B}(\mathbb{T}([0,1]^d))$ and $\xi_1\in \mathbb{T}([0,1]^d)$,
\begin{equation}\label{eq:kernel_M}
\mathcal{K}_M(\xi_1,A) = \frac{1}{2^d} \sum_{\xi_2\in\mathbb{T}([0,1]^d):\, 2\xi_2 = \xi_1} \psi(\xi_2) \mathbbm{1}_{\xi_2\in A} .
\end{equation}

Next, we show that the above kernel $\mathcal{K}_M$ is associated to the linear operator $M$. It is easy to see that for all $\xi_1\in\mathbb{T}([0,1]^d)$, $\mathcal{K}_M(\xi_1,\cdot)$ defines an $ L_1$ function on $\mathbb{T}([0,1]^d)$ that integrates to one.
Moreover, letting the multi-index set $\mathbb{K}_0$ be
$$\mathbb{K}_0 = \left\{ k\in\mathbb{N}^d: 0\le k[j]\le 1,\forall 1\le j \le d \right\},$$
for all $\rho \in  L_1(\mathbb{T}([0,1]^d))$,

\begin{align*}\label{eq:M 2 kernel}
 \int_{\mathbb{T}([0,1]^d)} \rho(\xi_1) \mathcal{K}_M(\xi_1,A) d\xi_1
= & \int_{\mathbb{T}([0,1]^d)}  \rho(\xi_1) 
  \frac{1}{2^d} \sum_{\xi_2\in\mathbb{T}([0,1]^d):\, 2\xi_2 = \xi_1} \psi(\xi_2) \mathbbm{1}_{\xi_2\in A} d\xi_1  \nonumber \\
= & \int_{\mathbb{T}([0,1]^d)} 
  \frac{1}{2^d} \sum_{\xi_2 = \frac{\xi_1+k}{2},\, k\in\mathbb{K}_0} \rho(2\xi_2) \psi(\xi_2) \mathbbm{1}_{\xi_2\in A}  \; d\xi_1\nonumber \\
= & \sum_{k\in\mathbb{K}_0} \int_{\frac{\mathbb{T}([0,1]^d) + k}{2}} \rho(2\xi_2) \psi(\xi_2) \mathbbm{1}_{\xi_2\in A}  \; d\xi_2\nonumber \\
= & \int_{A} \rho(2\xi_2) \psi(\xi_2) d\xi_2 = \int_A M(\rho)(\xi_2) d \xi_2,
\end{align*}
hence, $\mathcal{K}_M$ is indeed the kernel associated with the operator $M$. In particular, for $A = \mathbb{T}([0,1]^d)$,
\begin{equation}\label{eq:M 2 kernel}
\int_{\mathbb{T}([0,1]^d)} M(\rho)(\xi) d \xi = \int_{\mathbb{T}([0,1]^d)} \rho(\xi) \mathcal{K}_M(\xi,\mathbb{T}([0,1]^d)) d\xi.\\
\end{equation}

{\it \underline{Step 2}}
Note that, in view of assertion (\ref{eq:M 2 kernel}), if $\mathcal{K}_M(\xi,\cdot)$ were a probability measure, then the same reasoning as in Case 1 would apply. However, $\mathcal{K}_M(\xi,\cdot)$ may take complex values. 
To address this issue, we transform $M$ into a new operator whose associated kernel provides finite, non-negative measures.

First, for $\xi\in \mathbb{T}([0,1]^d)$, define the vector $ L_1$ norm on the torus as
$$\|\xi\|_{ L_1} = \sum_{j=1}^d \min\{\xi[j],1-\xi[j]\}.$$
Note that the $ L_1$ norm of $\xi$ is the $ L_1$ distance between $\xi$ and the origin on torus.
Furthermore, we define the operator $\Xi: L_1(\mathbb{T}([0,1]^d)) \to  L_1(\mathbb{T}([0,1]^d))$ such that for all $\xi \in \mathbb{T}([0,1]^d)$,
$$\Xi (\rho)(\xi) = \|\xi\|_{ L_1} |\rho(\xi)|.$$
Then, for any function $f\in L_1(\mathbb{Z}^d)$ and for all $z\in\mathbb{Z}^d$,
\begin{align*}
f(z) = &\int_{\mathbb{T}([0,1]^d)} e^{2\pi\iota \langle z,\xi\rangle} \mathcal{F}(f)(\xi) d\xi\\
= &\sum_{k\in\mathbb{N}^d:0\le k[j]\le 1} \int_{[0,1/2]^d+k/2} e^{2\pi\iota \langle z,\xi\rangle} \mathcal{F}(f)(\xi) d\xi \\
= &\sum_{k\in\mathbb{N}^d:0\le k[j]\le 1} \int_{[0,1/2]^d+k/2} e^{2\pi\iota \langle z,\xi-k\rangle} \mathcal{F}(f)(\xi) d\xi.
\end{align*}
By the definition of the inverse Fourier transform \eqref{eq:DTFT inverse}, $\mathcal{F}(f)$ is an $ L_\infty$ function on $\mathbb{T}([0,1]^d)$. Therefore, we can extend the domain of the function $f$ from $\mathbb{Z}^d$ to $\mathbb{R}^d$ as
$$f_{\mr{ext}}(x) = \sum_{k\in\mathbb{N}^d:0\le k[j]\le 1} \int_{[0,1/2]^d+k/2} e^{2\pi\iota \langle x,\xi-k\rangle} \mathcal{F}(f)(\xi) d\xi,$$
such that $f_{\mr{ext}}(z) = f(z), \forall z\in\mathbb{Z}^d$ and $f_{\mr{ext}}\in L_\infty(\mathbb{R}^d)$. Moreover, $f_{\mr{ext}}$ has uniformly bounded gradient, i.e. for all $x\in\mathbb{R}^d$
\begin{align*}
\|\nabla f_{\mr{ext}}(x)\|_2
= & \bigg\|\sum_{k\in\mathbb{N}^d:0\le k[j]\le 1} \nabla  \int_{[0,1/2]^d+k/2} e^{2\pi\iota \langle x,\xi-k\rangle} \mathcal{F}(f)(\xi) d\xi\bigg\|_2 \\
= & \bigg\|\sum_{k\in\mathbb{N}^d:0\le k[j]\le 1}   \int_{[0,1/2]^d+k/2} 2\pi(\xi-k)\iota e^{2\pi\iota \langle x,\xi-k\rangle} \mathcal{F}(f)(\xi) d\xi\bigg\|_2 \\
\le & \sum_{k\in\mathbb{N}^d:0\le k[j]\le 1}   \int_{[0,1/2]^d+k/2} \|2\pi(\xi-k)\iota e^{2\pi\iota \langle x,\xi-k\rangle} \mathcal{F}(f)(\xi)\|_2 d\xi \\
\le & 2\pi \sum_{k\in\mathbb{N}^d:0\le k[j]\le 1}   \int_{[0,1/2]^d+k/2} \|\xi\|_{ L_1}|\mathcal{F}(f)(\xi)| d\xi \\
= & 2\pi \int_{\mathbb{T}([0,1]^d)} \|\xi\|_{ L_1}|\mathcal{F}(f)(\xi)| d\xi \\
= & 2\pi \|\Xi(\mathcal{F}(f))\|_{1}.
\end{align*}
Furthermore, recall from the definition of \eqref{eq:fourier combine} that $P(f)(2z) = f(z), \;\forall z\in\mathbb{Z}^d$. Hence, for all $z\in\mathbb{Z}^d$, there exists a $z'\in\mathbb{Z}^d$ satisfying $\|z-z'\|_2\le\sqrt{d}$, such that $P(f)(z') = f(z'/2)$. Therefore, the supremum norm of $P(f)$ can be controlled as
\begin{equation}\label{eq:sup norm deriv}
\|P(f)\|_\infty \le \|f\|_\infty + 2\pi\sqrt{d} \|\Xi (\mathcal{F}(P(f)))\|_1 = \|f\|_\infty + 2\pi\sqrt{d} \|\Xi (M(\mathcal{F}(f)))\|_1.
\end{equation}

Note that assertion (\ref{eq:sup norm deriv}) provides a convenient way for controlling the supremum norm of $P^j(f)$ with the help of the $ L_1$ norm of $\Xi (M^j(\mathcal{F}(f)))$. However, the operator $M$ is not Markov and hence the arguments applied in Case 1 do not go through. Therefore, we define an auxiliary operator $\tilde{M}$, 
such that, for all $\rho\in  L_1(\mathbb{T}([0,1]^d))$ satisfying $\rho(k)=0,k\in\mathbb{K}_0$,
\begin{equation}\label{eq:op M t}
\tilde{M}(\rho)(\xi) = \frac{\|\xi\|_{ L_1}}{\|2\xi\|_{ L_1}}\rho(2\xi) |\psi(\xi)|,    
\end{equation}
where the function is interpreted as its limit when $\|2\xi\|_{ L_1}=0$. 
It is easy to see that
\begin{equation}\label{eq:tilde M}
\tilde{M}(\Xi(\rho))(\xi) = \Xi(M(\rho))(\xi).
\end{equation}
Note that in contrast to the operator $M$, which maps the Fourier transform of the function to the Fourier transform of its polynomial interpolation, the operator $\tilde{M}$ acts on $\Xi(\rho)$, which is associated with the gradient of the function.

Then, by applying (\ref{eq:sup norm deriv}) and (\ref{eq:tilde M}) recursively, we obtain for all $j\in\mathbb{N}$,
\begin{equation}\label{eq:sup norm recur 2}
\|P^j(f)\|_\infty \le \|f\|_\infty + 2\pi\sqrt{d} \sum_{j'=1}^j \|\tilde{M}^{j'} (\Xi(\mathcal{F}(f)))\|_1.
\end{equation}
Hence, instead of bounding the norm of $M^j\mathcal{F}(f)$ as in Case 1, we need to bound the sum of the norms of
$\tilde{M}^j\Xi(\mathcal{F}(f))$. Then, similarly to the definition of $\mathcal{K}_M$ in \eqref{eq:kernel_M}, for all $A=(a_1,a_2)^d\subset\mathbb{T}([0,1]^d)$, let
$$
\mathcal{K}_{\tilde{M}}(\xi_1,A) = \frac{1}{2^d} \sum_{\xi_2\in\mathbb{T}([0,1]^d): 2\xi_2 = \xi_1} \frac{\|\xi_2\|_{ L_1}}{\|\xi_1\|_{ L_1}} |\psi(\xi_2)| \mathbbm{1}_{\xi_2\in A}.
$$
Then, for any non-negative function $\rho\in L_1(\mathbb{T}([0,1]^d))$,
\begin{equation}\label{eq:M tilde 2 kernel}
\int_{\mathbb{T}([0,1]^d)} \tilde{M}(\rho)(\xi) d \xi = \int_{\mathbb{T}([0,1]^d)} \rho(\xi) \mathcal{K}_{\tilde{M}}(\xi,\mathbb{T}([0,1]^d)) d\xi.
\end{equation}
Note that for all $\xi\in \mathbb{T}([0,1]^d)$, $\mathcal{K}_{\tilde{M}}(\xi,\cdot)$ is a nonnegative finite measure.

Next, assume that the measure of the whole space is bounded away from one, i.e. $\exists c\in(0,1)$, such that
\begin{equation}\label{eq:measure bound}
\sup_{\xi\in\mathbb{T}([0,1]^d)} \mathcal{K}_{\tilde{M}}(\xi,\mathbb{T}([0,1]^d)) \le 1 - c.
\end{equation}
Then, in view of  $\Xi(\rho)\geq 0$, for all $f\in L_1({\mathbb{Z}^d})$,
\begin{align}\label{eq:M tilde l1}
\|\tilde{M}(\Xi(\mathcal{F}(f)))\|_1 
= & \int_{\mathbb{T}([0,1]^d)} \Xi(\mathcal{F}(f))(\xi) \mathcal{K}_{\tilde{M}}(\xi,\mathbb{T}([0,1]^d)) d\xi \nonumber\\
\le & \int_{\mathbb{T}([0,1]^d)} \Xi(\mathcal{F}(f))(\xi) (1-c) d\xi = (1-c) \|\Xi(\mathcal{F}(f))\|_1.
\end{align}
Plugging (\ref{eq:M tilde l1}) into the upper bound (\ref{eq:sup norm recur 2}), results in
\begin{align*}
\|P^j(f)\|_\infty \le & \|f\|_\infty + 2\pi\sqrt{d}\sum_{j'=1}^j(1-c)^{j'}\|\Xi(\mathcal{F})(f)\|_1 \\
\le & \|f\|_1+ 2\pi\sqrt{d}\frac{1-c}{c}\|\mathcal{F}(f)\|_\infty \le \big(1+2\pi\sqrt{d}\frac{1-c}{c}\big) \|f\|_1,
\end{align*}
proving that the right-hand side of \eqref{eq:poly interp infinite} is bounded by a constant and therefore verifying the statement of the lemma.

Therefore, it remains to verify the upper bound (\ref{eq:measure bound}), where the LHS is a function on $\mathbb{T}([0,1]^d)$. 
While the analytical derivation of this bound is highly challenging, the function is uniformly continuous in $\xi$, which allows for numerical evaluation to obtain the required upper bound for arbitrary combinations of the smoothness $\alpha$ and dimension $d$.  The results for $d=1,2$ and various values of $\alpha$ are displayed in Figure \ref{fig:poly interp}. One can observe that, for all regularity parameters considered, the constant $c$ in equation (\ref{eq:measure bound}) is bound away from $1$. 
Our approach could potentially be extended to arbitrary combinations of $\alpha$ and $d$. However, in such general settings, the operator $\tilde{M}$ defined in \eqref{eq:op M t} may no longer guarantee that the maximum transition measure does not exceed $1$. One could consider extending the torus $ L_1$ norm to other functions depending $\alpha, d$ and redefine the operator accordingly. Nonetheless, this would introduce additional technical complexity and is beyond the scope of the present paper.
\begin{figure}
    \centering
    \includegraphics[width=0.8\linewidth]{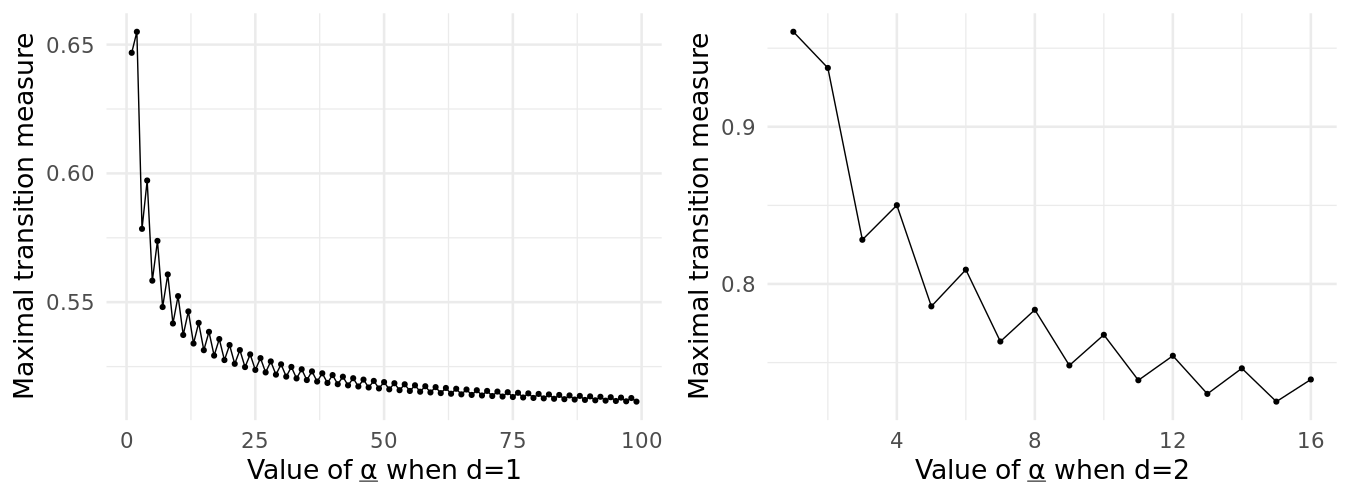}
    \caption{Evaluation of LHS of equation (\ref{eq:measure bound}) for $d=1,2$ and various regularity parameters.}
    \label{fig:poly interp}
\end{figure}
\end{proof}

\subsection{Proof of Lemma \ref{lem:GP recursive}}
\begin{proof}
We distinguish three cases according to the values of $j_1$ and $j_2$.

\textbf{Case 1} First, consider the case $\tau \gamma^{-j_2}\ge 1$. The key observation is that when the distance between $X_i$ and $\mr{pa}(X_i)$ is sufficiently large, the corresponding Gaussian interpolation term becomes negligibly small.

Since  $\tau \gamma^{-j_2}\ge 1$, we have $\tau \ge 1$ and  $\forall \; j_1\le j\le j_2$, $\tau \gamma^{-j}\ge 1$. For all $X_i\in \N_{j+1}$, define the set $A = \gamma^j \mr{pa}(X_i)$ and $x^* = \gamma^j X_i$. Then in view of Condition \ref{cond:layered}, the minimal distance among elements in the set $A$ is lower bounded by a constant $c_d$ and the minimal distance between set $A$ and $x^*$ is lower bounded by $c_d/\gamma$. Then in view of \eqref{def:op:G}, 
$$ G_j(f)(X_i) \le \|K_{\tau \gamma^{-j} A, \tau \gamma^{-j} A}^{-1}\|_2 \|K_{\tau \gamma^{-j} A, \tau \gamma^{-j} x*}\|_2 \|\boldsymbol{f}(\mr{pa}(X_i))\|_2.$$
Since $K$ is the Mat\'{e}rn covariance kernel, we can write $K_{x_1, x_2} = K_0(\|x_1-x_2\|_2)$, where $K_0(x)$ decays to zero as $x\to \infty$. Thus, for $\tau \gamma^{-j}$ sufficiently large, $K_{\tau \gamma^{-j} A, \tau \gamma^{-j} A}$ is  diagonally dominant, and 
we can bound its minimal eigenvalue using Gershgorin circle theorem (Theorem 6.1.1 of \citep{horn2012matrix}) as
$$ \lambda_{\min} (K_{\tau \gamma^{-j} A, \tau \gamma^{-j} A}) \ge 1-(m-1) \sup_{\substack{a_1,a_2\in A \\ a_1\ne a_2} } K_0(\|a_1-a_2\|_2)\ge 1- (m-1) K_0(\tau \gamma^{-j} c_d).$$
Thus
$$\|K_{\tau \gamma^{-j} A, \tau \gamma^{-j} A}^{-1}\|_2 \le   \big[ 1- (m-1) K_0(\tau \gamma^{-j} c_d) \big]^{-1}.$$

Upper bounding each element of the vector, we obtain the following $ L_2$-norm bound for the covariance between $X_i$ and $\mr{pa}(X_i)$:
$$ \|K_{\tau \gamma^j A, \tau \gamma^j x^*}\|_2 
\le \sqrt{ m \sup_{\substack{a_1,a_2\in A \\ a_1\ne a_2} } K_0(\|a_1-a_2\|_2)   }
\le \sqrt{m} K_0(\tau \gamma^{-j+1} c_d).$$
Therefore, for $X\in\N_{j+1}$, we have
$$ G_j(f)(X) \le [1- (m-1) K_0(\tau \gamma^{-j} c_d)]^{-1} \sqrt{m} K_0(\tau \gamma^{-j+1} c_d)\|\boldsymbol{f}(\mr{pa}(X))\|_2.$$
As $\tau \gamma^{-j} \to \infty$, $K_0(\tau \gamma^{-j} c_d)$ and $K_0(\tau \gamma^{-j+1} c_d)$ tend to zero and therefore, the above display also goes to zero. Thus, there exists a constant $c'$ independent of $j$, such that if $\tau \gamma^{-j}\ge c'$, $[1- (m-1) K_0(\tau \gamma^{-j} c_d)]^{-1} \sqrt{m} K_0(\tau \gamma^{-j+1} c_d)\le 1/\sqrt{m}$ and $G_j(f)(X) \le \sup_{X'\in\mr{pa}(X)}|f(X')|$. Since this holds for arbitrary $X_i\in\N_{j+1}$, in view of the definition of $G_j$ in \eqref{def:op:G},  $\|G_j\|\le 1$. 
If $\tau \gamma^{-j_2} \ge c'$, this proves the lemma. Otherwise, if $1\le \tau \gamma^{-j_2} < c'$, let us consider the index set 
$$\mathcal{J} = \{j: j_1\le j \le j_2, \tau \gamma^{-j} < c'\}.$$
Then the cardinality of $\mathcal{J}$ is bounded by an absolute constant
$\log_\gamma c'$. For all $j\in\mathcal{J}$, we use the norm bound 
$$G_j(f)(X) \le \|K_{\tau \gamma^j A, \tau \gamma^j A}\|^{-1}_2 \|K_{\tau \gamma^j A, \tau \gamma^j x^*}\|_2 \|f(\mr{pa}(X)\|_2 \le m \|f(\mr{pa}(X)\|_2,$$
resulting in $\|G_j\| \le m$. For all $j_1\le j\le j_2$ and $j\not\in\mathcal{J}$, we still have $\|G_j\| \le 1$, thus
$$\|G_{j_2} \cdots G_{j_1+1}  G_{j_1} \| \le m^{\log_\gamma c'}. $$
\vspace{0.1cm}

\textbf{Case 2} Next, we assume that $\tau \gamma^{-j_1}\le 1$, in which case $\tau \gamma^{-j}<1, \;\forall j_1\le j \le j_2$. The key idea is that as elements in $\mr{pa}(X_i)\cup\{X_i\}$ become increasingly close to each other, the flat-limit phenomenon of Lemma \ref{lem:GP flat} takes effect, and the Gaussian interpolation can be well approximated by polynomial interpolations.

For all $j\in\mathbb{N}$, let us denote by $E_j = G_j - P_j$ the difference between Gaussian and polynomial interpolation operators at layer $j$, with its supreme operator normed defined in the same way as $P_j$ in Section \ref{sec:rec:inter}. By Condition \ref{cond:norming}, for every $X_i\in\N_j$, the point $X_i$ together with its parent set is included in a cube $\C$, whose side length is at most $c_L\gamma^{-j}$. This configuration can be viewed as a rescaled version of the unit cube $[0,1]^d$ with scaling parameter $c_L\gamma^{-j}$. Combining this observation with Lemma \ref{lem:GP flat}, we obtain
\begin{equation*}
\|E_j\|\lesssim \gamma^{-2j(\alpha-\underline{\alpha})} + \gamma^{-j}.
\end{equation*}

The key idea is to show that the operator $E_j$ decays exponentially fast to zero as $j\rightarrow\infty$, which in particular implies that $\sum_{j=1}^{\infty}\|E_j\|<\infty$. With this summability at hand, we express the recursive Gaussian interpolation operator in terms of the polynomial interpolation operator together with the accumulated error contribution $E_j$. The finiteness of the series $\sum_{j=1}^{\infty}\|E_j\|$ will then allow us to control all error terms in the upper bound.

Let $\alpha^* = \min\{2(\alpha-\underline{\alpha}), 1\}$ and note that $\|E_j\|\lesssim\gamma^{-\alpha^* j}$. For $j=j_1$, write 
\begin{equation}
\label{eq:G init}
G_{j_1} = P_{j_1} + E_{j_1} = a_{j_1} P_{j_1} H_{j_1} + a_{j_1+1} H_{j_1+1}
\end{equation}
with
$a_{j_1}=1$, $a_{j_1+1}=\gamma^{-j_1\alpha^*}$, $H_{j_1}=I$(identity operator) and 
$$\|H_{j_1+1}\| = \|E_{j_1}/\gamma^{-j_1\alpha^*}\|\lesssim 1.$$
We prove below that for all $ j_1\le j\le j_2+1$, there exists an operator $\|H_{j+1}\|\lesssim 1$ and a constant $a_{j}\in [0,+\infty)$, such that 
\begin{equation}\label{induction:lem9}
G_j\cdots G_{j_1+1} G_{j_1} = \sum_{j'=j_1}^j a_{j'} P_{j} \cdots P_{j'+1}P_{j'} H_{j'} + a_{j+1} H_{j+1}.
\end{equation}
Moreover, the constants $a_{j}, j_1\le j\le j_2+1$ satisfy
$$a_{j+1} \lesssim \gamma^{-j\alpha^*} \vartheta_n \sum_{j'=j_1}^{j}a_{j'}.$$
Then, let us introduce the notation $S_j =\sum_{j'=j_1}^j a_{j'}$, and note that by the preceding display,
$$S_{j+1}-S_j \lesssim \gamma^{-j\alpha^*} \vartheta_n S_j.$$
After reorganizing the inequality and using logarithmic transformation, we arrive at
\begin{align*}
\ln S_j \leq & \ln \big(S_{j-1} (1+ c\vartheta_n \gamma^{-(j-1)\alpha^*})\big) \\
\leq & \ln S_{j_1} + \sum_{j'=0}^{j-1} \ln(1+c\vartheta_n\gamma^{-j'\alpha^*}) \le 
\ln a_{j_1} + c\vartheta_n\sum_{j'=0}^{\infty} \gamma^{-j'\alpha^*} \le \frac{c'\vartheta_n}{1-\gamma^{-\alpha^*}}, 
\end{align*}
for some universal constant $c'>0$. Therefore, in view of the definition of $\tilde{\vartheta}_n$,
$$\|G_{j_2} \cdots G_{j_1+1} G_{j_1}\| \le \vartheta_n S_{j_2} \le \vartheta_n \exp\left[\frac{c'\vartheta_n}{1-\gamma^{-\alpha^*}}\right] = \tilde{\vartheta}_n,$$
with $c=c'/(1-\gamma^{-\alpha*})$. This proves the claim of the lemma.

Hence it remains to prove \eqref{induction:lem9}. We proceed by induction.  By equation (\ref{eq:G init}), the claim holds for $j=j_1$. Now suppose the claim holds for $j=l$. Then, for $j=l+1$, we have
\begin{align*}
G_{l+1} G_l \cdots G_{j_1+1} G_{j_1}
= & (P_{l+1} + E_{l+1})\Big[ \sum_{j'=j_1}^l a_{j'} P_{l} \cdots P_{j'+1}P_{j'} H_{j'} + a_{l+1} H_{l+1} \Big] \\
= & \sum_{j'=j_1}^{l+1} a_{j'} P_{l+1} P_{l} \cdots P_{j'+1}P_{j'} H_{j'}\\
&\qquad+ E_{l+1} \Big[ \sum_{j'=j_1}^l a_{j'} P_{l} \cdots P_{j'+1}P_{j'} H_{j'} + a_{l+1} H_{l+1} \Big].
\end{align*}
Let $a_{l+2}=\|E_{l+1}\|\Big\| \sum_{j'=j_1}^l a_{j'} P_{l} \cdots P_{j'+1}P_{j'} H_{j'} + a_{l+1} H_{l+1} \Big\|$ and 
$$H_{l+2} =E_{l+1} \Big[ \sum_{j'=j_1}^l a_{j'} P_{l} \cdots P_{j'+1 }P_{j'} H_{j'} + a_{l+1} H_{l+1} \Big]/ a_{l+2},$$
then obviously $\|H_{l+2}\|\le 1$. 
Recalling the definition of $\vartheta_n$ and the induction assumption that $\|H_j\|_1\le 1, \forall j_1\le j\le l+1$, we have
$$a_{l+2} \lesssim \|E_{l+1}\|\vartheta_n\sum_{j=j_1}^{l+1} a_{j} \lesssim \gamma^{-(l+1)\alpha^*} \vartheta_n\sum_{j=j_1}^{l+1}a_{j},$$
concluding the proof of the lemma for Case 2.\\

\textbf{Case 3}  Finally, if $\tau \gamma^{-j_2} < 1 < \tau \gamma^{-j_1}$, then there exists $j_1< j_3 \le j_2$, such that $\tau \gamma^{-j_2} \le \tau \gamma^{-j_3} \le 1 < \tau \gamma^{-j_3+1} \le \gamma^{-j_1}$. For $j_1\le j\le j_3-1$, we can apply Case 2, whereas for $j_3\le j\le j_2$, we can apply Case 1. Combining both parts concludes the proof.
\end{proof}

\section{Proofs of the lemmas and theorems in Section \ref{sec:small ball} and \ref{sec:decentering}}
\subsection{Proof of Theorem \ref{lem:small ball}}
\begin{proof}
First, we prove the theorem for $s=1$ and abbreviate $\hat{Z}^{\tau,1}$ as $\hat{Z}^\tau$. Extending the result to a general $s>0$ is straightforward and can be accomplished by dividing $\epsilon$ by $s$.

Recall that $\eta(i)=j$ indicates that $X_i\in\mathcal{N}_j$. For every $X_i\in\mathcal{X}_n$, define $u_i$ as
$$u_i = \hat Z^\tau_{X_i}-\mathbb{E}(\hat Z^\tau_{X_i}|\hat Z^\tau_{\mr{pa}(X_i)}) = \hat{Z}^\tau_{X_i} - K_{\tau X_i,\tau\mr{pa}(X_i)} K_{\tau\mr{pa}(X_i),\tau\mr{pa}(X_i)}^{-1} \hat{Z}^\tau_{\mr{pa}(X_i)}.$$
Since $u_i$ is a linear transformation of a centered Gaussian random vector, it also has zero mean. Recall that in view of Theorem \ref{thm:var},  $\var(u_i) = \E[Z^\tau_{X_i}-\E(Z^\tau_{X_i}|Z^\tau_{\mr{pa}(X_i)})]^2 \asymp \tau^{2\alpha}\gamma^{-2\alpha\eta(i)}$.
Furthermore, let $j_0$ be the smallest integer satisfying $\gamma^{-\alpha j_0} \le \epsilon$ and define the sequence
\begin{equation*}
b_j = \Bigg\{
\begin{aligned}
& \tau^\alpha\gamma^{-(\frac{j_0-j}{2}+j_0)\alpha}, & \; j<j_0, \\
& \tau^\alpha\gamma^{-(\frac{j-j_0}{2}+j_0)\alpha}, & \; j\ge j_0.
\end{aligned}
\end{equation*}
Finally, for all $1\le i \le n$, define the event $A_i= \{|u_i|\le b_{\eta(i)}\}.$

The remainder of the proof consists of three parts. In the first part, we provide an upper bound for the $L_\infty$ norm of $\hat{Z}_{X_i}^\tau$ under the event $\cap_{i=1}^n A_i$. In the second part, we derive a lower bound for the probability of this event. In the final part, we combine the previous results and extend them to general scaling parameter $s>0$.\\

\textbf{Part 1:}\mbox{}
We restrict ourselves to the event $\cap_{i=1}^n A_i$. Then, for all $1\le j\le \eta(i)$, let us define the function $W_j: \cup_{j'=0}^j\N_{j'}\to\mathbb{R}$ as
$$W_j(X_i) = \Bigg\{ \begin{aligned}
    & u_i, \;\; & X_i \in \N_j,\\
    & 0, \;\; & X_i\not\in \N_j.
\end{aligned}$$
With a slight abuse of notation, for all $1\le j\le\eta(n)$, let $\hat{Z}^\tau_{\cup_{j'=0}^j\N_{j'}}$ denote the restriction of $\hat{Z}^\tau$ to the co-domain  $\cup_{j'=0}^j\N_{j'}$. Then, in view of \eqref{def:op:G}, we have for all $0\le j\le \eta(i)-1$,
$$\hat{Z}^\tau_{\cup_{j'=0}^{j+1}\N_{j'}} = W_{j+1}(\cup_{j'=0}^{j+1}\N_{j'}) + G_j \hat{Z}^\tau_{\cup_{j'=0}^{j+1}\N_{j'}}.$$
Therefore, by induction of the previous formula, and applying the triangle inequality, the definition of the event $A_i$, and Lemma \ref{lem:GP recursive}, we get 
\begin{align}\label{eq:on An bound}
\|\hat{Z}^\tau_{\X_n}\|_\infty = & \|\sum_{j=0}^{\eta(n)} G_{\eta(n)} G_{\eta(n)-1} \cdots G_{j+1} G_{j} W_{j} \|_\infty \nonumber\\
\le & \sum_{j=0}^{\eta(n)} \|G_{\eta(n)} G_{\eta(n)-1} \cdots G_{j+1} G_{j} W_{j}\|_\infty \nonumber\\
\lesssim & \tilde\vartheta_n  \sum_{j=0}^{\eta(n)} b_j
\le 2  \tilde\vartheta_n \tau^\alpha \frac{\gamma^{-\alpha j_0}}{1-\gamma^{-\alpha/2}} \lesssim \tilde\vartheta_n \tau^{\alpha} \epsilon.
\end{align}

\textbf{Part 2:}\mbox{}
We first prove by induction that $u_1, u_2, \cdots, u_n$ are independent random variables. A single random variable $u_1$ is trivially independent. Next, assuming that $u_1, u_2, \cdots, u_l$ are independent, we show that $u_{l+1}$ is independent of $(u_1, u_2, \cdots, u_l)$. By the definition of $u_i$, for $1\le i\le l+1$, we have
\begin{align*}
\mathbb{E}&\left\{ u_{l+1} [u_1,u_2,\cdots,u_l]^T \right\}\\
&\quad=  \mathbb{E} \left\{ \mathbb{E}\Big\{[\hat Z^\tau_{X_{l+1}}-\mathbb{E}(\hat Z^\tau_{X_{l+1}}|\hat Z^\tau_{\mr{pa}(X_{l+1})})] [u_1,u_2,\cdots,u_l]^T   \Big|Z^\tau_{X_i},1\le i\le l \Big\} \right\}=0,
\end{align*}
which implies independence, since the  $u_i$s are centered Gaussian random variables. Therefore,
\begin{align*}
\pr(\cap_{i=1}^n A_i) 
=  \prod_{i=1}^n \pr\Big(|U|\le b_{\eta(i)}/\sqrt{\var(u_i)} \Big) 
=  \prod_{j=1}^{\eta(n)} \left[\pr(|U|\le b_j \tau^{-\alpha}\gamma^{\alpha j})\right]^{|\mathcal{N}_j|},
\end{align*}
where $U\sim N(0,1)$. 

Next, recall some standard tail bounds for the standard normal distribution, i.e.
\begin{equation*}
\begin{aligned}
\pr(|U|\le u) &\ge u/3, \quad \forall \;0\le u\le 1,\\
\pr(|U|\le u) &\ge 1-e^{-u^2} \ge \exp\big(-2e^{-u^2/2}\big), \quad \forall \;u\ge 1,
\end{aligned}
\end{equation*}
see (7.2) and (7.3) in \cite{ledoux2006isoperimetry}, respectively. Then, noting that $b_j \tau^{-\alpha}\gamma^{\alpha j}\le 1$ if and only if $j\le j_0$ and that (\ref{eq:layered}) implies $ |\N_j|\lesssim \gamma^{dj}$, we obtain, for some constant $c_d>0$
\begin{align}\label{eq:prod}
\pr(\cap_{i=1}^n A_i)
\ge & \prod_{j=1}^{j_0} \big( \gamma^{-\frac{3}{2}\alpha(j_0-j)}/3\big)^{c_d\gamma^{dj}} 
\prod_{j=j_0+1}^{\eta(n)} \exp\big( -2c_d \gamma^{dj}  e^{- \gamma^{\alpha(j-j_0)}/2 } \big).
\end{align}
The logarithm of the first term in (\ref{eq:prod}) can be bounded as
\begin{align*}
-\ln\prod_{j=1}^{j_0} \big(\gamma^{-\frac{3}{2}\alpha(j_0-j)}/3\big)^{c_d\gamma^{dj}}
= & \sum_{j=1}^{j_0} c_d\gamma^{dj} \big[ \frac{3}{2}\alpha(j_0-j)\ln \gamma + \ln 3] \nonumber \\
= & c_d\gamma^{d j_0} \sum_{j=1}^{j_0} \gamma^{-d(j_0-j)} \big[ \frac{3}{2}\alpha(j_0-j)\ln 2\gamma + \ln 3]
\lesssim  \gamma^{d j_0}.
\end{align*}
The logarithm of the second term in (\ref{eq:prod}) can be bounded as
\begin{align*}\label{eq:log 2}
-\ln\prod_{j=j_0+1}^{\eta(n)}  \exp\big( -2 c_d\gamma^{dj}  e^{- \frac{1}{2}\gamma^{\alpha(j-j_0)} } \big)
= & \sum_{j=j_0+1}^{\eta(n)}  
2c_d \gamma^{dj}  e^{- \frac{1}{2}\gamma^{\alpha(j-j_0)} }
\lesssim  \gamma^{d j_0}.
\end{align*}
By combining the previous three displays and recalling that $\gamma^{-\alpha j_0}\asymp\epsilon$, we obtain
\begin{equation}\label{eq:event prob}
-\ln \pr(\cap_{i=1}^n A_i) \lesssim \epsilon^{-d/\alpha}.\\
\end{equation}

\textbf{Conclusion:}\mbox{}
In view of assertions (\ref{eq:on An bound}) and (\ref{eq:event prob}),
\begin{equation*}
-\ln \pr\left(\|\hat Z^\tau_{\mathcal{X}_n}\|_\infty\lesssim \tilde\vartheta_n \tau^\alpha \epsilon\right) \le -\ln \pr(\cap_{i=1}^n A_i) \lesssim \epsilon^{-d/\alpha}.
\end{equation*}
This implies, by replacing $\eps$ with $\eps\tau^{-\alpha}$, that
$$
-\ln \pr(\|\hat Z^\tau_{\mathcal{X}_n}\|_\infty\le \tilde\vartheta_n \epsilon) \lesssim \tau^d\epsilon^{-d/\alpha}.
$$
Finally, noticing that $\hat{Z}^{\tau,s} = s\hat{Z}^{\tau}$, provides the statement of the theorem, i.e.
$$-\ln \pr(\|\hat{Z}^{\tau,s}_{\mathcal{X}_n}\|_\infty <\tilde\vartheta_n \epsilon) \lesssim \tau^d s^{d/\alpha}\epsilon^{-d/\alpha},$$
where the constant in $\lesssim$ depends only on $d, \alpha$.
\end{proof}

\subsection{Proof of Lemma \ref{lem:RKHS min}}
\begin{proof}
Let $f\in\mathbb{H}$ be an arbitrary function satisfying $\bm{f}(A) = \bm{f}_0(A)$. Then, $f(x)-\sum_{i=1}^{|A|} a_i K(x_i,x) = 0$, $\forall x\in A$. In turn, the reproducing property implies that
$$\langle K(x,\cdot), f(\cdot)-\sum_{i=1}^{|A|} a_i K(x_i,\cdot) \rangle_{\mathbb{H}} = f(x)-\sum_{i=1}^{|A|} a_i K(x_i,x) = 0, \;\;\forall x\in A.$$
Therefore, $K(x,\cdot)$ is orthogonal to $f(\cdot)-\sum_{i=1}^{|A|} a_i K(x_i,\cdot)$ in $\mathbb{H}$, $\forall x\in A$. Thus,
\begin{align*}
\|f\|_{\mathbb{H}}^2 &= \|f(\cdot)-\sum_{i=1}^{|A|} a_i K(x_i,\cdot)\|_{\mathbb{H}}^2 + \|\sum_{i=1}^{|A|} a_i K(x_i,\cdot)\|_{\mathbb{H}}^2\\
& \ge \|\sum_{i=1}^{|A|} a_i K(x_i,\cdot)\|_{\mathbb{H}}^2 = \bm{f}_0(A)^T K_{A,A}^{-1} \bm{f}_0(A).
\end{align*}
The minimum 
is achieved if and only if $f(\cdot)-\sum_{i=1}^{|A|} a_i K(x_i,\cdot)=0$.
\end{proof}

\subsection{Proof of Lemma \ref{lem:Holder appro}}
\begin{proof}
In view of equation (\ref{eq:vecchia exp}),  without loss of generality, we can assume $s=1$ and simplify the notation by writing $Z^{\tau,1}$ as $Z^\tau$ and $\mathbb{H}^{\tau,1}$ as $\mathbb{H}^\tau$. 
The proof consists of three steps. In the first step, we control the interpolation error for constant functions. In the second step, we prove the lemma for $\beta\in(0,1]$. In the final step, we extend the proof to $\beta>1$.

Let us introduce the notation $B_i=\{X_i\}\cup\mr{pa}(X_i)$. By Condition \ref{cond:norming}, there exists a $d$-dimensional cube $\mathcal{C}\subset\X$ of side length $c_L\gamma^{-\eta(i)}$ that contains  $B_i$. Also, by Condition \ref{cond:layered}, the minimal $ L_\infty$ distance among elements in $\cup_{j'=0}^j \N_{j'}$ is $c_d \gamma^{-j}$. Since the cube $\C$ contains $X_i$ and its parent set, we have $c_d\le c_L$.

\textbf{Step 1:} Let $f_0\equiv c$ for some $c\in[-1,1]$. We prove that
\begin{equation}\label{eq:const interp}
\left|c-\mathbb{E}[Z^{\tau}_{X_{i}}|Z^{\tau}_{\mr{pa}(X_i)}=c\mathbbm{1}] \right| \lesssim \tau^{\alpha}\gamma^{-\eta(i)\alpha},
\end{equation}
where $\mathbbm{1}$ denotes the column vector of $1$s. Note, that in view of Lemma \ref{lem:RKHS min},
\begin{equation}\label{eq:H tau trans}
    \inf_{f\in\mathbb{H}^\tau,\bm{f}(B_i)=\bm{f}_0(B_i)} \|f\|_{\mathbb{H}^\tau}^2
= \bm{f}_0(B_i)^T K_{\tau B_i,\tau B_i}^{-1} \bm{f}_0(B_i) =
\inf_{f\in\mathbb{H},\bm{f}(\tau B_i)=\bm{f}_0(B_i)} \|f\|_{\mathbb{H}}^2.
\end{equation}
Next, let us consider a  function $\varphi_0(\cdot):\mathbb{R}\to\mathbb{R}$ such that $\varphi_0\in C^{\alpha+d/2}(\mathbb{R})_1$, $\varphi_0$ is supported on $[-2,2]$ and equals to a constant on $[-1,1]$. We define the function $\varphi_d(\cdot):\mathbb{R}^d\to\mathbb{R}$ as $\varphi_d(x) = \prod_{l=1}^d \varphi_0(x[l])$. Note that $\varphi_d\in C^{\alpha+d/2}(\mathbb{R}^d)_1$ is supported on $[-2,2]^d$, equals to a constant on $[-1,1]^d$, and has finite Sobolev norm $\|f\|_{W^{\alpha+d/2}}\lesssim 1$. We then split the proof into two cases based on whether $\tau \gamma^{-\eta(i)}\geq 1$.

If $\tau \gamma^{-\eta(i)}\le 1$, take the function
$$\tilde{f}(x) = \frac{c}{\varphi_d(0)}\varphi_d((x-\tau X_i)/c_L).$$
Note that $\tilde{f}(x)=c$ for all  $x\in \tau\mathcal{C}$.
In view of the proof of Lemma 11.36 of \cite{ghosal2017fundamentals}, the RKHS norm of the Mat\'{e}rn process is equivalent to the Sobolev norm. Therefore, by recalling \eqref{eq:H tau trans},
\begin{align}\label{eq:const RKHS 1}
\inf_{f\in\mathbb{H}^\tau,\, f(x)=c,\, x\in B_i} \|f\|_{\mathbb{H}^\tau}^2
\le & \|\tilde{f}\|_{\mathbb{H}}^2 \asymp  \|\tilde{f}\|_{W^{\alpha+d/2}}^2 \asymp \|\varphi_d\|_{W^{\alpha+d/2}}^2
\lesssim 1.
\end{align}

If $\tau \gamma^{-\eta(i)}\ge 1$, then $\tau\ge 1$. Take the function
$$\tilde{f}(x) = \sum_{w\in \tau B_i}\tilde{f}_w(x),\qquad \text{with $\tilde{f}_w(x)=\frac{c}{\varphi_d(0)}\varphi_d\Big(\frac{4\sqrt{d}}{c_d}(x-w)\Big)$}.$$
Note, 
$\tilde{f}_w(x)$ is supported on a cube centered at $w$ with side length $c_d/(2\sqrt{d})$, which belongs to a ball centered at $w$ with radius $ c_d/2\le c_d\tau \gamma^{-\eta(i)}/2 \le c_L\tau \gamma^{-\eta(i)}/2$. Thus, the supports of the functions $\tilde{f}_w(x)$, $w\in\tau B_i$, are 
not overlapping. These imply $\tilde{f}(x) = c, \forall x\in  \tau B_i$ and 
\begin{align}\label{eq:const RKHS 2}
\inf_{f\in\mathbb{H}^\tau,\, f(x)=c,\, x\in B_i} \|f\|_{\mathbb{H}^\tau}^2
\le & \|\tilde{f}\|_{\mathbb{H}}^2 \asymp \|\tilde{f}\|_{W^{\alpha+d/2}}^2\lesssim(m+1)\|\varphi_d\|_{W^{\alpha+d/2}}^2\lesssim 1.
\end{align}
Plugging in the bounds (\ref{eq:const RKHS 1}) and (\ref{eq:const RKHS 2}) into Lemma \ref{lem:GP interp} and using Theorem \ref{thm:var}, we get 
\begin{align}
\left|c-\mathbb{E}[Z^{\tau}_{X_{i}}|Z^{\tau}_{\mr{pa}(X_i)}=c\mathbbm{1}] \right| 
& \lesssim \sqrt{\var \big\{Z^{\tau}_{X_{i}}-\mathbb{E}[{Z}^{\tau}_{X_{i}}|{Z}^{\tau}_{\mr{pa}(X_{i})}]\big\} \inf_{f\in\mathbb{H}^\tau,\, f(x)=c,\, x\in B_i} \|f\|_{\mathbb{H}^\tau}^2 } \nonumber\\
& \lesssim \sqrt{\tau^{2\alpha}\gamma^{-2\eta(i)\alpha}\cdot 1 } = \tau^\alpha \gamma^{-\eta(i)\alpha}.\label{eq:lem13:case1}
\end{align}
\vspace{0.1cm}

\textbf{Step 2:} For $f_0\in C^{\beta}(\mathbb{R})_1$, $\beta\in(0,1]$, and for all $w_1, w_2\in [0,1]^d$,
$$\frac{|f_0(w_1)-f_0(w_2)|}{\|w_1-w_2\|_2^\beta}\le 1.$$
Since for all $w_1,w_2\in B_i$, $\|w_1-w_2\|\asymp\gamma^{-\eta(i)}$, for $c=f_0(X_i)$ we have
\begin{equation}\label{eq:beta holder}
|f_0(w)-c|\leq \|w-X_i\|_2^\beta \lesssim \gamma^{-\eta(i)\beta}, \;\;\forall w\in B_i.
\end{equation}
Let us introduce the notation $f_c(x) = f_0(x) -c$. Then, $|f_c(x)|\lesssim \gamma^{-\eta(i)\beta}$, for all $x\in B_i$.
Note that assertion (\ref{eq:H tau trans}) holds with $f_c$ in place of $f_0$. We next proceed by splitting the proof based on the value of $\tau.$

If $\tau \gamma^{-\eta(i)}\le 1$, let us consider the function
$$\tilde{f}(x) = \sum_{w\in  B_i} \tilde{f}_w(x),\qquad\text{with\,\, $ \tilde{f}_w(x)=\frac{f_c(w)}{\varphi_d(0)} \varphi_d\Big(\frac{x-\tau w}{c_d\tau\gamma^{-\eta(i)}/4}\Big)$}. $$
Noting that 
the minimal separation distance among elements in $B_i$ is $c_d\gamma^{-\eta(i)}$ and $\forall w\in B_i, \;\tilde{f}_w$ is supported on a cube with side length $c_d \tau \gamma^{-\eta(i)}/2$, the supports of $\tilde{f}_w,w\in B_i$ are disjoint. Since $\tilde{f}_w(\tau w) =f_c(w)$, denoting  $\hat{\varphi}_d(\cdot)$ as the Fourier transform of $\varphi_d(\cdot)$, we have by equation (\ref{eq:const interp}), 
\begin{align}\label{eq:H w01 1}
&\inf_{f\in\mathbb{H}^\tau,\bm{f}( B_i)= \bm{f}_c( B_i)} \|f\|_{\mathbb{H}^\tau}^2
\leq \|\tilde{f}\|_{\mathbb{H}}^2 \nonumber\\
&\qquad\lesssim   \frac{\sum_{w\in B_i}[f_c(w)]^2}{[\varphi_d(0)]^2}  \int \big|(\tau\gamma^{-\eta(i)})^d \hat{\varphi}_d(\lambda c_d\tau\gamma^{-\eta(i)})\big|^2 (1+\|\lambda\|_2^2)^{\alpha+d/2} d\lambda \nonumber\\
&\qquad\lesssim  \gamma^{-2\eta(i)\beta}\tau^d\gamma^{-\eta(i)d}\int \big|\hat{\varphi}_d(\lambda c_d\tau\gamma^{-\eta(i)})\big|^2 (1+\|\lambda\|_2^2)^{\alpha+d/2} d(\tau\gamma^{-\eta(i)}\lambda) \nonumber\\
&\qquad\lesssim  \gamma^{-2\eta(i)\beta}\tau^d\gamma^{-\eta(i)d} \sup_\lambda \frac{(1+\|\lambda\|_2^2)^{\alpha+d/2}}{(1+\|\lambda c_d\tau\gamma^{-\eta(i)}\|_2^2)^{\alpha+d/2}} \|\varphi_d\|^2_{W^{\alpha+d/2}} \nonumber\\
&\qquad\lesssim  \gamma^{-2\eta(i)\beta}\tau^d\gamma^{-\eta(i)d} (\tau \gamma^{-\eta(i)})^{-(2\alpha+d)} \nonumber \\
&\qquad=  \tau^{-2\alpha}\gamma^{2\eta(i)(\alpha-\beta)}.
\end{align}

If $\tau \gamma^{-\eta(i)}\ge 1$,
let us take
$$\tilde{f}(x) = \sum_{w\in B_i}\frac{f_c(w)}{\varphi_d(0)}\varphi_d\Big(\frac{4\sqrt{d}}{c_d}(x-w)\Big).$$
By the same arguments as in Step 1,
\begin{align}\label{eq:H w01 2}
\inf_{f\in\mathbb{H}^\tau,\bm{f}( B_i)=\bm{f}_c( B_i)} \|f\|_{\mathbb{H}^\tau}^2 
&\lesssim \sum_{w\in B_i} [f_c(w)]^2 \|\varphi_d\|^2_{W^{\alpha+d/2}}\nonumber\\
& \asymp \sum_{w\in B_i} [f_c(w)]^2
\lesssim \gamma^{-2\eta(i)\beta},
\end{align}
where the last step follows from (\ref{eq:beta holder}). Plugging in the bounds \eqref{eq:H w01 1} and \eqref{eq:H w01 2} into  Lemma \ref{lem:GP interp}, yields
\begin{align}\label{eq:H w01}
\left|f_c(0)-\mathbb{E}[Z^{\tau}_{X_{i}}|Z^{\tau}_{\mr{pa}(X_i)}=\bm{f}_c(\mr{pa}(X_i))] \right|
 \lesssim & \sqrt{ \tau^{2\alpha} \gamma^{-2\eta(i)\alpha} \cdot [\tau^{-2\alpha}\gamma^{2\eta(i)(\alpha-\beta)} + \gamma^{-2\eta(i)\beta}] } \nonumber\\
\lesssim & \gamma^{-\eta(i)\beta} + \tau^\alpha \gamma^{-\eta(i)(\alpha+\beta)}.
\end{align}
Then in view of (\ref{eq:H w01}) and (\ref{eq:const interp}), the linearity of the expectation on the condition, and  $\gamma^{\eta(i)\beta}\ge 1$,
\begin{align}
\left|f_0(X_{i})-\mathbb{E}[Z^{\tau}_{X_{i}}|Z^{\tau}_{\mr{pa}(X_i)}=\bm{f}_0(\mr{pa}(X_i))] \right|
& \lesssim  \tau^\alpha \gamma^{-\eta(i)\alpha}+  \gamma^{-\eta(i)\beta} + \tau^\alpha \gamma^{-\eta(i)(\alpha+\beta)}\nonumber\\
& \lesssim \tau^\alpha \gamma^{-\eta(i)\alpha}+  \gamma^{-\eta(i)\beta},\label{eq:lem13:case2}
\end{align}
proving the lemma for $\beta\in (0,1]$.\\

\textbf{Step 3:} Finally, we extend the proof for $\beta>1$.
We use proof by induction on $\beta$. The case when $\beta\in(0,1]$ is proved in Step 2. Suppose the  lemma holds for $\beta\in(0,r]$ for some $r\in\mathbb{N}_{+}$. Now consider the case when $\beta\in (r,r+1]$. Because $r\ge 1$, $f_0$ is a differentiable function. For arbitrary choice of $\zeta=(\zeta[1],\zeta[2], \cdots , \zeta[d])^T\in\mathcal{C}$,
$$
f_0(x) = f_0(\zeta) + \sum_{l=1}^d \int_{\zeta[l]}^{x[l]} \frac{\partial}{\partial t} f_0(\zeta[1],\ldots, \zeta[l-1], t, x[l+1],\ldots,x[d]) dt.
$$
Define the operator $\mathcal{L}$ such that $\forall f:\X\to \mathbb{R}, \;x\in\X$ , 
$\mathcal{L}(f)(x)=f(x)-\mathbb{E}[Z^{\tau}_{x}|Z^{\tau}_{\mr{pa}(X_i)}=\bm{f}(\mr{pa}(X_i))]$. 
Note that $\mathcal{L}$ is a linear operator and can exchange orders with summation and integration (as long as integration is finite). 
Then, in view of the induction assumption and inequality \eqref{eq:lem13:case2}, for all  $x\in \mathcal{C}$,
\begin{align*}
& \left|f_0(x)-\mathbb{E}[Z^{\tau}_{x}|Z^{\tau}_{\mr{pa}(X_i)}=\bm{f}_0(\mr{pa}(X_i))] \right| \\
\le & |\mathcal{L}(f_0)(\zeta)| + \sum_{l=1}^d \int_{\zeta[l]}^{x[l]}\mathcal{L}\Big( \frac{\partial}{\partial t} f_0(\zeta[1],\ldots, \zeta[l-1], t, x[l+1],\ldots,x[d])\Big) dt\\
\lesssim & \tau^\alpha \gamma^{-\eta(i)\alpha} + d \gamma^{-\eta(i)} \big[ \tau^\alpha \gamma^{-\eta(i)\alpha} + \gamma^{-\eta(i)(\beta-1)} \big] \\
\lesssim & \tau^\alpha \gamma^{-\eta(i)\alpha} +  \gamma^{-\eta(i)\beta}.
\end{align*}
We conclude the proof by taking $x=X_i$.
\end{proof}

\subsection{Proof of Lemma \ref{lem:decentering}}
\begin{proof}[Proof of Lemma \ref{lem:decentering}]
The proof consists of three steps. In the first step, we derive a summation formula regarding the RKHS norms of $l(l\le n)$-dimensional marginals of the process $\hat{Z}^{\tau,s}$. In the second step, we obtain a decomposition for all functions in $ L_1(\X_n)$ using the difference formula in Lemma \ref{lem:Holder appro}. Finally, we combine the results of the first two steps to complete the proof.\\

\textbf{Step 1}
Recall the formula for the RKHS of $\hat{Z}^{\tau,s}$ on $\X_n$, see (\ref{f:RKHS decom}). For all $l\le n$, letting $\X_l =\{X_1, X_2, \cdots, X_l\}$, we prove that
\begin{equation}\label{eq:Dm induction}
\bm{f}_0(\X_l)^T (\hat{K}^{\tau,s}_{\X_l,\X_l})^{-1} \bm{f}_0(\X_l) = \sum_{i=1}^l \frac{\big\{f_0(X_{i})-\mathbb{E}[{Z}^{\tau,s}_{X_{i}}|{Z}^{\tau,s}_{\mr{pa}(X_{i})}=\bm{f}_0(\mr{pa}(X_{i}))]\big\}^2 }{
\var \big\{Z^{\tau,s}_{X_{i}}-\mathbb{E}[{Z}^{\tau,s}_{X_{i}}|{Z}^{\tau,s}_{\mr{pa}(X_{i})}]\big\}
}.
\end{equation}
For $i=1$, note that $X_1$ has no parents, implying (\ref{eq:Dm induction}) directly. Next, we assume that equation (\ref{eq:Dm induction}) holds for $l=1,2, \cdots, l'$ and we show below that it holds also for $l=l'+1$. Define the matrix $U\in\mathbb{R}^{(l'+1) \times (l'+1)}$ as
$$
U = \begin{bmatrix}
I & 0 \\ -(\hat K^{\tau,s}_{\X_{l'},X_{{l'}+1}})^T (\hat{K}^{\tau,s}_{\X_{l'},\X_{l'}})^{-1} & 1
\end{bmatrix},
$$
and note that
\begin{align*}
U \hat{K}^{\tau,s}_{\X_{l'+1},\X_{l'+1}} U^T = &
\begin{bmatrix}
\hat{K}^{\tau,s}_{\X_{l'},\X_{l'}} & 0 \\
0 & \hat{K}^{\tau,s}_{X_{l'+1},X_{l'+1}} - (\hatK{\X_{l'},X_{l'+1}})^T (\hatK{\X_{l'},\X_{l'}})^{-1} \hatK{\X_{l'},X_{l'+1}} 
\end{bmatrix}.
\end{align*}
Taking the inverse of the above formula, and noting that $ \hat{K}^{\tau,s}_{X_{l'+1},X_{l'+1}}=s^2$,  the precision matrix on $\X_{l'+1}$ can be computed as
\begin{align}\label{eq:precision decom}
(\hat{K}^{\tau,s}_{\X_{l'+1},\X_{l'+1}})^{-1} 
& = 
U^T \begin{bmatrix}
(\hat{K}^{\tau,s}_{\X_{l'},\X_{l'}})^{-1} & 0 \\
0 & \big[s^2- (\hat K^{\tau,s}_{\X_{l'},X_{l'+1}})^T (\hat{K}^{\tau,s}_{\X_{l'},\X_{l'}})^{-1} \hat{K}^{\tau,s}_{\X_{l'},X_{{l'}+1}})\big]^{-1} 
\end{bmatrix} U.
\end{align}
For the conditional expectation, using the DAG structure of the Vecchia approximation and  equation (\ref{eq:vecchia 2}), 
\begin{align*}
(\hat K^{\tau,s}_{\X_{l'},X_{l'+1}})^T (\hat K^{\tau,s}_{\X_{l'},\X_{l'}})^{-1} \bm{f}_0(\X_{l'}) 
= & \mathbb{E} \big[\hat{Z}^{\tau,s}_{X_{l'+1}}|\hat{Z}^{\tau,s}_{\X_{l'}} = \bm{f}_0(\X_{l'})\big] 
 \nonumber \\
= & \mathbb{E} \big[\hat{Z}^{\tau,s}_{X_{l'+1}}|\hat{Z}^{\tau,s}_{\mr{pa}(X_{l'+1})} = \bm{f}_0(\mr{pa}(X_{l'+1}))\big] 
 \nonumber \\
= & \mathbb{E} \big[{Z}^{\tau,s}_{X_{l'+1}}|{Z}^{\tau,s}_{\mr{pa}(X_{l'+1})} = \bm{f}_0(\mr{pa}(X_{l'+1}))\big].
\end{align*}
The variance can be rewritten along the same lines
\begin{align}\label{eq:K induct 2}
s^2- (\hat K^{\tau,s}_{\X_{l'},X_{l'+1}})^T (\hat{K}^{\tau,s}_{\X_{l'},\X_{l'}})^{-1} \hat{K}^{\tau,s}_{\X_{l'},X_{{l'}+1}} = & 
\var\big[\hat{Z}^{\tau,s}_{X_{l'+1}}- \mathbb{E} \big[\hat{Z}^{\tau,s}_{X_{l'+1}}|\hat{Z}^{\tau,s}_{\X_{l'}} \big] \big] \nonumber\\
=  & \var\big[Z^{\tau,s}_{X_{l'+1}}- \mathbb{E} \big[{Z}^{\tau,s}_{X_{l'+1}}|{Z}^{\tau,s}_{\mr{pa}(X_{l'+1})} \big] \big]. 
\end{align}
By combining the preceding displays and applying the induction assumption,
\begin{align*}
 & \bm{f}_0(\X_{l'+1})^T (\hat{K}^{\tau,s}_{\X_{l'+1},\X_{l'+1}})^{-1} \bm{f}_0(\X_{l'+1}) \\
= &  \Big[ U \bm{f}_0(\X_{l'+1})\Big]^T \begin{bmatrix}
(\hat{K}^{\tau,s}_{\X_{l'},\X_{l'}})^{-1} & 0 \\
0 & \big[s^2- (\hat K^{\tau,s}_{\X_{l'},X_{l'+1}})^T (\hat{K}^{\tau,s}_{\X_{l'},\X_{l'}})^{-1} \hat{K}^{\tau,s}_{\X_{l'},X_{{l'}+1}})\big]^{-1} 
\end{bmatrix} U \bm{f}_0(\X_{l'+1}) \\
= & \frac{\big\{f_0(X_{{l'}+1})-\mathbb{E}[{Z}^{\tau,s}_{X_{{l'}+1}}|{Z}^{\tau,s}_{\mr{pa}(X_{{l'}+1})}=\bm{f}_0(\mr{pa}(X_{{l'}+1}))]\big\}^2}{
\var \big\{Z^{\tau,s}_{X_{{l'}+1}}-\mathbb{E}[{Z}^{\tau,s}_{X_{{l'}+1}}|{Z}^{\tau,s}_{\mr{pa}(X_{{l'}+1})}]\big\}
}  
 + \bm{f}_0(\X_{l'})^T (K^{\tau,s}_{\X_{l'},\X_{l'}})^{-1} \bm{f}_0(\X_{l'}) \\
=  & \sum_{i=1}^{{l'}+1} \frac{\big\{f_0(X_{i})-\mathbb{E}[{Z}^{\tau,s}_{X_{i}}|{Z}^{\tau,s}_{\mr{pa}(X_{i})}=\bm{f}_0(\mr{pa}(X_{i}))]\big\}^2 }{
\var \big\{Z^{\tau,s}_{X_{i}}-\mathbb{E}[{Z}^{\tau,s}_{X_{i}}|{Z}^{\tau,s}_{\mr{pa}(X_{i})}]\big\}
}.
\end{align*}

\textbf{Step 2}
For all $0\le j \le \eta(n)$, we define the operator $\Theta_j:  L_1(\X_n) \to  L_1(\cup_{j'=0}^{j}\mathcal{N}_{j'})$, such that for all $f\in L_1(\X_n)$,
\begin{equation}\label{def:Theta_j}
\Theta_j(f)(X_i) = \Bigg\{ \begin{aligned}
& f(X_{i})-\mathbb{E}[{Z}^{\tau,s}_{X_{i}}|{Z}^{\tau,s}_{\mr{pa}(X_{i})}=\bm{f}(\mr{pa}(X_{i}))], && X_i\in \mathcal{N}_j, \\
& 0, && X_i\in \cup_{j'=0}^{j-1}\mathcal{N}_{j'}.
\end{aligned}
\end{equation}
The operator $\Theta_j$ ``extracts'' the differences between the function values on $\N_j$ and the conditional expectations on $\N_j$ based on the DAG structures. Our objective is to prove the following equation
\begin{equation}\label{eq:Theta decom}
f = \Theta_{\eta(n)}(f) + \sum_{j=0}^{\eta(n)-1} G_{\eta(n)-1}  \cdots G_{j+1} G_{j} \Theta_j(f).
\end{equation}
We prove by induction that for all $0\le j\le \eta(n)$,
\begin{equation}\label{eq:Theta decom induc}
f|_{\cup_{j'=0}^{j} \N_{j'}} = \Theta_{j}(f) + \sum_{j'=0}^{j-1} G_{j-1}  \cdots G_{j'+1} G_{j'} \Theta_{j'}(f),
\end{equation}
where $f|_{\cup_{j'=0}^{j} \N_{j'}}$ denotes the restriction of the function $f$ to the set $\cup_{j'=0}^{j} \N_{j'}$. 
For $j=0$, $\cup_{j'=0}^{j} \N_{j'}=\{X_1\}$ and
$$f(X_1) = \Theta_0(f)(X_1) = f(X_1) - 0 =f(X_1).$$
Suppose that (\ref{eq:Theta decom induc}) holds for $j=\tilde{j}$ and consider the case that $j=\tilde{j}+1$. Note that for all $X_i\in \cup_{j'=0}^{\tilde{j}+1} \N_{j'}$,
\begin{align*}
f(X_i) =  \Theta_{\tilde{j}+1}(f)(X_i) + G_{\tilde{j}}(f|_{\cup_{j'=0}^{\tilde{j}} \N_{j'}})(X_i).
\end{align*}
Then, in view of (\ref{eq:Theta decom induc}) for $j=\tilde{j}$,
\begin{align*}
f|_{\cup_{j'=0}^{\tilde{j}+1} \N_{j'}}
= & \Theta_{\tilde{j}+1}(f) + G_{\tilde{j}} \sum_{j'=0}^{\tilde{j}} G_{\tilde{j}-1}  \cdots G_{j'+1} G_{j'} \Theta_{j'}(f) =
\sum_{j'=0}^{\tilde{j}+1} G_{\tilde{j}}  \cdots G_{j'+1} \Theta_{j'}(f),
\end{align*}
concluding the proof of (\ref{eq:Theta decom induc}) and hence also of (\ref{eq:Theta decom}).\\

\textbf{Step 3}
In view of  Theorem \ref{thm:var} and Lemma  \ref{lem:Holder appro},
$$\var \big\{Z^{\tau,s}_{X_{i}}-\mathbb{E}[{Z}^{\tau,s}_{X_{i}}|{Z}^{\tau,s}_{\mr{pa}(X_{i})}]\big\} \asymp \tau^{2\alpha} s^2 \gamma^{-2\alpha\eta(i)},$$
$$\big| f_0(X_{i})-\mathbb{E}[{Z}^{\tau,s}_{X_{i}}|{Z}^{\tau,s}_{\mr{pa}(X_{i})}=\bm{f}_0(\mr{pa}(X_{i}))] \big| \lesssim \tau^\alpha\gamma^{-\alpha\eta(i)}+ \gamma^{-\beta\eta(i)}.$$
The above two displays together with  equation (\ref{eq:Dm induction}) and  $l \le \gamma^{\eta(l)d}$ imply,
\begin{align}
\bm{f}_0(\X_l)^T (\hat{K}^{\tau,s}_{\X_l,\X_l})^{-1} \bm{f}_0(\X_l) &  \lesssim \sum_{i=1}^l \frac{ \tau^{2\alpha}\gamma^{-2\alpha\eta(i)}+\gamma^{-2\beta\eta(i)}}{\tau^{2\alpha}s^2\gamma^{-2\alpha\eta(i)}} \nonumber\\
&\lesssim   s^{-2}\gamma^{\eta(l)d}+ \tau^{-2\alpha} s^{-2} \gamma^{(2\alpha-2\beta+d)\eta(l)}.\label{UB:RKHS_norm}
\end{align}
Let $\X_l = \cup_{j=0}^{\eta(l)}\N_j$, i.e., $\X_n$ contain all elements up to layer $\eta(l)$.
Let us take $f_1\in\mathbb{H}_n^{\tau,s}$ such that $f_1(x)=f_0(x),\forall x\in \X_l$ and $\Theta_j(f_1)=0, \forall \eta(l)+1\le j \le \eta(n)$.
In view of equation (\ref{eq:Theta decom}),
\begin{align*}
& \sup_{1\le i\le n} |f_1(X_i)-f_0(X_i)| \\
= & \Big\|\sum_{j=0}^{\eta(l)} G_{\eta(n)-1}  \cdots G_{j+1} G_{j} \Theta_j(f_1)  - \Theta_{\eta(n)}(f_0) - \sum_{j=0}^{\eta(n)-1} G_{\eta(n)-1}  \cdots G_{j+1} G_{j} \Theta_j(f_0)\Big\|_\infty\\
= & \bigg\| - \Theta_{\eta(n)}(f_0) - \sum_{j=\eta(l)+1}^{\eta(n)-1} G_{\eta(n)-1} \cdots G_{j+2} G_{j+1}G_{j} \Theta_j(f_0) \bigg\|_\infty  \\
\le & \|\Theta_{\eta(n)}(f_0)\|_\infty + \sum_{j=\eta(l)}^{\eta(n)-1} \|G_{\eta(n)-1} \cdots G_{j+2} G_{j+1}G_{j}\|\cdot  \|\Theta_j(f_0)\|_\infty\\
\lesssim & \tilde\vartheta_n \sum_{j=\eta(l)}^{\eta(n)-1} (\tau^\alpha \gamma^{-j\alpha} + \gamma^{-j\beta}) \lesssim  \tilde\vartheta_n \big( \tau^\alpha \gamma^{-\eta(l)\alpha} + \gamma^{-\eta(l)\beta}\big),
\end{align*} 
where the last inequality utilizes Lemmas \ref{lem:GP recursive} and \ref{lem:Holder appro}.
Let us denote by $c$ the constant multiplier in the upper bound $\lesssim$. Then,  by choosing $\eta(l)$ to be the smallest integer satisfying $c(\tau^\alpha\gamma^{-\eta(l)\alpha}+\gamma^{-\eta(l)\beta})\le \epsilon$, we get 
$$ \|f_1-f_0\|_{\infty,n}\le \tilde\vartheta_n \epsilon.$$
If $\tau^\alpha\gamma^{-\eta(l)\alpha} \ge \gamma^{-\eta(l)\beta}$, then $\epsilon\asymp \tau^\alpha\gamma^{-\eta(l)\alpha}$, which together with \eqref{UB:RKHS_norm} implies that
$$
\inf_{f\in\mathbb{H}_n^{\tau,s}:\|f-f_0\|_{\infty,n}\le \tilde\vartheta_n\epsilon} \|f\|_{\mathbb{H}_n^{\tau,s}}^2 \le  \|f_1\|_{\mathbb{H}_n^{\tau,s}}^2\lesssim s^{-2}\gamma^{\eta(l)d} \asymp s^{-2}\tau^d \epsilon^{-d/\alpha}.
$$  
Otherwise, for $\tau^\alpha\gamma^{-\eta(l)\alpha} < \gamma^{-\eta(l)\beta}$,  $\epsilon\asymp  \gamma^{-\eta(l)\beta}$, hence
$$
\inf_{f\in\mathbb{H}_n^{\tau,s}:\|f-f_0\|_{\infty,n}\le \tilde\vartheta_n \epsilon} \|f\|_{\mathbb{H}_n^{\tau,s}}^2 \le  \|f_1\|_{\mathbb{H}_n^{\tau,s}}^2
\lesssim  s^{-2} \epsilon^{-d/\beta}+ \tau^{-2\alpha} s^{-2} \epsilon^{-\frac{2(\alpha-\beta) +d}{\beta}} \lesssim \tau^{-2\alpha} s^{-2} \epsilon^{-\frac{2(\alpha-\beta) +d}{\beta}}.
$$  
Combining these two cases results in the stated upper bound.

\end{proof}


\section{Proofs of the theorems and corollaries in Section \ref{sec:BNP}}\label{sec:proof BNP}

\subsection{Classical results for posterior contractions}
Deriving contraction rates for the posterior distribution is a well-established topic. Here, we recall two fundamental results from the literature. The first provides sufficient conditions for Bayesian posterior contraction rates, with the corresponding notation adapted for the present paper.

\begin{lemma}
\label{lem:3cond}
For arbitrary fixed $\X_n$, for $\|\cdot\|_{2,n}$ the empirical $ L_2$ norm for functions defined on $\X_n$, if there exists a sequence of sets $F_n \subset  L_1(\X_n)$, such that 
\begin{equation}\label{eq:c rmass}
\ln \pr (\hat{Z}_{\X_n} \not\in F_n) \lesssim -n\epsilon_n^2,
\end{equation}
\begin{equation}\label{eq:c close}
-\ln \pr (\|f_0-\hat{Z}_{\X_n}\|_{2,n}<2\epsilon_n) \lesssim n\epsilon_n^2,
\end{equation}
\begin{equation}\label{eq:c entropy}
\ln N(3\epsilon_n, F_n, \|\cdot\|_{2,n}) \lesssim n\epsilon_n^2,
\end{equation}
then for all $M_n\to \infty$,
$$\Pi(\|f-f_0\|_{2,n} \ge M_n\epsilon_n|\D_n)\to 0 $$ in $P_{f_0}^n$ probability.
\end{lemma}
Lemma \ref{lem:3cond} is a simplified version of Theorem 8.19 in \cite{ghosal2017fundamentals}, with their equations (8.22) and (8.23) replacing (8.20) and (8.21), respectively. Moreover, under the nonparametric regression settings \eqref{eq:regression}, the KL divergence of functions on $\X_n$ is equivalent to empirical $ L_2$ norm on $\X_n$.

The next result is exactly Theorem 2.1 of \cite{van2008rates}, which links the three conditions of general Bayesian posterior contraction with the concentration function of Gaussian processes.

\begin{lemma}
\label{lem:3cond GP}
If $\phi_{f_0,n}^{\tau,s}(\epsilon_n) \le n\epsilon_n^2$ for some sequence $\epsilon_n\to 0$, then there exists a sequence of subsets $F_n\subset\mathbb{R}^n$, such that
 equations \eqref{eq:c rmass} to \eqref{eq:c entropy} hold.
\end{lemma}

\subsection{Proof of Lemma \ref{lem:GP contraction}}
The statement is obtained by straightforward combination of Lemma \ref{lem:3cond} and \ref{lem:3cond GP}.

\subsection{Proof of Theorem \ref{thm:contraction}}
\begin{proof}
In view of Lemma \ref{lem:GP contraction}, it is sufficient to find an $\epsilon_n$ satisfying the inequality
$$\phi^{\tau,s}_{f_0,n}(\epsilon_n) \lesssim n\epsilon_n^2.$$
Since $\phi^{\tau,s}_{f_0,n}(\epsilon_n)$ is upper bounded in Lemma \ref{thm:small:decent}, solving the inequalities for each term of this upper bound yields
$$
\left\{  
\begin{aligned}
& s^{-2}\tau^d (\epsilon_n/\tilde{\vartheta}_n)^{-\frac{d}{\alpha}} \lesssim n\epsilon_n^2 \\
& s^{-2}\tau^{-2\alpha} (\epsilon_n/\tilde{\vartheta}_n)^{-\frac{2(\alpha-\beta)+d}{\beta}} \lesssim n\epsilon_n^2 \\
& s^{d/\alpha} \tau^d (\epsilon_n/\tilde{\vartheta}_n)^{-\frac{d}{\alpha}} \lesssim n\epsilon_n^2 
\end{aligned}
\right. \Rightarrow
\left\{ 
\begin{aligned}
& \epsilon_n \gtrsim (ns^2\tau^d)^{-\frac{\alpha}{2\alpha+d}} \tilde{\vartheta}_n^{\frac{d}{2\alpha+d}}\\
& \epsilon_n \gtrsim n^{-\frac{\beta}{2\alpha+d}} (s\tau^\alpha)^{-\frac{\beta}{2\alpha+d}} \tilde{\vartheta}_n^{\frac{2(\alpha-\beta)+d}{2\alpha+d}} \\
& \epsilon_n \gtrsim n^{-\frac{\alpha}{2\alpha+d}}(s\tau^\alpha)^{\frac{d}{2\alpha+d}} \tilde{\vartheta}_n^{\frac{d}{2\alpha+d}}
\end{aligned}
\right.
$$
Taking the maximum of the lower bounds on the right-hand side provides the stated contraction rate.
\end{proof}

\subsection{Proof of Corollary \ref{cor:rates}}
\begin{proof}
The case that both $\tau$ and $s$ are constants can be directly obtained from Theorem \ref{thm:contraction}. For $\tau^\alpha s \asymp n^{\frac{\alpha-\beta}{2\beta+d}}$, substituting this into equation (\ref{eq:epsilon choice}) yields
$$n^{-\frac{\beta}{2\alpha+d}} (\tau^\alpha s)^{-\frac{2\beta}{2\alpha+d}} = n^{-\frac{\alpha}{2\alpha+d}} (\tau^\alpha s)^{\frac{d}{2\alpha+d}} = n^{-\frac{\beta}{2\beta+d}}.$$
Because $s\gtrsim 1$, we have
$$(ns^2)^{-\frac{\beta}{2\beta+d}} \lesssim n^{-\frac{\beta}{2\beta+d}} \text{ and } (ns^2\tau^{-d})^{-\frac{\alpha}{2\alpha+d}}  \lesssim n^{-\frac{\alpha}{2\alpha+d}} (\tau^\alpha s)^{\frac{d}{2\alpha+d}}.$$
Therefore, the posterior contraction rate is $n^{-\frac{\beta}{2\beta+d}}$.
\end{proof}

\subsection{Proof of Theorem \ref{thm:adaptation}}
We first prove a (nearly) nested property of the RKHS balls corresponding to the Vecchia GPs  $\hat{Z}^{\tau,s}$, with respect to the scaling parameters $s$ and $\tau$. For the space scaling parameter $s$, the inclusions are strictly monotone and straightforward. However, for the time rescaling $\tau$, the analysis is substantially more challenging. To begin with, let us denote the unit balls in the RKHS and the Banach space $(\mathcal{X}_n,\|\cdot\|_{n,\infty})$  by
$$\mathbb{H}_n^{\tau,s}(1) = \{f\in\mathbb{H}_n^{\tau,s}: \|f\|_{\mathbb{H}_n^{\tau,s}}\le 1\},\quad \mathbb{B}_n(1)=\{f: \mathcal{X}_n\rightarrow \mathbb{R},\, \|f\|_{n,\infty}\leq 1\},$$
respectively.

\begin{lemma}\label{lem:RKHS continuity}
For arbitrary $c>0$, there exists a universal constant $C>0$, such that for all $0<\tau, \Delta \tau, s, \Delta s$ with $\Delta \tau\le 1$, $c\le (\tau-\Delta  \tau)/\tau $, 
$$\mathbb{H}^{\tau,s}_n(1) \subset  \mathbb{H}_n^{\tau-\Delta\tau,s-\Delta s}(C) + \mathbb{B}_n\big( C [\Delta \tau + \Delta \tau^{2(\alpha-\underline{\alpha})} ] \tau^\alpha s \ln n  \big).$$
\end{lemma}
\begin{proof}
The proof consists of two parts. First, we compare $\mathbb{H}_n^{\tau-\Delta\tau,s}(1)$ and $\mathbb{H}_n^{\tau,s}$ and then $\mathbb{H}_n^{\tau-\Delta\tau,s-\Delta s}$ and $\mathbb{H}_n^{\tau-\Delta\tau,s}$.\\

\textbf{Step 1:}
For $f\in\mathbb{H}^{\tau,s}_n(1)$, in view of (\ref{eq:f RKHS def}), (\ref{eq:Dm induction}) and \eqref{def:Theta_j}, 
\begin{align}\label{eq:fn RKHS}
\|f\|_{\mathbb{H}^{\tau,s}_n}^2 &= \sum_{i=1}^n \frac{\big\{f(X_{i})-\mathbb{E}[{Z}^{\tau,s}_{X_{i}}|{Z}^{\tau,s}_{\mr{pa}(X_{i})}=\bm{f}(\mr{pa}(X_{i}))]\big\}^2 }{
\var \big\{Z^{\tau,s}_{X_{i}}-\mathbb{E}[{Z}^{\tau,s}_{X_{i}}|{Z}^{\tau,s}_{\mr{pa}(X_{i})}]\big\}
} \nonumber\\
&= \sum_{i=1}^n \frac{ [\Theta_{\eta(i)}(f)(X_i)]^2 }{
\var \big\{Z^{\tau,s}_{X_{i}}-\mathbb{E}[{Z}^{\tau,s}_{X_{i}}|{Z}^{\tau,s}_{\mr{pa}(X_{i})}]\big\}
}.
\end{align}
Furthermore, recall the Gaussian interpolation operator $G_j: \cup_{j'=0}^{j-1}\mathcal{N}_{j'} \to \mathcal{N}_j$ in \eqref{def:op:G}. For clarity we denote by $G^\tau$ the Gaussian interpolation operator used to construct $\hat{Z}^{\tau,s}$. In view of equation (\ref{eq:Theta decom}),
$$f = \sum_{j=0}^{\eta(n)-1} G^\tau_{\eta(n)-1} \cdots G^\tau_{j+1} G^\tau_j \Theta_j(f).$$
For $\tilde{\tau}=\tau-\Delta \tau$, we define the function $\tilde{f}$ on $\X$ as
$$\tilde f = \sum_{j=0}^{\eta(n)-1} G^{\tilde\tau}_{\eta(n)-1} \cdots G^{\tilde\tau}_{j+1} G^{\tilde\tau}_j \Theta_j(f).$$

By Lemma \ref{lem:GP recursive} and the assumption $\tilde\vartheta_n=O(1)$, there exists a constant $C$, such that for all $j_1, j_2\in\mathbb{N}$, 
\begin{equation}\label{eq:GP recurisve}
\| G^{\tau}_{j_2}\cdots G^{\tau}_{j_1}\| \lesssim 1, \;\; \| G^{\tilde\tau}_{j_2}\cdots G^{\tilde\tau}_{j_1}\| \lesssim 1.
\end{equation}
Since $0<\Delta\tau<1$, by Lemma \ref{lem:GP flat}, for all $j\in\mathbb{N}$, 
\begin{equation}\label{eq:GP interp diff}
\|G_j^\tau - G_j^{\tilde{\tau}}\|_\infty \lesssim \Delta \tau + \Delta \tau^{2(\alpha-\underline{\alpha})}.
\end{equation}
Combing the previous two displays with the formulas for  $f(\X_n)$ and $\tilde{f}(\X_n)$ yields
\begin{align*}
\|f(\X_n) - \tilde{f}(\X_n)\|_\infty 
& = \Big\|\sum_{j=0}^{\eta(n)-1} (G^\tau_{\eta(n)-1} \cdots G^\tau_{j+1} G^\tau_j-G^{\tilde\tau}_{\eta(n)-1} \cdots G^{\tilde\tau}_{j+1} G^{\tilde\tau}_j)\Theta_j(f) \Big\|_\infty \\
& = \Big\|\sum_{j=0}^{\eta(n)-1} \sum_{l=j}^{\eta(n)-1} G_{\eta(n)-1}^{\tilde{\tau}}\cdots G_{l+1}^{\tilde{\tau}} (G_l^\tau -G_l^{\tilde{\tau}}) G_{l-1}^\tau \cdots G_j^\tau \Theta_j(f) \Big\|_\infty \\
& \lesssim \eta(n) \big[\Delta \tau + \Delta \tau^{2(\alpha-\underline{\alpha})} \big] \sum_{j=0}^{\eta(n)-1} \|\Theta_j(f)\|_\infty.
\end{align*}

Furthermore, in view of  Cauchy-Schwarz inequality, the definition of $\Theta_j(f)$, assertion (\ref{eq:fn RKHS}) and Theorem \ref{thm:var}, 
\begin{align*}
\Big(\sum_{j=0}^{\eta(n)-1} \|\Theta_j(f)\|_\infty\Big)^2
& \le \bigg[\sum_{j=0}^{\eta(n)-1} \sqrt{\sum_{X_i\in \mathcal{N}_j} [\Theta_{\eta(i)}(f)(X_i)]^2} \bigg]^2\\
& \le \bigg[\sum_{j=0}^{\eta(n)-1}
s^2\tau^{-2\alpha}\gamma^{2\alpha j} \sum_{X_i\in\N_j} [\Theta_{\eta(i)}(f)(X_i)]^2\bigg]
\bigg[\sum_{j=0}^{\eta(n)-1} s^2\tau^{2\alpha}\gamma^{-2\alpha j}\bigg] \\
& \lesssim \|f\|_{\mathbb{H}^{\tau,s}_n}^2 \bigg[\sum_{j=0}^{\eta(n)-1}s^2 \tau^{2\alpha}\gamma^{-2\alpha j}\bigg]\lesssim  s^2\tau^{2\alpha}.
\end{align*}
Therefore, by combining the previous two displays, we obtain
\begin{align}\label{eq:f diff inf}
\|f(\X_n) - \tilde{f}(\X_n)\|_\infty &\lesssim \eta(n) \big[\Delta \tau + (\Delta \tau)^{2(\alpha-\underline{\alpha})} \big] s\tau^{\alpha} \nonumber\\
& \lesssim  \big[\Delta \tau + (\Delta \tau)^{2(\alpha-\underline{\alpha})} \big] s\tau^{\alpha}  \ln n. 
\end{align}

We proceed to bound the $\mathbb{H}^{n,\tilde{\tau}}$ norm of the function $\tilde{f}$. By the definition of $\tilde{f}$,  equation (\ref{eq:fn RKHS}) and Theorem \ref{thm:var}
\begin{align}\label{eq:H norm alt}
\|\tilde{f}\|_{\mathbb{H}^{\tilde\tau,s}_n}
& = \sum_{i=1}^n \frac{[\Theta_{\eta(i)}(f)(X_i)]^2 }{
\var \big\{Z^{\tilde\tau,s}_{X_{i}}-\mathbb{E}[{Z}^{\tilde\tau,s}_{X_{i}}|{Z}^{\tilde\tau,s}_{\mr{pa}(X_{i})}]\big\} } \nonumber\\
& = \sum_{i=1}^n \frac{\big\{f(X_{i})-\mathbb{E}[{Z}^{\tau,s}_{X_{i}}|{Z}^{\tau,s}_{\mr{pa}(X_{i})}=\bm{f}(\mr{pa}(X_{i}))]\big\}^2 }{
\var \big\{Z^{\tilde\tau,s}_{X_{i}}-\mathbb{E}[{Z}^{\tilde\tau,s}_{X_{i}}|{Z}^{\tilde\tau,s}_{\mr{pa}(X_{i})}]\big\} } \nonumber\\
& \le \|f\|_{\mathbb{H}^{\tau,s}_n} \sup_{X_i\in \X} \frac{ \var \big\{Z^{\tau,s}_{X_{i}}-\mathbb{E}[{Z}^{\tau,s}_{X_{i}}|{Z}^{\tau,s}_{\mr{pa}(X_{i})}]\big\} }{ \var \big\{Z^{\tilde\tau,s}_{X_{i}}-\mathbb{E}[{Z}^{\tilde\tau,s}_{X_{i}}|{Z}^{\tilde\tau,s}_{\mr{pa}(X_{i})}]\big\} } \nonumber\\
& \lesssim  \|f\|_{\mathbb{H}^{\tau,s}_n} \sup_{1\le j \le \eta(n)} \frac{s^2(\tau \gamma^{-j})^{2\alpha}}{s^2(\tilde\tau \gamma^{-j})^{2\alpha}} \lesssim  \|f\|_{\mathbb{H}^{\tau,s}_n}.
\end{align}
Note that in view of (\ref{eq:f diff inf}) and (\ref{eq:H norm alt}), for all  $f\in\mathbb{H}_n^{\tau,s}(1)$, we have $\tilde{f}\in \mathbb{H}_n^{\tau-\Delta\tau,s}(C)$ and  $f-\tilde{f}\in\mathbb{B}_n\big(  [\Delta \tau + \Delta \tau^{2(\alpha-\underline{\alpha})} ] \tau^\alpha s \ln n  \big)$, for some constant $C>0$.
Therefore, by writing $f=\tilde{f}+f-\tilde{f}$,
\begin{equation}\label{eq:RKHS embed 1}
\mathbb{H}_n^{\tau,s}(1) \subset  \mathbb{H}_n^{\tau-\Delta\tau,s}(C) + \mathbb{B}_n\big(  C[\Delta \tau + \Delta \tau^{2(\alpha-\underline{\alpha})} ] \tau^{\alpha} s \ln n  \big). 
\end{equation}
\vspace{0.1cm}

\textbf{Step 2:}
By the definition of the space scaling parameter $s$, (\ref{eq:fn RKHS}) can be reformulated as
$$
\|f\|_{\mathbb{H}^{\tau,s}_n}^2 = \sum_{i=1}^n \frac{\big\{f(X_{i})-\mathbb{E}[{Z}^{\tau,1}_{X_{i}}|{Z}^{\tau,1}_{\mr{pa}(X_{i})}=\bm{f}(\mr{pa}(X_{i}))]\big\}^2 }{s^2
\var \big\{Z^{\tau,1}_{X_{i}}-\mathbb{E}[{Z}^{\tau,1}_{X_{i}}|{Z}^{\tau,1}_{\mr{pa}(X_{i})}]\big\}
}.
$$
Therefore, there exist constants $C,C'>0$ such that
\begin{equation}\label{eq:RKHS embed 2}
\mathbb{H}_n^{\tau-\Delta\tau,s}(1) \subset \mathbb{H}_n^{\tau-\Delta \tau,s-\Delta s}\Big(C\frac{1/s^2}{1/(s-\Delta s)^2}\Big) = \mathbb{H}_n^{\tau-\Delta \tau,s-\Delta s}(C').
\end{equation}
Finally, we finish the proof by combining (\ref{eq:RKHS embed 1}) and (\ref{eq:RKHS embed 2}).
\end{proof}

Furthermore, recall the fundamental Theorem 2.1 of \cite{van2008rates}, which links the concentration function with the three conditions for Bayesian posterior contraction rates.

We turn our attention to the proof of Theorem \ref{thm:adaptation}.

\begin{proof}[Proof of Theorem \ref{thm:adaptation}]
Let $\epsilon_n \asymp  n^{-\frac{\beta}{2\beta+d}}$ and define the set 
$$T = \{(\tau,s): \tau^\alpha s \le n^{\frac{\alpha-\beta}{2\beta+d}}, s\ge n^{-1}, \tau\ge n^{-\frac{\beta}{2\beta+d}}\ln^2 n\}.$$
Furthermore, let us consider the sieve
\begin{equation}\label{eq:sieve}
F_n = \cup_{(\tau,s)\in T}  \mathbb{H}_n^{\tau,s}(a_n) +  \mathbb{B}_n(\epsilon_n),
\end{equation}
where $a_n$ is chosen as 
$$ a_n = \inf\left\{y\in\mathbb{R}: y\ge-2\Phi^{-1}(\pr(\|\hat{Z}^{\tau,s}_{\X_n}\|_\infty<\epsilon_n)), \;\;\forall \tau^\alpha s \le n^{\frac{\alpha-\beta}{2\beta+d}} \right\},$$
with $\Phi$ denoting the standard normal cdf. Then, in view of $\Phi^{-1}(u)\asymp -\sqrt{\ln (1/u)}$, for $u\in(0,1/2)$, see Lemma K.6 of \cite{ghosal2017fundamentals}, and recalling Theorem \ref{lem:small ball} and  the assumption $\tilde\vartheta_n=O(1)$,
$$a_n^2 \asymp \sup_{ \tau^\alpha s \le n^{\frac{\alpha-\beta}{2\beta+d}} }-\ln\pr(\|\hat{Z}^{\tau,s}_{\X_n}\|_\infty<\epsilon_n) \lesssim \big(n^{\frac{\alpha-\beta}{2\beta+d}}\big)^{d/\alpha} \big(n^{-\frac{\beta}{2\beta+d}}\big)^{-d/\alpha}  \asymp n^{\frac{d}{2\beta+d}}=n\eps_n^2.$$

Then, in view of Lemma \ref{lem:3cond},
 it suffices to prove that conditions (\ref{eq:c rmass})-(\ref{eq:c entropy}) hold for $\hat{Z}^{\pr}_{\X_n}$, the restriction of the process $\hat{Z}^{\pr}$ to $\X_n$. \\

\textbf{Proof of (\ref{eq:c rmass}):}\mbox{}
In view of the definition of $\hat{Z}_{\X_n}^{\mr{Pr}}$ and Borell's inequality, see Theorem 3.1 of \cite{borell1975brunn}, 
\begin{align}
\label{eq:1e}
\pr(\hat{Z}^{\mr{Pr}}_{\X_n} \not\in F_n) 
\le & \mr{Pr}((\tau,s)\not\in T) + \int_{(\tau,s)\in T}\pr(\hat{Z}^{\tau,s}_{\X_n} \not\in  \mathbb{H}_n^{\tau,s}(a_n) + \mathbb{B}_n(\epsilon_n ) ) p(\tau,s) d\tau ds \nonumber\\
\le & \mr{Pr}((\tau,s)\not\in T) + \sup_{(\tau,s)\in T}\pr(\hat{Z}^{\tau,s}_{\X_n} \not\in  \mathbb{H}_n^{\tau,s}(a_n) + \mathbb{B}_n(\epsilon_n ) ) \nonumber\\
\le & \mr{Pr}(\tau^\alpha s> n^{\frac{\alpha-\beta}{2\beta+d}}) + \mr{Pr}(\tau<n^{-\frac{\beta}{2\beta+d}}\ln^2 n) + \mr{Pr}(s<1/n) \nonumber\\
& + \sup_{\tau^\alpha s \le n^{\frac{\alpha-\beta}{2\beta+d}}}\Big[1- \Phi\Big(\Phi^{-1}\big(\pr(\|\hat{Z}^{\tau,s}_{\X_n}\|_\infty<\epsilon_n)\big) + a_n\Big)\Big].
\end{align}
Then, by the definition of $a_n$,
\begin{align*}
\sup_{\tau^\alpha s \le n^{\frac{\alpha-\beta}{2\beta+d}}}&\Big[1- \Phi\Big(\Phi^{-1}\big(\pr(\|\hat{Z}^{\tau,s}_{\X_n}\|_\infty<\epsilon_n)\big) + a_n\Big)\Big]\\
&\le 1- \Phi ( a_n/2) \asymp a_n^{-1}\exp(-a_n^2/8).
\end{align*}
Note that minus the logarithm of the right-hand side is $O(n^{\frac{d}{2\beta+d}})$.

Recalling Condition \ref{cond:hyperprior}, we have
\begin{equation}\label{eq:3e}
\ln \left[\mr{Pr}(\tau^\alpha s> n^{\frac{\alpha-\beta}{2\beta+d}}) + \mr{Pr}(\tau<n^{-\frac{\beta}{2\beta+d}}\ln^2 n) + \mr{Pr}(s<1/n) \right] \lesssim - n^{\frac{d}{2\beta+d}}.
\end{equation}
The combination of the above displays verifies the remaining mass condition
$$\ln(\pr(\hat{Z}^{\pr}_{\X_n}\not\in F_n)) \lesssim -n^{\frac{d}{2\beta+d}}\asymp -n\epsilon_n^2.$$

\textbf{Proof of (\ref{eq:c close}):}\mbox{}
Let us consider the set $\tilde{T} = \big\{\tau^\alpha s \in \big[n^{\frac{\alpha-\beta}{2\beta+d}}, 2n^{\frac{\alpha-\beta}{2\beta+d}}\big]\big\}\cap\{s\ge 1 \}$. Then, in view of \eqref{eq:conc:ball},
\begin{align*}
-\ln \pr (\|f_0-\hat{Z}^{\pr}_{\X_n}\|_{\infty,n}<\epsilon_n)
\le & -\ln \int_{\tilde T} e^{-\phi_{f_0,n}^{\tau,s}(\epsilon_n)} p(\tau,s) d\tau ds\\
\le & -\ln \big[\pr(\tilde T) \big] + \sup_{(\tau,s)\in \tilde{T}} \phi_{f_0,n}^{\tau,s}(\epsilon_n) \\
\lesssim & n^{\frac{d}{2\beta+d}} \asymp n\epsilon_n^2,
\end{align*}
where in the last inequality, $\ln \big[\pr(\tilde T) \big]$ is controlled by Condition \ref{cond:hyperprior} and the concentration function by Corollary \ref{cor:rates}, providing the small ball probability condition.

\textbf{Proof of (\ref{eq:c entropy}):}\mbox{}
Define $\Delta \tau\in(0,1)$ such that
$$\Delta \tau = \frac{1}{2}\min\left\{\frac{\epsilon_n }{a_n n^{\frac{\alpha-\beta}{2\alpha+d}} \ln n},  \Big(\frac{\epsilon_n }{a_n n^{\frac{\alpha-\beta}{2\alpha+d}} \ln n}\Big)^{\frac{1}{2(\alpha-\underline{\alpha})}} \right\}.$$
Then we have $\forall (\tau,s)\in T$, 
$$ a_n \big[ \Delta \tau + \Delta \tau^{2(\alpha-\underline{\alpha})} \big] \tau^\alpha s \ln n \le a_n \big[ \Delta \tau + \Delta \tau^{2(\alpha-\underline{\alpha})} \big] n^{\frac{\alpha-\beta}{2\alpha+d}} \ln n \le \epsilon_n. $$
Recalling the definition of $T$ and $\Delta \tau$, for all $(\tau,s)\in T$, we have 
$$\tau \ge n^{-\frac{\beta}{2\beta+d}} \ln^2 n \ge  \epsilon_n /[a_n n^{\frac{\alpha-\beta}{2\alpha+d}} \ln n]  \ge \Delta \tau,$$ 
$$\tau \le \Big( n^{\frac{\alpha-\beta}{2\beta+d}} / n^{-1} \Big)^{\frac{1}{\alpha}} = n^{\frac{\alpha+\beta+d}{(2\beta+d)\alpha}}.$$ 
Now define a finite set $T_n\subset\mathbb{R}_+\times\mathbb{R}_+$ as
$$T_n = \left\{ \Big(l  \Delta \tau, n^{-1}\Big): l\in\mathbb{N},\; 1\le l \le n^{\frac{\alpha+\beta+d}{(2\beta+d)\alpha}} (\Delta \tau)^{-1} \right\}.$$
Thus,  for all $(\tau,s)\in T$, there exists $(\tau',n^{-1})\in T_n$, such that $0<\tau'<\tau$, $1/2 < \tau'/\tau$ and $\tau - \tau'<\Delta \tau$. Then, by Lemma \ref{lem:RKHS continuity} (with $c=1/2$, $\Delta s=s-n^{-1}$), there exists a constant $C>0$, such that
\begin{align*}
\mathbb{H}_n^{\tau,s}( a_n) &\subset  \mathbb{H}_n^{\tau',n^{-1}}(Ca_n) + \mathbb{B}_n\big(Ca_n\big[\Delta \tau + \Delta \tau^{2(\alpha-\underline{\alpha})} \big]\tau^{\alpha}s \ln n\big)\\
& \subset \mathbb{H}_n^{\tau',n^{-1}}(Ca_n) + \mathbb{B}_n(C\epsilon_n).
\end{align*}
Recalling the definition of $F_n$ in equation (\ref{eq:sieve}), the above display implies that 
$$ F_n \subset \cup_{(\tau,s)\in T_n} \mathbb{H}_n^{\tau,s}(Ca_n) + \mathbb{B}_n(C\eps_n), $$
for some universal constant $C>0$. Therefore
$$\ln N(\epsilon_n, F_n, \|\cdot\|_{\infty,n}) \lesssim \ln|T_n| + \sup_{(\tau,s) \in T_n} \ln N\big(\epsilon_n, \mathbb{H}^{\tau,s}_n(Ca_n) + \mathbb{B}_n(C\eps_n), \|\cdot\|_{\infty,n}\big).$$
Following the proof of Theorem 2.1 of \cite{van2008rates} in the middle of page 21, the above $\epsilon_n$-entropy can be bounded by $(Ca_n)^2$ plus the centered small ball probability. Combining this with the centered small ball probability derived in Theorem \ref{lem:small ball} yields 
\begin{align*}
\sup_{(\tau,s) \in T_n}\ln N(\epsilon_n,  \mathbb{H}^{\tau,s}_n(Ca_n) +  \mathbb{B}_n(C\epsilon_n), \|\cdot\|_{\infty,n})
\lesssim &(a_n)^2 + \sup_{(\tau,s) \in T_n} \phi_{0,n}^{\tau,s}(\epsilon_n)\\
\lesssim & (C a_n)^2 + \sup_{(\tau,s)\in T_n} \tau^d s^{d/\alpha} \epsilon^{-d/\alpha} \\
\lesssim & a_n^2 + (n^{\frac{\alpha-\beta}{(2\beta+d)\alpha}})^d (n^{-\frac{\beta}{2\beta+d}})^{-\frac{d}{\alpha}} \\
\lesssim & n^{\frac{d}{2\beta+d}}\asymp n\epsilon_n^2.
\end{align*}
Therefore,
\begin{align*}
\ln N(\epsilon_n, F_n, \|\cdot\|_{\infty,n}) \lesssim \ln|T_n| + n\epsilon_n^2 \lesssim \ln \left[n^{\frac{\alpha-\beta}{(2\beta+d)\alpha}}(\Delta\tau)^{-1} \right] + n\epsilon_n^2\lesssim n\epsilon_n^2,
\end{align*}
concluding the proof of the entropy condition and hence the proof of the theorem.
\end{proof}

\subsection{Proof of Corollary \ref{cor:tau adapt}}
\begin{proof}
We only have to verify that Condition \ref{cond:hyperprior} holds. Once this is established, Theorem \ref{thm:adaptation} directly implies the result.

Since the prior is supported on the set $\{(\tau,s):s=1,\tau\ge 1\}$, the inequalities $\ln\pr (s<n^{-1})\lesssim -n^{\frac{d}{2\beta+d}}$ and $\ln\pr (\tau<n^{-\frac{\beta}{2\beta+d}}\ln^2n )\lesssim -n^{\frac{d}{2\beta+d}}$ are straightforward. Then, in view of the assumption on the hyper-prior $p(\tau)$, there exist constants $c_1,c_2>0$ such that
\begin{align*}
\pr\Big(\tau^\alpha > n^{\frac{\alpha-\beta}{2\beta+d}} \Big) \le & \int_{n^{\frac{\alpha-\beta}{\alpha(2\beta+d)}}}^\infty c_2 \exp\big( -c_1n^{\frac{d}{2\alpha+d}} \tau^{\frac{2\alpha d}{2\alpha+d}}\big) d\tau.
\end{align*}
Let $u = c_1n^{\frac{d}{2\alpha+d}} \tau^{\frac{2\alpha d}{2\alpha+d}}$ and 
$u_n = c_1n^{\frac{d}{2\alpha+d}} (n^{\frac{\alpha-\beta}{\alpha(2\beta+d)}})^{\frac{2\alpha d}{2\alpha+d}}= n^{\frac{d}{2\beta+d}},$
we have
\begin{align*}
\pr\Big(\tau^\alpha > n^{\frac{\alpha-\beta}{2\beta+d}} \Big) \le & 
\int_{u_n}^\infty c_2 c_1^{-\frac{2\alpha+d}{2\alpha d}} n^{-\frac{1}{2\alpha}} e^{-u} u^{\frac{2\alpha+d}{2\alpha d}-1} du.
\end{align*}
Further, let $v =u- u_n/2$. Noticing $u/v\in (1,2]$ for all $u\ge u_n$, 
\begin{align*}
\pr\Big(\tau^\alpha > n^{\frac{\alpha-\beta}{2\beta+d}} \Big) \le & 
c_2 c_1^{-\frac{2\alpha+d}{2\alpha d}} n^{-\frac{1}{2\alpha}} \max\{1, 2^{\frac{2\alpha+d}{2\alpha d}-1} \} e^{-u_n/2} \int_{u_n/2}^{\infty} e^{-v} v^{\frac{2\alpha+d}{2\alpha d}-1} dv\\
\le & c_2 c_1^{-\frac{2\alpha+d}{2\alpha d}} n^{-\frac{1}{2\alpha}} \max\{1, 2^{\frac{2\alpha+d}{2\alpha d}-1} \} e^{-u_n/2} \;\Gamma\Big(\frac{2\alpha+d}{2\alpha d}\Big).
\end{align*}
Taking the logarithm, we have
$$
\ln \pr\Big(\tau^\alpha > n^{\frac{\alpha-\beta}{2\beta+d}} \Big) \lesssim -u_n/2 - (\ln n)/(2\alpha) \lesssim - n^{\frac{d}{2\beta+d}}.
$$
Finally, for $\tau^\alpha\in\big[ n^{\frac{\alpha-\beta}{2\beta+d}}, 2n^{\frac{\alpha-\beta}{2\beta+d}} \big]$, we have
$\ln p(\tau) \asymp -n^{\frac{d}{2\beta+d}}$ and thus
$$ 
\ln\pr\Big(\tau^\alpha \in [ n^{\frac{\alpha-\beta}{2\beta+d}},2n^{\frac{\alpha-\beta}{2\beta+d}}] \Big) \asymp -n^{\frac{d}{2\beta+d}} + \ln \big( n^{\frac{\alpha-\beta}{2\beta+d}} \big) \asymp -n^{\frac{d}{2\beta+d}}.
$$
\end{proof}

\end{appendix}

\end{document}